\pdfoutput=1
\documentclass{article}
\usepackage[english]{babel}
\usepackage{amsmath,amssymb,enumerate,bbm}
\usepackage{amsthm}
\usepackage{enumitem}

\newcommand{\tmop}[1]{\ensuremath{\operatorname{#1}}}

\newtheorem{theorem}{Theorem}[section]
\newtheorem{lemma}[theorem]{Lemma}
\theoremstyle{definition}
\newtheorem{definition}[theorem]{Definition}

\theoremstyle{remark}
\newtheorem{remark}[theorem]{Remark}

\usepackage{enumitem}

\usepackage{anysize}
\marginsize{2.5cm}{2.5cm}{2.5cm}{3.4cm}

\title{Free modules with isomorphic duals}
\author{Theodoros Ioannis Kyriopoulos }

\begin{document}

\maketitle

\begin{abstract}
Let $M,N$ be free modules over a Noetherian commutative ring  $R$  and let  $F$ be a field such that $\mathrm{card}(F)$ does not exceed the continuum. Then : \\
$\mathbf{(1)}$ The  assertion  that [Any two 
 $F$-vector spaces with isomorphic duals are isomorphic] is equivallent to the  $\mathsf{ICF}$ (Injective continium function) hypothesis and  it is a non-decidable statement  in $\mathsf{ZFC}$. \\
$\mathbf{(2)}$ If the dual of $M$ is a projective $R$-module  and  
$\mathrm{rank}_R(M)$ is infinite then the ring $R$  is Artinian.\\
$\mathbf{(3)}$ If $R$ is Artinian and $\tmop{card}(R)$  does not exceed the continuum  then the the dual of $M$ is  free.\\
$\mathbf{(4)}$ Assume that  $R$ is a non-Artinian ring that is either   Hilbert or countable. Then $:$  \\  $\mathrm{(a)}$ If $M,\, N$ have isomorphic duals then they are themselves isomophic \\ 
 $\mathrm{(b)}$ Any free direct summand of the dual of $M$ is finitely gernerated, if $\mathrm{Card}(R)$ is not $\omega$-measurable.\\
 $\mathrm{(c)}$ If  $R$ is connected and both
 $\mathrm{Card}(R)$, $\mathrm{rank}_R(M)$ are not $\omega$-measurable then [Any direct summand of the dual of $M$ that is not finitely generated is a dual of a 
free $R$-module].\\
$\mathbf{(5)}$ If $R$ is a non-local domain then  $R$ is a half-slender ring.\\
 $\mathbf{(6)}$ If $R$ is  Artinian ring and it's cardinality $\mathrm{card}(R)$  does not exceed the continuum then the assertion that [any two free $R$-modules with isomorphic duals are isomorphic] is  non-decidable in  $\mathsf{ZFC}$.\\
$\mathbf{(7)}$ If $R$ is a domain and $\mathrm{rank}_R(M)$ is infinite then the Goldie dimension of the dual of $M$ is equal to it's cardinality .   \\
$\mathbf{(8)}$ If $R$  is a complex affine algebra whose corresponding affine variety has no isolated points then [any two projective
 $R$-modules with isomorphic duals are themselves isomorphic]. \\
$\mathbf{(9)}$ Let $V$ be an $F$-vector space of infinite dimension $\lambda$. The dimension of it's dual is $2^{\lambda}$.  \\

\end{abstract}

\section{ Basic concepts and notation}
\subsection{Basic Set theory} \label{sets}
\begin{enumerate}[label=(\roman*)]
\item  \underline{All set theory in this article assumes $ \mathsf{ZFC}$}.  Let $X$ be a set. By $|X|=\mathrm{card}(X)$, we denote it's cardinality. Also  $\, \mathcal{P}(X)$ is it's poweset and   $ \mathcal{P}_{\text{fin}}(X)$ 
 $=\{ Y\in  \mathcal{P}(Y)\ |\, Y$ is finite$\}$.  By $ \Phi_X$, we denote the map  $\Phi_X:\mathcal{P}(X)\rightarrow (\mathbb{F}_2)^X $,  defined  by the formula  $ \Phi_X(A)(t)=[\, 1$ if $t\in A$, $0$ if 
$t \notin A\,]$.
It is trivial that $ \Phi_X$ is a bijection and  $[\,\Phi_X(A\oplus B)=\Phi_X(A)+\Phi_X( B)\,]$,                  where $A\oplus B $ is the symmetric difference of sets defined by 
$A\oplus B 
=(A\smallsetminus B)\cup (B\smallsetminus A)$.
For 
$\Delta \subseteq \mathcal{P}(X)$, we define $\cap \Delta= \cap\{T \,| \, T\in \Delta\}$ and 
 $\cup \Delta= \cup\{T \,| \, T\in \Delta\}$. 
It is trivial that  
$(\, \Delta \subseteq \Lambda \,)\Rightarrow(\, \cap \Delta \supseteq \cap \Lambda \,)$ and that $(\, \Delta \subseteq \Lambda\,) \Rightarrow (\, \cup \Delta \subseteq \cup \Lambda \,)$. Also  $\mathbb{R},\mathbb{Z},\mathbb{N},$ are the sets of  reals, integers, natural numbers correspondingly (note that 
$0\in \mathbb{N}$), and $\mathrm{T_n}=\{1,2,\dots n \}$.

\item \label{introducrtion2}   By  \underline{$\omega=|\mathbb{N}|$}, we denote the infinite countable cardinal and by $\mathrm{Card}$  we denote class of all cardinals. 
 It is well known  that $|\mathbb{R}|=2^{\, \omega}$ and we call this cardinal as the \emph{continuum}, which we denote by $\mathfrak{c}$. So $\mathfrak{c}=\big|\mathbb{R}\big|=2^{\, \omega}$. 
 Let  $\mathcal{L}$ be a definable class of cardinals and $\lambda$ be a cardinal. When we, somehow abusively, write that $[\, \lambda \in \mathcal{L} \,]$, we do not mean that  this is a valid formula in $\mathsf{ZFC}$ (it can't be valid unless the class $\mathcal{L}$ is a set) but that $[\, \lambda$ satisfies the defining property of the class      $\mathcal{L} \, ]$.  We define  $\mathrm{Card}_{\, > \lambda} =\{ \mu \in \mathrm{Card}\, | \, \mu > \lambda \}$ and  $\lambda^{+}=\mathrm{min}(\mathrm{Card}_{\,> \lambda})$. The sequences  $(\aleph_n)_{n\in \omega}$,  $(\lambda_{\cdot n})_{n \in \omega}$ are  defined inductively by the formulas
$[\, \aleph_0=\omega,\, \aleph_{n+1}=(\aleph_n )^{+}  \,]$,  $[\, \lambda_{\cdot 0}=\lambda \, ,\,  \lambda_{\cdot(n+1)}=2^{\, \lambda_n}\,]$. 
  Also $\lambda^{\bullet}=\underset{n\in \omega} \Sigma \lambda_{\cdot n} $. Note that the \emph{beth} notation $\beth_{\,n}(\lambda)=\lambda_{\cdot n}$ and $\beth_{\,\omega}(\lambda)=\lambda^{\bullet}$ is commonly used.
 \item Note that  $\mathrm{Card}$ is not a set. If $\kappa,\lambda$ are cardinals then $[\,(\kappa<\lambda)\lor  (\kappa=\lambda) \lor (\kappa> \lambda) \,] $. Let
    $\mathcal{L},\, \mathcal{K} $ be  two classes of cardinals. Then  $\mathcal{L}'\overset{\text{def}}=\mathrm{Card}\smallsetminus \mathcal{L}$. If  $\mathcal{L} \neq \emptyset$ then $\mathcal{L} $ has a minimum element which we denote by $\mathrm{min} (\mathcal{L})$. Note that $ [\,\mathcal{L}\subseteq \mathcal{K}\,] \Rightarrow [\, \mathrm{min} (\mathcal{K}) \leq \mathrm{min} (\mathcal{L})\,]$.
 \item  Let $A,\,B$ be sets and $f\in \mathrm{Hom}_{\text{sets}}(A,B)$ (i.e $f:A\rightarrow B$ is a function). Then $:$ \\
$\mathrm{(i)}\,\mathrm{(a)}$ $A$ is  finite iff $|A| <\omega$. 
$\mathrm{(b)}$ $A$ is  countable  iff $|A| \leq \omega$. 
$\mathrm{(c)}$ $A$ is  uncountable   iff $|A| > \omega$. \\ 
$\mathrm{(ii)}\,\mathrm{(a)}$ $[\,f$ is  $1$-$1\,]\Rightarrow[\,|A|\leq |B|\,] $
$\mathrm{(b)}$ $[\,f$ is  onto$\,]\Rightarrow[\,|A|\geq |B|\,] $
$\mathrm{(c)}$  $[\,f$ is  bijection$\,]\Rightarrow[\,|A|= |B|\,] $\\
$\mathrm{(iii)}$
$[\,|A|\leq |B|\,]\Leftrightarrow[\,\, \exists \, g\in \mathrm{Hom}_{\text{sets}}(A,B):g$ is $1$-$1\,]$. \\
$\mathrm{(iv)}$
$[\,|A|\geq |B|\,]\Leftrightarrow[\,\, \exists \, g\in \mathrm{Hom}_{\text{sets}}(A,B):g$ is onto$\,]$.  \\ 
$\mathrm{(v)}$
$[\,|A|= |B|\,]\Leftrightarrow[\,\, \exists \, g\in \mathrm{Hom}_{\text{sets}}(A,B):g$ is bijection$\,]$. \\
$\mathrm{(vi)}$ $\mathrm{(a)}\,\, [\,f$ is $1$-$1\,]\Leftrightarrow[\, f$ is left invertible$\,]$ \quad\quad\quad\quad\quad\quad\quad
 $\mathrm{(b)}\,\, [\,f$ is onto $\,]\Leftrightarrow[\, f$ is right invertible$\,]$.\\ Let also $g\in \mathrm{Hom}_{\text{sets}}(B,C)$. Consider the composition  $g\circ f\in \mathrm{Hom}_{\text{sets}}(A,C)$. Then $:$ \\
$\mathrm{(vii)}$ $\mathrm{(a)}\,[\, (f$ is $1$-$1)$, $(g$ is $1$-$1)\,]\Rightarrow$    $(g\circ f$ is $1$-$1)\,]$ \quad
 $\mathrm{(b)}\,[\, (f$ is onto), $(g$ is  onto)$\,] \Rightarrow$    $(g\circ f$ is onto$\,]$.\\
$\mathrm{(viii)}$ $\mathrm{(a)}$ $[\,g\circ f$ is $1$-$1\,]\Rightarrow$    $[\, f$ is $1$-$1)\,]$ \quad\quad\quad\quad\quad\quad
 $\mathrm{(b)}$  $[\,g\circ f$ is  onto$\,] \Rightarrow[\,g$ is
  onto$\,]$.

\end{enumerate}

\begin{theorem} \label{cardinals1}
 Let $\kappa, \lambda ,\mu$ be cardinals, $I$ be a set and
  $[\, A_i  \in \mathrm{Sets},\, \{ \nu_i,\, \xi_i  \} \subseteq \mathrm{Card}\,]$  $\forall  \in \,  I$. Then $:$ \\
 $\mathrm{(i)}$  $\kappa\cdot \lambda=\kappa+\lambda=\max \{  \kappa,\lambda \}$, if $[\, ( \kappa \geq \omega) \lor (\lambda \geq \omega) \,]$. \quad
   $\mathrm{(ii)}$  $\kappa^{\lambda}=\kappa$    if $[\, \kappa \geq \omega,\,\,  1 \leq \lambda <\omega\,]$.  \\ 
   $\mathrm{(iii)} $ 
  $(\kappa^{\lambda})^{\, \mu}=\kappa^{\lambda\cdot \mu}$
  \quad  $\mathrm{(iv)}$ $\kappa^{\lambda +\mu}=\kappa^{\lambda }\cdot \kappa^{\mu}$         
  \quad $\mathrm{(v)}$ $ [\, \kappa \ge \lambda \,] \Rightarrow [\, \kappa^{\,\mu}\ge \lambda^{\, \mu} \,]$  \ \quad $\mathrm{(vi)}$ 
  $[\, \lambda \ge 2,\, \kappa \geq 1\,]\Rightarrow [\, \lambda^{\kappa}\ge \lambda\,]$.  \\  $\mathrm{(vii)}$  $[\, \kappa \le \lambda  \le \kappa \,]\Rightarrow [\, \kappa=\lambda \,]$ $($Schröder–Bernstein Theorem$)$. \space $\mathrm{(viii)}$ $2^{\, \kappa} > \kappa$  $($Cantor's Theorem$)$.   \\                                   
$\mathrm{(ix)}$ $|   \mathcal{P}_{\text{fin}}(I) |=|  
 I |$,
  if  $|  I | \geq
   \omega$ , 
 where  $ \mathcal{P}_{ \text{fin}}(I)=$ 
 $\{ X\in  \mathcal{P}(I) :  |X |< \omega   \}$.\\ 
$\mathrm{(x)}$ If $[\,  \nu_i < \xi_i\,]\,  \forall \,  i \in I $   then  $ \underset{i\in I }\Sigma  \nu_i < \underset{i\in I} \Pi  \xi_i$ $($Kőnig's inequality$)$.\quad\quad\quad\quad $\mathrm{(xi)}$  $|\underset{i\in I} \cup A_i |$
 $\le \underset{i\in I} \Sigma
   | A_i | $.
\end{theorem}

\begin{definition}\label{def.filter} Let $X$ be a set, $\mathcal{D}\subseteq \mathcal{P}(X)$ and $\kappa, $  $\mu$ be  cardinals.\\
 $\mathrm{(i)}\,\, \kappa$ is regular $\Leftrightarrow \{$
$[\, (\lambda_i<\kappa)\, \forall \, i\in I,\, \big| I \big|<k \,]\Rightarrow$
$[\, \underset{i\in I}\Sigma \lambda_i <\kappa \,] \}$. \\
 $\mathrm{(ii)}\,\,\kappa$  is strongly inaccessible $\Leftrightarrow \{ \mathrm{(a)} \, \kappa >\omega, \,\mathrm{(b)} [\,(\lambda<\kappa) \Rightarrow (2^{\lambda}<\kappa)\,]\, \forall \, \lambda \in \mathrm{Card}, $
$\mathrm{(c)} \, \kappa$ is regular$\}$.\\
 $\mathrm{(iii)}\,\,\mathcal{D}$ is a filter on $X \Leftrightarrow$  $ \{\, \, \mathrm{(a)}\,[\, \emptyset\notin \mathcal{D},\, X \in \mathcal{D}\,]$,  $\,\,\mathrm{(b)}\,[\,  A\in \mathcal{D}$, $A\subseteq Y\subseteq X \,]\Rightarrow Y\in D$,  
$\,\,\mathrm{(c)}\,[\, A\in \mathcal{D} $, $B\in \mathcal{D} \,]\Rightarrow [\,(A\cap B)\in  \mathcal{D}\,] \, \}$ \\
$\mathrm{(v)}\,\,\mathcal{D}$ is ultrafilter on $X \Leftrightarrow$
 $\mathcal{D}$ is a maximal (with respect to inclusion) filter  
  on $X$. \\
 $\mathrm{(vi}\,\, \mathcal{D}$ is a principal filter on $X\Leftrightarrow \,[\, \exists \, A	\in \mathcal{P}(X):\mathcal{D}=\{ Y\in \mathcal{P}(X)\big|  \, A\subseteq Y \subseteq X \}\,]$.  \\
$\mathrm{(vii)}\,\,\mathcal{D}$ is a $\kappa$-closed $\Leftrightarrow \{ $ 
 $[\,( A_i \in  \mathcal{D})\, \forall \, i \, \in \kappa \,]\Rightarrow [\, (\cap \{A_i \,| \, i\in \kappa  \})\in \mathcal{D}\,]\, \} $  \\
$\mathrm{(viii)}\,\, \mathcal{D}$ is a $\kappa$-complete  $\Leftrightarrow [ \, \mathcal{D}$  is $\lambda$-closed for any cardinal $\lambda<\kappa \,]$. \\
$\mathrm{(ix)}\,\,\kappa$ is a measurable cardinal $\Leftrightarrow
[$There is a $\kappa$-complete ultrafilter on $\kappa$
that is not principal and $\kappa > \omega\,]$.\\
$\mathrm{(x)}\,\,\kappa$ is an $\omega$-measurable cardinal $\Leftrightarrow[$There is an $\omega$-closed ultrafilter  on $\kappa$ that is not principal and $\kappa > \omega\,]$. \\
$\mathrm{(xi}\,\,\mathcal{L}_1$, $\mathcal{L}_2$, $\mathcal{L}_3$ are the classes of $\omega$-measurable, measurable and  strongly inaccessible cardinals.\\
$\mathrm{(xii)}\,\,\mathcal{L}'_1=\mathrm{Card}\smallsetminus  \mathcal{L}_1$, $\mathcal{L}'_2=\mathrm{Card}\smallsetminus  \mathcal{L}_2$ are the classes of non-$\omega$-measurable, non- measurable cardinals.
\end{definition}
\begin{remark} \label{remark.advanced}Let $X$ be a set, $\mathcal{D}\subseteq \mathcal{P}(X)$ and $\kappa, $  $\mu$ be  cardinals.\\
$\mathrm{(i)}$  It is trivial that $[\, \mathcal{L}_2\subseteq \mathcal{L}_1\,]\,(1)$   and it is not so trivial that  $[\, \mathcal{L}_2\subseteq \mathcal{L}_3\,]\,(2)$. See pages 26-27 of \cite{eklof}.\\
$\mathrm{(ii)}$ Assume that \underline{$\mathcal{L}_2 \neq \emptyset$}. It follows from $\mathrm{(i)}$ that   $[\, \mathcal{L}_2 \neq \emptyset$, $\mathcal{L}_3 \neq \emptyset\,]$.  Let 
$\mu_i=\mathrm{min}(\mathcal{L}_i)$,  for $\i \in \{1,2,3\}$.
It follows from $\mathrm{(i)}$ that $[\, \mu_2 \geq \mu_1,\,  \mu_2 \geq \mu_3\,]\,(3)$. It is  known that 
$[\, \mu_1 \in \mathcal{L}_2 \,]\,(4)\Rightarrow$  
$[\, \mu_1 \geq\mu_2\,]\Rightarrow $
$[ \, \mu_1=\mu_2 \, ]$. We define   $[\, \zeta=\mu_1=\mu_2\,]\,(5),\,$ as the first $\omega$-measurable cardinal (if $\mathcal{L}_2 \neq \emptyset$).\\
$\mathrm{(iii)}$ When we say that $[\, \kappa$ is less that the first $\omega$-measurable cardinal$\,]$, we do not necessarily assume that $\mathcal{L}_2 \neq \emptyset$. We simply  mean that   $[\, (\mathcal{L}_2 = \emptyset)$  or  $(\, \mathcal{L}_2 \neq \emptyset$, $\kappa < \zeta)\,]$. In this case $\kappa$ is not $\omega$-measurable and any smaller cardinal that $\kappa$ is not  $\omega$-measurable. \\
$\mathrm{(iv)}$  It is known that $\mathrm{(a)} [\,( \mathcal{L}_1 \neq \emptyset)  \Leftrightarrow    
(\mathcal{L}_2 \neq \emptyset)\,]$, and  that $\mathrm{(b)}\, \mathsf{ZFC} \nvdash  [\,\mathcal{L}_1 \neq\emptyset  \,]$ (see pages 26-27 of \cite{eklof}).
$\mathrm{(v)}$ It follows from $\{\ref{def.filter}\mathrm{.ii},\,(2),\,(4),\,(5)    \}$ that  $[\, (\kappa < \zeta)\Rightarrow(2^{\kappa} < \zeta)\,]$ and that $[\, \zeta$ is regular$\,]$. \end{remark}

\begin{lemma}\label{simple} Let $ \lambda\in \mathrm{Card}$ and $\zeta$ be the first $\omega$-measurable cardinal. Then  $[\, (\lambda < \zeta)\Rightarrow (\lambda^{\bullet} < \zeta)\,]$.
\end{lemma}

\begin{proof}It follows from  $\ref{remark.advanced}\mathrm{.v}$ that $[\, \zeta$ is regular$\,]\,(2)$ Let $A=\{n\in \mathbb{N}\, | \, \lambda_{\cdot n} < \zeta  \}$. Then $0\in A$ since $\lambda_{\cdot 0}=\lambda < \zeta$, by assumnption. $\cdot$ If $[\, n\in A\,]$ then $[ \lambda_{\cdot n} < \zeta\,]\overset{\ref{remark.advanced}\mathrm{.v}}\Longrightarrow$
$[\, 2^{\, \lambda_{\cdot n}}< \zeta \,]\Rightarrow[\, \lambda_{\cdot n+1} <\zeta\,]\Rightarrow[\, (n+1) \in A\,]$. It follows by induction that $A=\mathbb{N}$. Hence $[\,( \lambda_{\cdot n} < \zeta\,)\, \forall \, n \in \mathbb{N}\,]\,(3)\Rightarrow$
$[\, \lambda^{\bullet}=\underset{n\in \mathbb{N}}\Sigma \lambda_{\cdot n}\leq \underset{n\in \mathbb{N}}\Sigma \zeta=\omega \cdot \zeta=\zeta\,]\Rightarrow[\,\lambda^{\bullet} \leq \zeta\,]\,(4) $. It also follows from the regularity of $\zeta$ that  $[\,   \zeta  \neq \lambda^{\bullet}=\underset{n\in \mathbb{N}}\Sigma \lambda_{\cdot n}        \,]\,(5)$ and from $\{(4),\,(5)\}$ that $[\, \lambda^{\bullet} < \zeta\,]$. \end{proof}

 \subsection{Duals of free modules and extensions of modules}
 All rings $R$ in this article are unital and non-trivial, i.e they  have an identity element $1_{R}\neq 0_R$ (So $\mathrm{card}(R)\ge 2$) and all modules are left and unitary. The term  \underline{r.h$_{}$} means unital ring homomorphism  i.e ring homomorphism $\phi:R\rightarrow S$ such that  $\phi(1_R)=1_S$.  The terms \underline {f.g,$_{}$ f.p} are abreviations of the terms finitely generated, finitely presented  respectively.  By $R$-$\mathrm{Mod}$, we denote the class of (left) $R$-modules and by  $R$-$\mathrm{Mod}_{\, \mathrm{f.g}}$,  $R$-$\mathrm{Mod}_{\, \mathrm{f.p}}$, we denote the classes of f.g, f.p $R$-modules respectively. Also  $\hat{R}$
is the opposite ring of $R$ and  $\mathrm{Mod}$-$R=\hat{R}$-$\mathrm{Mod}$ is the class of right $R$-modules. If  $M\in R$-$\mathrm{Mod}$ then   $\mathrm{Sub}_R(M)$ is the set of  all $R$-submodules of $M$ and   $\mathrm{Sub}_R(R)$ is the set of all left ideals of $R$.

\begin{definition}\label{def.cofree}Let $R=(R,+,\cdot)$ be a ring and $X$ be 
 set and $M\in R$-$\mathrm{Mod}$. Then $:$ \\
$\mathrm{(i)}\,$ $R^X=\mathrm{Hom}_{\mathrm{sets}}(X,R)=$ The set of all   functions from $X$ to $R$. We often denote $y\in R^X$ by
 $y=y(i)_{i\in X}=(y_i)_{i\in X}$. Also define $\bold{+}:R^X \times R^X \rightarrow R^X$ by $[\, (y_i)_{i\in X} \bold{+} (z_i)_{i\in X}=(y_i+z_i)_{i\in X}\,]$ and  $\bold{\cdot}:R \times R^X \rightarrow R^X$ by $[\, r \bold{\cdot}(y_i)_{i\in X} =(r \cdot y_i)_{i\in X}\,]$.
 Trivially the abelian group $(R^X,\bold{+})$ becomes an $R$-module by the action of $\bold{\cdot}$. So $R^X\in R$-$\mathrm{Mod}$. \\
$\mathrm{(ii)}\,$  $R^{[X]} =\{y\in R^X$ 
$ \big| \, \mathrm{supp}(y)\,\, \text{is finite} \}$, where 
$\mathrm{supp}(y)= \{i\in R^X | \, y_i \neq 0 \}$. Trivially $R^{[X]}$ is an $R$-submodule of $R^X$. For  $ R^{[X]}\ni y=(y_i)_{i\in X}$, we define $ \hat{\xi}_y \in \mathrm{Hom}_R(R^{\lambda},R)$   by $\hat{\xi}_y(x)= \mathcal{h} x, y \mathcal{i}=\underset{i\in X} \Sigma
x_i \cdot y_i$. \\
$\mathrm{(iii)}\,$   Let  $i\in X$. Then  $e_{i}=e^	X_{i,R} \in R^{[X]}$ is defined by $\{ e_i(j)=1 \,\, \text{for}\,\, j=i \,\,\text{and}\,\,e_i(j)=0 \,\,\text{for} \,\, 
 j\ne i \}$. We can write  $y=(y_i)_{i\in X}$ as $ y=\underset{i\in \kappa}  \Sigma y_i\, e_i \in R^{[X]} $, since $\mathrm{supp}(y)$ is finite. \\
 $\mathrm{(iv)}$ The operation  $\hat{\cdot}:R \times R\rightarrow R$
is defined by $r \, \hat {\cdot}\, s= s\cdot t$. The triple $R=(R,+,\hat{\cdot})$ forms a ring, which we call as the oppossite ring of $R$ and denote by  $\hat{R}=R^{\, \text{op}}$.
 \end{definition}

\begin{definition}  \label{def.restriction} Let  
  $\phi: R\rightarrow S$ be a r.h. Let also  
 $\{M,\,N\}\subseteq$  $S$-$\mathrm{Mod}$ and $g\in\mathrm{Hom}_R(M,N)$. Then $:$ \\
$\mathrm{(i)}$ The $\mathbb{Z}$-module $M$ becomes an   $R$-module by the action $r \ast m=f(r)\cdot m$, which we denote by
 $_{\phi}M$. \\
$\mathrm{(ii)}$ The map $_{\phi}g:\, _{\phi}M \rightarrow \, _{\phi}N$ is defined by the formula $_{\phi}g(m)=m$. Trivially $_{\phi}g\in \mathrm{Hom}_R(_{\phi}M,\, _{\phi}N)$.\\
$\mathrm{(iii)}$ We call \underline{$_{\phi}M$ as the restriction of $M$ along $\phi$} and  \underline{$_{\phi}g$ as the restriction of 
$g$ along $\phi$}.\\
$\mathrm{(iv)}$ The Functor $K_{\phi}:S$-$\mathrm{Mod}\rightarrow$
$R$-$\mathrm{Mod}$ is defined by $K_{\phi}(M)=$ $ _{\phi} M$ and
 $K_{\phi}(g)=$ $ _{\phi}g$. \\
$\mathrm{(v)}$ We say that $\phi$ is f.p, f.g, flat, projective, injective iff  $K_{\phi}(S)=$ $ _{\phi} S\in R$-$\mathrm{Mod}$ is so. 
 \end{definition}
\begin{definition}  \label{def.ext} Let  $R,\,S,\,T$ be  rings, $\{M,\,N\}\subseteq$  $R$-$\mathrm{Mod}$, $f\in
 \mathrm{Hom}_R(M,N)$ and  $\phi: R\rightarrow S$ be a r.h. Let also
 $\{M',\,N'\}\subseteq R$-$\mathrm{Mod}$, $f'\in
 \mathrm{Hom}_R(M',N')$ and  $\psi: S\rightarrow T$ be a r.h. We define $:$ \\
 $\mathrm{(i)}$  Then $S$ is   an $(R$-$R)$-bimodule by the the action $r_1 \ast s \ast r_2 = \phi(r_1) \cdot s \cdot \phi(r_2)$. Hence the abelian group  $ S \otimes_R M$    becomes an $S$-module            by the action  $s*(\overset{n}{\underset{k=1}\Sigma}s_k \otimes m_k)=\overset{n}{\underset{k=1}\Sigma}(s\cdot s_k) \otimes m_k$, which we call    as\underline{ the extension of $M$ along $\phi$} and denote by $M_{\phi}$. So  $[\, M_{\phi} =S \otimes_R M\in S$-$\mathrm{Mod}\,]$.   We also define define the extension of $f$ along 
 $\phi$ by $f_{\phi}\overset{\text{def}}=\mathrm{id}_S \otimes f:M_{\phi}\rightarrow N_{\phi}$.  The Functor $\Lambda_{\phi}:$
  $R$-$\mathrm{Mod}\rightarrow$ $S$-$\mathrm{Mod}$ defined by 
  $[\, \Lambda_{\phi}(M)=M_{\phi},\,  \Lambda_{\phi}(f)=f_{\phi}\,] $       is the  the extension functor along $\phi$. \\
 $\mathrm{(ii)}$ Let $J \in \mathrm{Sub}_R(R)$ i.e $J$ a left ideal of $R$. It it trivial that $\phi(J)\in \mathrm{Sub}_S(S)$. Also if
 $J=\langle W \rangle_R$ then  $\phi(J)=\langle \phi(W) \rangle_S$. We call the left ideal $\phi(J)$ of $S$ as \underline{the  extension of $J$ via $\phi$} and denote it by $J^{\phi}$. Note that $J_{\phi}  \neq J^{\phi}$. \\
$\mathrm{(iii)}$  $[\, f \underset{R}\sim f' \,]\overset{\text{def}}\Leftrightarrow$
$[\, \exists\, h\in \mathrm{Hom}_R(M,M')\,]$
$[\, \exists\, h'\in \mathrm{Hom}_R(N,N')\,]:\{\,  [\, h'\circ f=f'\circ h\,] \, $  and   $[\, \text{both}\,  f ,f'$ are bijections$\,]\, \}$. We then say that $f$ is $R$-equivallent to $f'$. \\
$\mathrm{(iv)}\,\mathrm{(a)}$  $[\,M\underset{R}\rightarrowtail N\,] \overset{\text{def}}\Leftrightarrow$
$[\, \exists \, f\in \tmop{Hom}_R(M,N):f \,\,\text{is 1-1}\,]$.\space  
$\mathrm{(b)}\,\,[\,M\underset{R} \twoheadrightarrow N\,] \overset{\text{def}}\Leftrightarrow$
$[\, \exists\,  f\in \tmop{Hom}_R(M,N):f \,\,\text{is onto}\, ]$.
 \end{definition}\label{remark.ext}

\begin{remark} \label{remark.ext} Consider the setting  of  \ref{def.ext}.
Let also $g,\, g'$ be morphisms in $R$-$\mathrm{Mod}$ such that 
$f\circ g$, 
$f'\circ g'$ are defined. Then $:$  \\
$\mathrm{(i)}$ Note that $M_{\phi}$ is  an $R$-module as well. Also  $f_{\phi} \in  \mathrm{Hom}_S(M_{\phi},N_{\phi})$ and  $f_{\phi} \in  \mathrm{Hom}_R(M_{\phi},N_{\phi})$. \\ 
$\mathrm{(ii)}$ It is trivial that $\underset{R}\sim$ is an equivallence relation on the the class of morphisms of  $R$-$\mathrm{Mod}$. Also  $:$  \\
 $\mathrm{(a)}\,\,[\, f \underset{R}\sim f'\,] \Rightarrow [\, \mathrm{Ker}(f)\underset{R}\cong  \mathrm{Ker}(f'),\,\,  \mathrm{Im}(f)\underset{R}\cong  \mathrm{Im}(f'),\, \,  \mathrm{Cok}(f)\underset{R}\cong  \mathrm{Cok}(f')\,]$ \quad $\mathrm{(b)}\,\, (f\circ g)_{\phi} \underset{S}\sim f_{\phi}\circ g_{\phi} $ \\
 $\mathrm{(c)}\,\,[\, f \underset{R}\sim f',\, \, g \underset{R}\sim g'] \Rightarrow$ $[\, f\circ g \underset{R}\sim
f'\circ g'  \,]$ \quad\quad $\mathrm{(d)}\,\, (f_{\phi})_{\psi} \underset{T}\sim f_{\psi \circ \phi}$.\\
 $\mathrm{(iii)}$   The Functors $\Lambda_{\phi}:$
  $R$-$\mathrm{Mod}\rightarrow$ $S$-$\mathrm{Mod}$  and  $K_{\phi}$:
  $S$-$\mathrm{Mod}\rightarrow$ $R$-$\mathrm{Mod}$ are additive. \\
 $\mathrm{(iv)}$  Trivially  $\mu_R(M)=\mathrm{min}\{\kappa \in \mathrm{Card}\, |\,
R^{[\kappa]}\underset{R} \twoheadrightarrow M  \}$. Also  $M$ is a  f.g $R$-module iff $\mu_R(M) < \omega$. \\
$\mathrm{(v)}$ $\mathrm{Hom}_R(M,R)$ becomes an $\hat{R}$-module by the action $(r\ast f) (m)=f(r\cdot m)$. 
\end{remark}

 \begin{definition}  \label{def.dual} Let $R=R$ be a ring, $\{M,\, N\}\subseteq R$-$\mathrm{Mod}$ and $f\in \mathrm{Hom}_R(M,N)$. We define  $:$  \\
 $\mathrm{(i)}$ The $R$-dual of $M$ by $M^{\ast} = \mathrm{Hom} (M, R)\in \hat{R}$-$\mathrm{Mod}$. \\ $\mathrm{(ii)}$ The $R$  bidual of $M$    
 by $M^{\ast\ast}=(M^{\ast})^{\ast} \in R$-$\mathrm{Mod}$. \\
$\mathrm{(iii)}$  $\mathrm{S}_{M} : M \rightarrow M^{\ast \ast}$ is  defined by $ \mathrm{S}_{M} (m) (g) = g (m)$, for $m\in M$ and $g\in M^{\ast}$. \\
$\mathrm{(iv)}\,\mathrm{(a)} \, [\, M$ is $R$-torsionless$\,]\Leftrightarrow [\, \mathrm{S}_{M}$ is 1-1$\,]$ \quad\quad\quad
$\mathrm{(b)} \, [\, M$ is $R$-reflexive $\,]\Leftrightarrow [\, \mathrm{S}_{M}$ is a bijection$\,]$ \\
$\mathrm{(v)}$ The map $f^{\ast}=\mathrm{Hom} (f, R):N^{\ast}\rightarrow M^{\ast}$ is defined by $f^{\ast}(h)=h\circ f$. Trivially 
$f^{\ast}\in  \mathrm{Hom}_{\hat{R}}(N^{\ast},M^{\ast})$.
\end{definition}
 
\begin{remark} \label{remark.dual} Consider the setting of  Definition \ref{def.dual}. \\
$\mathrm{(i)}$  By $\mathrm{Refl}(R)$, we denote the class of reflexive $R$-modules. The term  $\sigma_{M}$ is used instead of $\mathrm{S}_{M}$.
in \cite{fuller}  \\
$\mathrm{(ii)}$ It follows from $\ref{remark.ext}\mathrm{.v}$ that 
$M^{\ast} \in \hat{R}$-$\mathrm{Mod}$. Hence $M^{\ast\ast} \in 
\hat{\hat{R}}$-$\mathrm{Mod}=R$-$\mathrm{Mod}$. Also $\mathrm{S}_M\in \mathrm{Hom}(M,M^{\ast\ast})$. \\
 Hence $[\,M \in \mathrm{Refl}(R)\,]\Leftrightarrow [\,\mathrm{S}_M$ is $R$-isomorphism$\,]\Rightarrow [\,M^{\ast\ast} \cong_R M\,]$.  \\
$\mathrm{(iii)}$ Note that  $f^{\ast}(h)=h\circ f$. So $[\, f$ onto$\,]\Rightarrow$  $[\, f^{\ast}$ is 1-1$\,]$. Therefore 
$[\,M\underset{R} \twoheadrightarrow N\,] \Rightarrow$  
$[\,N^{\ast}\underset{\hat{R}} \rightarrowtail M^{\ast}\,]$.    \\  
$\mathrm{(iii)}$ Let $ \kappa\in \mathrm{Card}$. Trivially $[\, R^{\ast}  \cong_{\hat{R}} \hat{R} \,]\,(1)$. So $[\, (R^{[\kappa]})^{\ast}=\mathrm{Hom}(\underset{i \in \kappa}\Sigma R ,R)\cong_{\hat{R}}$
$\underset{i \in \kappa}\Pi \mathrm{Hom}(R,R)$
 $\overset{(1)} \cong_{\hat{R}}  (\hat{R})^{\kappa}\,]\Rightarrow$  
 $[\, (R^{[\kappa]})^{\ast} \cong_{\hat{R}}  (\hat{R})^{\kappa}\,]$. Hence $[\, ((\hat{R})^{[\kappa]})^{\ast} \cong_{\hat{R}}  R^{\kappa}\,]$. In particular   $[\, (R^{[\kappa]})^{\ast} \cong_{\hat{R}}  R^{\kappa}\,]$, when $R$ is commutative.  \\ 
 $\mathrm{(iii)}$ Assume that $[\, M$ is $R$-torsionless$\,]\,(4)$ and  $[\lambda=\mu_R(M^{\ast})\,]$. Therefore we have  $[\,(\hat{R})^{[\lambda]}\underset{\hat{R}} \twoheadrightarrow M^{\ast}\,]\overset{\mathrm{(iii)}}\Longrightarrow$
$[\,   M^{\ast\ast}  \underset{R}   \rightarrowtail   ((\hat{R})^{[\lambda]})^{\ast}  \underset{R} \cong R^{\lambda} \,]\Rightarrow$
$[\,   M^{\ast\ast}  \underset{R}   \rightarrowtail  R^{\lambda} \,]\,(5)$. Also  $(4)\overset{\ref{def.dual}\mathrm{.iv}}\Longrightarrow$
$[\,  M  \underset{R}   \rightarrowtail  M^{\ast\ast} \,]\,(6)$. It follows from $\{(5),\, (6)  \}$ that $[\,  M  \underset{R}   \rightarrowtail R^{\lambda}  \,]$. Therefore  \underline{any torsionless $R$-module can be embedded in a cofree one.} 

\end{remark}

\begin{lemma}Let  $R$ be a ring, $\{M,N,\}\subseteq R$-$\mathrm{Mod}$, $f\in
 \mathrm{Hom}_R(M,N)$ and  $\phi: R\rightarrow S$ be a r.h.
  Then $:$ \\
$\mathrm{(i)}$ $[f$ is $R$-left invertible$]\Rightarrow$               $[f_{\phi}$ is  $S$-left invertible$\,]$ 
$\mathrm{(ii)}$ $[f$ is $R$-right invertible$]\Rightarrow$               $[f_{\phi}$ is  $S$-right invertible$\,]$  \\
$\mathrm{(iii)}$ $[f$ is onto $]\Rightarrow$    $[f_{\phi}$ is  $S$-onto$\,]$  \quad\quad\quad\quad\quad\quad\quad $\mathrm{(iv)}$ $[f$ is $R$-isomorphism $]\Rightarrow$ $[f_{\phi}$ is  
$S$-isomorphism$\,]$  \\\
$\mathrm{(v)}$    $\,[f$ is $1$-$1$, $\phi$ is  flat r.h $]\Rightarrow$ $[f_{\phi}$ is $1$-$1\,]$ \\
 
\end{lemma}

\begin{lemma}\label{rh}Let $R,\,S$ be  rings, $\lambda$ be a cardinal $($or a set$)$ and $\phi: R\rightarrow S$ be  a r.h.\\
$\mathrm{(i)}\,(R^{[\lambda]})_{\phi}\underset{S}\cong S^{[\lambda]}$ \quad\quad
$\mathrm{(ii)}$  $(R^{\lambda})_{\phi}\underset{S}\cong S^{\lambda}$,
when  $\phi $ is f.p \quad\quad
$\mathrm{(iii)}$
$(R^{\lambda})_{\phi}\underset{S} \twoheadrightarrow S^{\lambda}$, when $\phi $ is f.g .\end{lemma}
\begin{proof} Follows from Propositions $3.22$,$3.23$ of \cite{enochs}. \end{proof}

\begin{lemma} \label{trivial1}Let   $R$ be a ring  and $\kappa, \lambda$ be cardinals and $X,Y$ be sets. Then $:$ \\
$\mathrm{(i)\,}$  If $|X|=|Y|$ then $:$ \quad\quad $\mathrm{(a)}\,\, R^{X}\cong_R R^{Y} $     \quad\quad $\mathrm{(b)}\,\, R^{[X]}\cong_R R^{[Y]} $         \\
$\mathrm{(ii)\,}R^{\kappa+\lambda}\cong_R (R^{\kappa} \oplus R^{\lambda})$ \quad\quad $\mathrm{(iii)\,}R^{[\kappa+\lambda]}\cong_R (R^{[\kappa]} \oplus R^{[\lambda]})$    \quad\quad 
$\mathrm{(iv)\,}$ If $k \leq \lambda$ then $\exists\, \mu \in \mathrm{Card}:[\,\lambda=\kappa +\mu\,]$.  \\
$\mathrm{(v)}$ If $\kappa \leq \lambda$ then $:$ $\mathrm{(a)}\,[\,  R^{\kappa} \underset{R} \rightarrowtail R^{\lambda}\,]$ \quad
$\mathrm{(b)\,}[\,  R^{[\kappa]} \underset{R} \rightarrowtail R^{[\lambda]}\,]$ \quad $\mathrm{(c)\,}[\,  R^{\lambda} \underset{R} \twoheadrightarrow R^{\kappa}\,]$   \quad $\mathrm{(d)\,}[\,  R^{[\lambda]} \underset{R} \twoheadrightarrow R^{[\kappa]}\,]$ \\
\end{lemma}
\begin{proof}\underline{$\mathrm{Proof\, of\, (i)}$} : $|X|=|Y|\Rightarrow \exists \, f\in \mathrm{Hom}_{sets}(X,Y):[\, f$ is bijection$\,]$. Let   $H_{f}:R^{X} \rightarrow R^{Y}$ be the function defined by   the formula  $[\, H_{f}((r_x)_{x \in X})=(f(r_x))_{x\in X}\,]$. Trivially $H_{f}$ is $R$-linear and $H_{f}$
is a bijection since $f$ is so. So $R^{X}\cong_R R^{Y} $ since $H_{f}$ is $R$-isomorphism. Similarly we get that $R^{[X]}\cong_R R^{[Y]}$.  \\ 
 In a similar way, we can prove  $\mathrm{(ii)}$. Also    
 $\mathrm{(iii)}$  and $\mathrm{(v)}$  are trivial.         \\
 \underline{$\mathrm{Proof\, of\, (iv)}$} : Note that $ [\,M  \underset{R}\rightarrowtail ( M\oplus N)\,]\,(5)$ and  $ [\, ( M\oplus N) \underset{R} \twoheadrightarrow M\,]\,(6)$, hold for any $R$-modules $M,N$.
Hence  $\{\mathrm{(iii), (i),(5)}  \}\Rightarrow \mathrm{(iv.a)}$,
$\{\mathrm{(iii), (ii),(5)}  \}\Rightarrow \mathrm{(iv.b)}$, $\{\mathrm{(iii), (i),(6)}  \}\Rightarrow \mathrm{(iv.c)}$ and $\{\mathrm{(iii), (ii),(6)}  \}\Rightarrow \mathrm{(iv.d)}$.
\end{proof}

\begin{lemma} \label{rankiswelldefined}      Any two bases of a free 
$R$-module over a commutative ring  $R$ have equal cardinalities. Equivallently for any two cardinals $\kappa,\, \lambda$ it holds that 
 $[\, R^{[\kappa]}\cong_R R^{[\lambda]}\,]\,(1) \Leftrightarrow $
 $ [\,\kappa=\lambda \,]\,(2)$.
\end{lemma}
\begin{proof}  $\cdot $ The trivial direction  $(2)\Rightarrow (1)$   is $\ref{trivial1}\mathrm{.iii}$. \quad\quad\quad\quad\quad\quad\quad\quad\quad\quad \\
\underline{$\cdot\, \mathrm{Proof\, of\, }$   $(1)\Rightarrow (2)$} :   
Let $\mathfrak{m}\in \tmop{max}(R)\neq \emptyset$ 
 and let $\phi=\pi^R_\mathfrak{m}:R\rightarrow \dfrac{R}{\mathfrak{m}}$  
given by $\phi(r)=r+\mathfrak{m}$. Then $\dfrac{R}{\mathfrak{m}}=F$ is a field. Note that $\phi$ may not be  f.p, since $R$ may not be Noetherian, but it will not matter. Therefore
  $[\, R^{[\kappa]} {\cong}_R R^{[\lambda]}\,] \Rightarrow$
 $[\, (R^{[\kappa]})_{\phi} {\cong}_F (R^{[\lambda]})_{\phi}\,] \overset{  \ref{rh}\mathrm{.i}} \Leftrightarrow$
 $[\, F^{[\kappa]} \cong_F F^{[\lambda]}\,] \Leftrightarrow $  
 $[\, \tmop{dim}_F (F^{[\kappa]})= \mathrm{dim}_F(F^{[\lambda]})\,]\Leftrightarrow$
 $[\, \kappa=\lambda\,]$. 
\end{proof}
Note that  \ref{rankiswelldefined} is not valid without the commutativity assumption on $R$ (See example 1.4 of \cite{lam}).\\
By applying the ring homomorphism $\phi$ like in 
\ref{rankiswelldefined} and using the identities of \ref{rh}, we get the  next. 
\begin{lemma}\label{rh2}Let $R,S$ be commutative rings, $\phi: R\rightarrow S$ be r.h and $\kappa, \lambda$ be cardinals $($or sets$)$ 
\end{lemma}
\begin{enumerate} [label=(\roman*)]
\item  $[\, R^{[\kappa]}\cong_R R^{[\lambda]}\,] \Rightarrow  [\, S^{[\kappa]}\cong_S S^{[\lambda]}\,],$ under no  additional assumptions for $\phi$. \\
$\mathbf{Proof:}$
 $[\, R^{[\kappa]}\cong_R R^{[\lambda]}\,] \Rightarrow [\, (R^{[\kappa]})_{\phi}\cong_S (R^{[\lambda]})_{\phi}\,]\overset{\ref{rh}\mathrm{.i}}\Longrightarrow$
$ [\, S^{[\kappa]}\cong_S S^{[\lambda]}\,]$
\item  $[\,R^{\kappa}\cong_R R^{\lambda}\,]  \Rightarrow $
$[\,  S^{\kappa}\cong_S S^{\lambda}        \, ]$, when   $\phi $  is f.p, i.e $S$ is a finitely presented 
$R$-module.  \\
$\mathbf{Proof:}$
 $[\, R^{\kappa}\cong_R R^{\lambda}\,] \Rightarrow [\, (R^{\kappa})_{\phi}\cong_S(R^{\lambda})_{\phi}\,]\overset{\ref{rh}\mathrm{.ii}}\Longrightarrow$
$ [\, S^{\kappa}\cong_S S^{\lambda}\,]$.
\item  $[\,R^{\kappa}\cong_R R^{[\lambda]}\,] \Rightarrow $
$[\,  S^{\kappa}\cong_S S^{[\lambda]}        \, ]$, when   $\phi $  is f.p, i.e $S$ is a finitely presented 
$R$-module.\\
$\mathbf{Proof:}$
 $[\, R^{\kappa}\cong_R R^{[\lambda]}\,] \Rightarrow [\, (R^{\kappa})_{\phi}\cong_S (R^{[\lambda]})_{\phi}\,]\overset{\ref{rh}}\Longrightarrow$
$ [\, S^{\kappa}\cong_S S^{[\lambda]}\,]$
\item  $[\,R^{[\kappa]}\underset{R}\cong R^{\lambda} \,]\Rightarrow$.
$ [\, S^{[\kappa]}\underset{S} \twoheadrightarrow S^{\lambda}        \, ]$, when is   $\phi $  is f.g, i.e $S$ is a finitely generated  $R$-module.\\
$\mathbf{Proof:}$
 $[\, R^{[\kappa]}\underset{R}\cong R^{[\lambda]}\,] \Rightarrow [\, (R^{[\kappa]})_{\phi}\underset{S}\cong (R^{[\lambda]})_{\phi}\,]\overset{\ref{rh}\mathrm{.i}}\Longrightarrow$
$ [\, S^{[\kappa]}\underset{S}\cong (R^{\lambda})_{\phi} \underset{S}{\overset{\ref{rh}\mathrm{.iii}}\twoheadrightarrow} S^{\lambda}\,]\Rightarrow$
$ [\, S^{[\kappa]}\underset{S} \twoheadrightarrow S^{\lambda}        \, ]$.
 \item $\mathrm{(a)}\,\,[\, R^{[\kappa]}\underset{R} \twoheadrightarrow R^{[\lambda]}\,] \Rightarrow $
$ [\, S^{[\kappa]}\underset{S} \twoheadrightarrow S^{[\lambda]}\,]$
\quad\quad
 $\mathrm{(b)}\,\,[ R^{[\kappa]}\underset{R} \rightarrowtail R^{[\lambda]}\,] \Rightarrow $
$ [\, S^{[\kappa]}\underset{S} \rightarrowtail S^{[\lambda]}\,]$,
if $\phi$ is  flat. \\
$\mathbf{Proof}:$ The proof is similar to the proof of $\mathrm{(vi)}$ below.
 \item $\mathrm{(a)}\,[\, R^{\kappa}\underset{R} \twoheadrightarrow R^{\lambda}\,] \Rightarrow $
$ [\, S^{\kappa}\underset{S} \twoheadrightarrow S^{\lambda}\,]$, if $\phi$ is f.p.
\space
 $\mathrm{(b)}\,\,[\, R^{\kappa}\underset{R} \rightarrowtail R^{\lambda}\,] \Rightarrow $
$ [\, S^{\kappa}\underset{S} \rightarrowtail S^{\lambda}\,]$,
if $\phi$ is flat and f.p. \\
\underline{$\mathrm{Proof\,of \,(a)}$}: $[\, R^{\kappa}\underset{R} \twoheadrightarrow R^{\lambda}\,] \Rightarrow $
 $[\, \exists \, f\in \mathrm{Hom}_R(R^{\kappa}, R^{\lambda})$
 : $f$ is onto$\,]\Rightarrow$  $[\, f_{\phi}=\mathrm{id}_S \otimes f$   is onto$\,]\,(6)$. It follows from $\ref{rh2}\mathrm{.ii}$ that 
 $[\, (R^{\kappa})_{\phi}\underset{S}\cong S^{\kappa}  \, ,\, (R^{\lambda})_{\phi}\underset{S}\cong S^{\lambda}\,]\,(7)$ since
 $\phi$ is f.p. It follows from (7) that $[\,  \exists \, h\in \mathrm{Hom}_S(S^{\kappa}, S^{\lambda}):(\,  h \underset{S}\sim f_{\phi}\,)\,]\,(8)\overset{(6)}\Rightarrow[\, h$ is onto$\,]\Rightarrow$
$ [\, S^{\kappa}\underset{S} \twoheadrightarrow S^{\lambda}\,]$.\\ 
\underline{$\mathrm{Proof\,of \,(b)}$}: $[\, R^{\kappa}\underset{R} \rightarrowtail R^{\lambda}\,] \Rightarrow $
 $[\, \exists \, f\in \mathrm{Hom}_R(R^{\kappa}, R^{\lambda})$
 : $f$ is 1-1$\,]\overset{\phi \text{ is\,\,flat}}\Longrightarrow$  $[\, f_{\phi}=\mathrm{id}_S \otimes f$   is 1-1$\,]\,(6')$.
 It follows from $\ref{rh2}\mathrm{.ii}$ that 
 $[\, (R^{\kappa})_{\phi}\underset{S}\cong S^{\kappa}  \, ,\, (R^{\lambda})_{\phi}\underset{S}\cong S^{\lambda}\,]\,(7')$ since
 $\phi$ is f.p.  
  It follows from $(7')$ that $[\,  \exists \, h\in \mathrm{Hom}_S(S^{\kappa}, S^{\lambda}):(\,  h \underset{S}\sim f_{\phi}\,)\,]\,(8')\overset{(6')}\Rightarrow[\, h$ is 1-1$\,]\Rightarrow$
$ [\, S^{\kappa}\underset{S} \rightarrowtail S^{\lambda}\,]$.\\
\end{enumerate}

\section{Vector spaces}
\subsection{The dimension of the dual of a vector space}
\label{section4}
We start this section with a proof of the \emph{Erdős  Kaplansky Theorem}, which is  Theorem $\ref{EKextension}\mathrm{.i}$, since many  results in this  article depend heavily on it. Note that $\omega=\mathrm{card}(\mathbb{N})$.

\begin{lemma}\label{cardinallemma1}
Let $\kappa, \lambda$ be cardinals and $R$ be abelian group such that  $\big| R \big|=\kappa$. If $[\, \kappa\geq 2,\, \lambda \geq 1\,]\,(1)$  and $[\, \kappa\geq \omega$ or $\lambda\geq \omega\,]\,(2)$  then $| R^{[\lambda]} |=\max \{ \kappa ,  \lambda  \}$.
\end{lemma}
\begin{proof}
  Assumption (1) implies  that   $[\, | R^{[\lambda]} |\geq \kappa\,]$ (3) and  
 $[\, | R^{[\lambda]} |\geq  \lambda\,]$ (4). So  $[\, | R^{[\lambda]} | \geq \max \{ \kappa ,  \lambda  \}\,]$ (5).
 Let $[\, \lambda'=\mathcal{P}_{\text{fin}}(\lambda)\smallsetminus \{\emptyset\} \,]\,(6)$  and $\xi\in \lambda' $. Then $|\{f\in R^{[\lambda]}\,|\,   \mathrm{Supp}(f)\subseteq \xi\,\}|=\kappa^{\, \xi}$. So $[\,\, | R^{[\lambda]} | \leq \underset{\xi \in \lambda'} \Sigma \kappa^{\, \xi}\,]$ (7). Let's  consider the three following cases of the assumptions (1), (2) separately .\\              
 \underline{$\mathrm{Case \,1_{}}$}: $[\, \kappa\geq \omega,\, \lambda\geq \omega \,]$ .  Then  $[\, | \lambda '|=\lambda\,]$ (8). If $\xi\in \lambda'$ then 
 $[\, | \kappa^{\, \xi}|=\kappa\,]$ (9), since $[\,\kappa\geq \omega,\, 1\leq\xi< \omega  \,]$. 
So $(7)\Rightarrow$ 
 $[\,\, | R^{[\lambda]} | \leq \underset{\xi \in \lambda'} \Sigma \kappa^{\, \xi} \underset{(9)}\leq |\lambda' |\cdot \kappa \underset{(8)}=\lambda \cdot \kappa\underset{(2)}= \max \{ \kappa ,  \lambda  \}\,]\Rightarrow$ 
 $[\, | R^{[\lambda]} | \leq \max \{ \kappa ,  \lambda  \}\,]$\\
 \underline{$\mathrm{Case \,2_{}}$}: $[\, \kappa\geq \omega,\, \lambda < \omega \,]$ .  Then  $[\, | \lambda '| <\omega\,]\,(8')$ . If $\xi\in \lambda'$ then 
 $[\, | \kappa^{\, \xi}|=\kappa\,]\,(9')$ , since $[\,\kappa\geq \omega,\, 1\leq\xi< \omega  \,]$. 
So $(7)\Rightarrow$ 
 $[\,\, | R^{[\lambda]} | \leq \underset{\xi \in \lambda'} \Sigma \kappa^{\, \xi} \underset{(9')}\leq |\lambda' | \cdot \kappa \underset{(8')} \leq \omega \cdot \kappa\underset{(2)}=\kappa \underset{(2)}=\max \{ \kappa ,  \lambda  \}\,]\Rightarrow$ 
 $[\, | R^{[\lambda]} | \leq \max \{ \kappa ,  \lambda  \}\,]$\\ 
\underline{$\mathrm{Case \,3_{}}$}: $[\, \kappa < \omega,\, \lambda \geq \omega \,]$ .  Then  $[\, | \lambda '| =\lambda\,]\,(8'')$ . If $\xi\in \lambda'$ then 
 $[\, | \kappa^{\, \xi}| < \omega\,]\,(9'')$ , since $[\,\kappa< \omega,\, \xi< \omega  \,]$. \\
So $(7)\Rightarrow$ 
 $[\,\, | R^{[\lambda]} | \leq \underset{\xi \in \lambda'} \Sigma \kappa^{\, \xi} \underset{(9'')}\leq |\lambda' | \cdot \omega \underset{(8'')} \leq  \lambda \cdot  \omega \underset{(2)}=\lambda=  \max \{ \kappa ,  \lambda  \}\,]\Rightarrow$ 
 $[\, | R^{[\lambda]} | \leq \max \{ \kappa ,  \lambda  \}\,]$\\ 
 So we proved that  $[\, | R^{[\lambda]} | \leq \max \{ \kappa ,  \lambda  \}\,]$ (11), in any case. Hence $\{(5),\,(11)  \}\Rightarrow$
 $[\, | R^{[\lambda]} |=\max \{ \kappa ,  \lambda  \}\,]$.     
\end{proof}

\begin{lemma}\label{shortEK}
Let $V$ be a infinite-dimensional vector space over a  field $F$. Then 
$\boxed{\big| V \big|= \max \{ \big|F\big|, \dim _F (V) \}}$,
 i.e $[\, \, \big| F^{[\lambda]}\big|= \max \{ \big| F\big|, \lambda \}$, if $\lambda\geq \omega \,]$.
\end{lemma} 

\begin{proof} By assumption $V$ has an infinite basis. Hence  $V\cong_F F^{[\lambda]}$ for  $[\,\omega\leq \lambda=\dim_F(V)\,]$ (2). Also  $[\,\big|F \big|\ge 2\,] $ (3) by convention. So 
  $ \big| V\big|=\big| F^{[\lambda]}\big| $
$\underset {(2),(3)}{\overset{\ref{cardinallemma1}} \Longrightarrow} \big| V\big|=\max\{\big| F\big| ,\lambda \} \Rightarrow$
$ \big| V\big|$
$=\max\{\big| F\big| ,\big| \dim_F (V)\big| \}$.
\end{proof} 

Note that for  the proof of  Theorem \ref{EKconsequence} ahead  we do not  actually need  the full \emph{Erdős  Kaplansky Theorem} but the lemma below.
\begin{lemma}\label{finiteEK}
Let $V$ be a infinite-dimensional vector space over a  finite field $F$. Then  
$\big| V \big|=  \dim _F (V)$.  \\ 
 So if   $\lambda$ is an infinite cardinal then  \space  $\dim_F(F^{\lambda})=\big| F^{\lambda} \big|=\big| F\big|^{\lambda}$. Hence
$\dim_{\mathbb{F}_2} (\mathbb{F}_2)^{\lambda}  =2^{\lambda}$
\end{lemma} 
$\mathbf{Proof:}$ We have that  $[\, \big| F\big|< \omega$, $\dim_F (V)\ge  \omega \,]$ (3) by assumption. Also  it follows  from \ref{shortEK} above  that
 $[\, \big| V\big|$
$=\max\{\big| F\big| ,\dim_F (V) \}\overset{(3)}=\dim_F (V)\,] \Rightarrow [\, \big|V \big|=\dim_F (V) \,] $.

\begin{lemma}\label{EKconsequence}
Let $V$ be an infinite-dimensional vector space over a  finite field $F$ and $\kappa,\, \lambda$ be cardinals. \\ 
$\mathrm{(i)}\, | V |=  \dim _F (V)$ \space 
$\,\mathrm{(ii)}\, \dim_{\mathbb{F}_2} (\mathbb{F}_2)^{\, \lambda}  =2^{\, \lambda},$ when  $\lambda \geq \omega$. \space 
$\mathrm{(iii)}\,(\mathbb{F}_2)^{\kappa}\cong_{\mathbb{F}_2}(\mathbb{F}_2)^{\, \lambda}\,]\,(1)\, \Leftrightarrow  \,[\, 2^ {\, \kappa} =2^{\lambda}\,]\,(2)\,. $ 

\end{lemma} 
\begin{proof} $\ref{EKconsequence}\mathrm{.i}$ follows from \ref{shortEK} (since $|F|\leq \dim_F(V)$) and $\ref{EKconsequence}\mathrm{.ii}$ is a special case of $\ref{EKconsequence}\mathrm{.i}$ for $V=(\mathbb{F}_2)^{\, \lambda}$. 
\underline{Proof of $\mathrm{(iii)}$} :
 $ [(1)\Rightarrow (2)] :$ \quad  $[\, (\mathbb{F}_2)^{\kappa}\cong_{\mathbb{F}_2}(\mathbb{F}_2)^{\lambda}\,]\Rightarrow$
$[\, \big|(\mathbb{F}_2)^{\kappa}\big|= \big|(\mathbb{F}_2)^{\lambda}\big|\,] \Rightarrow$
$[\, 2^{\,\kappa}=2^{\,\lambda}\,]$. \\
 $[ (2)\Rightarrow (1)] :$ \quad 
$\cdot$ Let $[\, (\kappa< \omega)\lor (\lambda< \omega) \,]$ (3).
Then $\{(1),\,(3),\, \ref{cardinallemma}\mathrm{.vi}  \}\Rightarrow$
$[\,\kappa=\lambda  \,]\Rightarrow$
 $[\, (\mathbb{F}_2)^{\kappa}\cong_{\mathbb{F}_2}(\mathbb{F}_2)^{\lambda}\,]$. \\
 $\cdot$ Let $[\, \kappa\geq \omega,\, \lambda\geq \omega \,]$ (4).
 Then $[\, \dim  (\mathbb{F}_2)^{\kappa}\overset{\mathrm{(ii)}}=2^{\, \kappa}\overset{(2)}= 2^{\, \mu} $
 $\overset{\mathrm{(ii)}}=\dim  (\mathbb{F}_2)^{\, \lambda}\,]\Rightarrow$
 $[\, (\mathbb{F}_2)^{\kappa}\cong_{\mathbb{F}_2}(\mathbb{F}_2)^{\lambda}\,]$.  \end{proof}

 \begin{lemma} \label{vandermonde}
 Let $F$ be any field and $\lambda$ be a cardinal. If  $\lambda\geq \omega$ then  $ \dim_F (F^{\lambda})\ge \big| F \big|$.
\end{lemma}
\begin{proof} For given $a\in F$,  consider the element $\theta_a \in F^{\mathbb{N}}$
defined by $\theta_a (n)=a^n $ i.e $\theta_a=(1,a,a^2,a^3,...)$ 
(here $0^0\overset{\text{def}}=1$). Let $L=\{\theta_a \,|\, a\in F \}$ \space\space\space  \\
\underline{$\mathrm{Claim \, 1_{}}$}:  $L$ is an $F$-linearly independent subset of $F^{\mathbb{N}}$. \\
$\mathrm{Proof:}$ Let $[\, a_1,a_2 \dots,a_n $ be distinct$\,]$ (1) such that  $[\, \overset{n} {\underset{k=1} \Sigma} \lambda_k \cdot \theta_{a_k}=0\,]$, for some
 $\vec{\lambda}=(\lambda_1,\lambda_2,..,\lambda_n)^{\top}\in F^{\, n\times 1}$  
   $\Rightarrow$ 
 $[\, \overset{n} {\underset{k=1} \Sigma} \lambda_k \cdot(1,a^1_k,a^2_k,...,a^{n-1}_k)=0\,]\Rightarrow [\, B\cdot \vec{\lambda}=0 \,]$ (2),  where  $B\in F^{n\times n}$  has $(i,j)$ entry   $B_{i,j}=a^{\, j-1}_i$ i.e it is a Vandermonde matrix. It is known that
 $\det(B)=\underset{i>j}\prod (a_i-a_j)\underset{\,(1)} \neq 0 \Rightarrow$  
  $[\, \det(B) \neq 0\,]$ (3). Therefore follows from 
 $\{ (2),(3)  \} $ that  $[\,  \vec{\lambda}=0\,]$. Therefore $L$ is an $F$-linearly independent subset of $F^{\mathbb{N}}$. \\
  Also the map 
 $\theta:F\rightarrow F^{\mathbb{N}}$ defined by $\theta(a)=\theta_a$ is trivially 1-1 hence $\big| L \big| = \big| F \big|$.  Since  linearly independent  subsets  extend to  bases, it follows that 
 $[\,  \dim_F (F^{\mathbb{N}})\ge \big| L \big|=\big| F \big| \,] \Rightarrow$
 $[\,\dim_F (F^{\mathbb{N}})\ge \big| F \big| \,]$. So 
 $[\, \dim_F(F^{\omega})\geq  \big| F \big|\,]$ (4). We have $[\lambda\geq \omega\,]\Rightarrow[\,  F^{\omega}\underset{F}  \rightarrowtail F^{\lambda}\,]\Rightarrow$
 $[\, \dim_F (F^{\lambda})\geq  \dim_F  (F^{\omega})\underset{(4)} \geq \big|F \big|\,]\Rightarrow$
  $[\,  \dim_F (F^{\lambda})\ge \big| F \big| \,]$.
 \end{proof}
 \begin{remark} \label{vandermonde2} We have actually shown  in \ref{vandermonde} that if $\lambda \geq \omega$ then $F^{\lambda}$ has an $F$-linearly independent set of cardinality  $\big| F \big|$, even when $F$ is a commutative domain. 
 \end{remark}

\begin{lemma} \label{cardinallemma}
Let  $\beta,\, \gamma,\, \kappa,\lambda, \mu$ be cardinals such  $[\,\lambda\geq \omega\,]\,(1)$. Then $:$ \\ 
$\mathrm{(i)}\,\,\lambda^{\, \lambda}=2^{\, \lambda}$\quad 
$\mathrm{(ii)}\,\, [\,2 \le \kappa \le\lambda\,] \Rightarrow [\, \kappa^{\, \lambda}=2^{\, \lambda}\,]$ $\mathrm{(iii)}\,\,[\, 2\le  \kappa \le \mathfrak{c} \,]\Rightarrow[\, \kappa^{\, \lambda}=2^{\, \lambda}\,]$ 
 $\mathrm{(iv)}\,\, \kappa^{ \, \lambda} \leq \mathrm{max} \{  2^{\,\lambda},\,  2^{\,\kappa} \} $  \
  $\mathrm{(v)}\, (\kappa<\omega) \Leftrightarrow (2^{\, \kappa}<\omega)\,]$ \quad\quad\quad\quad  $\mathrm{(vi)}$ If  $[\, (\kappa< \omega)\lor (\mu< \omega)\,]$ then   
  $[\,(2^{\, \kappa}=2^{\, \mu}) \Rightarrow (\kappa=\mu) \,]$\\
  $\mathrm{(vii)}$ $[\,\beta^{\,\lambda}=\gamma^{\,\lambda},\, \kappa \geq \lambda\geq \omega \,]\Rightarrow[\, \beta^{\,\kappa}=\gamma^{\,  \kappa}\,]$\\
 $\mathrm{(viii)}\,\, \mathrm{(a)}\, \, \lambda^{\bullet} >\lambda _{\cdot n} $. \quad\quad\ $\mathrm{(b)}\,\, \lambda _{\cdot n} \geq \omega$.  \quad\quad
$ \mathrm{(c)}\,\, (\lambda_{\cdot n})^{\, \omega}=\lambda_{\cdot n}$ , when  $n\geq 1$. \quad\quad
$ \mathrm{(d)} \,\,(\lambda^{\bullet})^{\, \omega} > \lambda^{\bullet}$. \end{lemma}
\begin{proof}  Note that  $\lambda_{\cdot n}=\beth_{\,n}(\lambda)$ and $\lambda^{\bullet}=\beth_{\,\omega}(\lambda)$, where defined in $\ref{sets}$ and that  $\{\mathrm{v},\,\mathrm{vi} \}$ are trivial. \\
 \underline{$\mathrm{Proof:}$ of $\mathrm{(i)}$} : \space  
$[\, 2\le \lambda \le 2^{\, \lambda}\,]\Rightarrow$
$[\,2^{\, \lambda}\le \lambda^{\lambda}\le (2^{\, \lambda})^{\lambda}=2^{\, {\lambda}\cdot \lambda}\overset{(1)}=2^{\,\lambda}\,]\Rightarrow$
$[\,2^{\, \lambda}\le \lambda^{\, \lambda} \le 2^{\, \lambda}\,]\Rightarrow$
$[\,\lambda^{\,\lambda}=2^{\,\lambda}\,]$\\
\underline{$\mathrm{Proof:}$ of $\mathrm{(ii)}$} : \space
 $[\,  2 \le \kappa \le\lambda\,] \Rightarrow [\, 2^{\lambda}\le \kappa^{\, \lambda}\le \lambda^{\, \lambda}\overset{\mathrm{(i)}}=2^{\,\lambda}\,]\Rightarrow$
$[\, 2^{\,\lambda}\le \kappa^{\, \lambda}\le 2^{\, \lambda}\,]\Rightarrow$
$[\, \kappa^{\, \lambda}= 2^{\,\lambda}\,]$ \\
\underline{$\mathrm{Proof:}$ of $\mathrm{(iii)}$} : \space
 We consider the following four cases :\\ 
$\mathrm{(a)}\,\,\kappa=\omega\, :$ Then $[\, 2\leq \omega \le \lambda\,]\overset{\mathrm{(ii)}}\Rightarrow [\, \omega^{\,\lambda}=2^{\,\lambda}  \,]\Rightarrow[\, \kappa^{\, \lambda}=2^{\, \lambda}\,]$.\\
$\mathrm{(b)\,\,}\kappa<\omega\, :$. Then  $[\, 2\le \kappa \le \omega\overset{\mathrm{(ii)}}\,]\Rightarrow$
$[\, \kappa^{\, \lambda}={\omega}^{\, \lambda  \,}\overset{\mathrm{(iii.a)}}=\, 2^{\, \lambda}\,]\Rightarrow[\, \kappa^{\, \lambda}=2^{\,\lambda}\,]$.\\
 $\mathrm{(c)}\,\, \kappa=\mathfrak{c}\,:$ Then $[\,\kappa=2^{\,\omega}\,]\Rightarrow$  $[\, \kappa^{\, \lambda}= (2^{\, \omega})^{\lambda}=2^{\, \omega\cdot \lambda}=2^{\, \max\{\omega,\lambda \}}=2^{\, \lambda}\,]\Rightarrow$ 
$[\, \kappa^{\, \lambda}=2^{\, \lambda}\,]$. \\
$\mathrm{(d)}\,\,\omega \le\kappa \le\mathfrak{c}\,:$ Then 
$[\, 2^{\, \lambda}\overset{\mathrm{(iii.a)}}=\omega^{\, \lambda}\le \kappa^{\lambda}\le \mathfrak{c}^{\, \lambda}\overset{\mathrm{(iii.c)}}=2^{\, \lambda}\,] $
$\Rightarrow[\,  \kappa^{\, \lambda}=2^{\, \lambda}\,]$.\\
\underline{$\mathrm{Proof:}$ of $\mathrm{(iv)}$} : \space
 The claim is trivial for $ \kappa \leq 1$ (Note that 
$0^0=1$). Assume that 
$[\, \kappa\geq 2\,]$.\\
$\cdot$ Let   $[\,2\leq \kappa \leq \lambda \,]\overset{\mathrm{(ii)}} \Rightarrow$ $[\,  \kappa^{\,\lambda}=2^{\,\lambda}\,]$. In this case 
$\mathrm{max} \{  2^{\,\lambda},\,  2^{\,\kappa} \}=2^{\, \lambda} $. Therefore
 $[\,  \kappa^ { \,\lambda} \leq \mathrm{max} \{  2^{\, \lambda},\,  2^{\, \kappa} \} \,]$.   \\ 
$\cdot$ Let  $[\, \kappa > \lambda \,] \Rightarrow$
 $[\,  \kappa^{\, \lambda} \leq \kappa^{\,\kappa}\overset{\mathrm{(ii)}}=2^{\,\kappa}\,]$,
 since $[\, \kappa>\lambda>\omega\,]$. Here
$\mathrm{max} \{  2^{\lambda},\,  2^{\kappa} \}=2^{\lambda}$. Hence  
 $[\,  \kappa^  {\, \lambda} \leq \mathrm{max} \{  2^{\, \lambda},\,  2^{\, \kappa} \}\,] $.
 So we proved the claim for any case. \\
\underline{$\mathrm{Proof:}$ of $\mathrm{(vii)}$} : By assumption
$[\, \omega \leq \lambda \leq \kappa\,]\Rightarrow$ 
 $[\, \lambda \cdot \kappa=\mathrm{max}\{\kappa,\,\lambda\}=\kappa \,]\, (5)$. Also by assumption 
 $[\, \beta^{\,\lambda}=\gamma^{\, \lambda}\,]\,(6)$. Hence  
 $[\,\beta^{\, \kappa}\overset{(5)}=\beta^{\,  \lambda \cdot \kappa}=(\beta^{\, \lambda})^{\,\kappa}\overset{(6)}=(\gamma^{\,\lambda})^{\,\kappa}=\gamma^{\, \lambda \cdot \kappa}\overset{(5)}=\gamma^{\,\kappa}\,]\Rightarrow$ 
 $[\, \beta^{\,\kappa}=\gamma^{\,  \kappa}\,]$.
\underline{Proof of $\mathrm{(viii.a)}$} :
$[\, \lambda^{\bullet} \overset{2.\mathrm{ii}}=\underset{i\in \omega} \Sigma \lambda_{\cdot i} \geq \lambda_{\cdot (n+1)}=2^{\, \lambda_{\cdot n}}> \lambda_{\cdot n}\,]\Rightarrow$ 
$[\,\lambda^{\bullet} >\lambda _{\cdot n}  \,] $. \\ 
\underline{Proof of $\mathrm{(viii.b)}$} :
 Follows  trivially by induction. Note that  $\lambda_{\cdot n}=\beth_{\,n}(\lambda)$ and $\lambda^{\bullet}=\beth_{\,\omega}(\lambda)$.   \\
\underline{Proof of $\mathrm{(viii.c)}$} :
$[\, n\geq 1\,]\Rightarrow[\, \exists\,  m\in \omega : (n=m+1)\,]\Rightarrow$
$[\,(\lambda_{\cdot n})^{\, \omega}=(2^{\,\lambda _{\cdot m} })^{\, \omega}=2^{\, \omega \cdot \lambda _{\cdot m} }=2^{\, \max \{ \omega ,\,      \lambda _{\cdot m}\}  } $
$\overset{\mathrm{(b)}}=2^{\, \lambda _{\cdot m} }=\lambda _{\cdot n} \,]\Rightarrow$
$[\, (\lambda_{\cdot n})^{\, \omega}=\lambda_{\cdot n} \,]$.\\
\underline{Proof of $\mathrm{(viii.d)}$} :
$[\, (\lambda^{\bullet})^{\, \omega}=\underset{i \in \omega}\Pi \lambda \underset{\mathrm{(a)}} > \underset{i \in \omega} \Sigma \,  \lambda_{\cdot i}= \lambda^{\bullet}\,]\Rightarrow$
$[\,(\lambda^{\bullet})^{\, \omega} > \lambda^{\bullet} \,]$ (\emph{Kőnig's inequality} was used). \end{proof}

\begin{theorem} \label{longproof}
$^{\ast}$Let $X,Y$ be two  sets such that   $\big| \mathcal{P}(X)\big| =\big|\mathcal{P}(Y)\big|$. There is a function  
$\hat{f}:\mathcal{P}(X)\rightarrow \mathcal{P}(Y)$  such that
 $ \mathrm{(i)}\,  \hat{f}(\emptyset)=\emptyset$. \quad    $\mathrm{(ii)} \,\hat{f}$ is a bijection. \quad
$\mathrm{(iii)} \,[\, \hat{f}(A\oplus B)=\hat{f}(A)\oplus\hat{f}(B)\,]$,
when $\{A,B\}\subseteq \mathcal{P}(X)$.
\end{theorem}
\begin{proof} Let  $X'=(\mathbb{F}_2)^X$ and  $Y'=(\mathbb{F}_2)^Y$. 
By assumption  $\big| \mathcal{P}(X)\big| =\big|\mathcal{P}(Y)\big|\overset{\ref{EKconsequence}\mathrm{.iii}}\Rightarrow$   
 $X' \underset{\mathbb{F}_2}\cong Y'  \Rightarrow $. There is a map  $f:X'\rightarrow Y'$  such that 
 $[\, f$ is $\mathbb{F}_2$-linear $]\,(1)$ and $[\, f$ is bijection$\,]\,\mathrm{(ii)}$. It is trivial that $[\,\Phi_X,\, (\Phi_Y)^{-1}$
 are bijections$\,]\,(2)$, where     $\Phi_X:\mathcal{P}(X)\rightarrow X'$ and $(\Phi_Y)^{-1}: Y'\rightarrow \mathcal{P}(Y)$  were defined in     $\mathrm{2.i}$.  
Let $\hat{f}:\mathcal{P}(X)\rightarrow \mathcal{P}(Y)  $  defined by  $[\,\hat{f}=(\Phi_Y)^{-1}\circ f\circ\Phi_X\,]\,(3)$. It follows from
$\{ (1), (3) \}$ that 
$[\,  f(\emptyset)= \emptyset\,]\, \mathrm{(i)}$ and from
$\{ (2), (3) \}$ that  $[\, \hat{f}$ is a bijection$\,]\mathrm{(ii)}$. 
Note  
$[\, \Phi_Z(A\oplus B)=\Phi_Z(A)+\Phi_Z(B) \,]\,(4)$ and   $[\, (\Phi_Z)^{-1}(u+v)=(\Phi_Z)^{-1}(u) \oplus(\Phi_Z)^{-1}(v)\,]\,(5)$.    Finally $\mathrm{(iii)}$  follows from  $\{ (1),(3),(4),(5) \}$. \end{proof}

\begin{theorem} \label{EKextension}
Let $V$ be an infinite-dimensional vector space over a field $F$ 
$($i.e $\dim_F (V \geq \omega)$. Then $:$
 \end{theorem} 
 \begin{enumerate}[label=(\roman*)]
\item  $\boxed{| V^{\ast} |=\dim_F (V^{\ast}) }$ , i.e $[\,\dim_F (F^{\lambda})=| F^{\lambda} |$, if $\lambda \geq \omega \,]$ (\emph{Erdős  Kaplansky Theorem}).\\
$\mathbf{Proof:}$  By assumption  $V\underset{F}\cong F^{[\lambda]}$ for  $[\,\omega\leq \lambda=\dim_F(V)\,]$. So
$V^{\ast}\underset{F} \cong F^{\lambda}$ and $[\, \dim_F(V^{\ast})\geq \dim_F(V)\geq \omega \,]$ (1).
Therefore   $| V^{\ast} |= | F^{\lambda} |   \underset{(1)}{\overset{\ref{shortEK}} =}$
$\mathrm{max}\{ \dim_F (F^{\lambda}),\, | F |  \}$
$\overset{\ref{vandermonde}}=\dim_F (F^{\lambda})=
\dim_F(V^{\ast})$.\\\\
$\bullet$ Let now  $[\, \lambda=\mathrm{dim}_F(V) \geq \omega \,]$ (2) and  $[\, \kappa=| F |\geq 2\,]$ (1),
for the proofs of $\mathrm{(ii)}$,$\mathrm{(iii)}$,$\mathrm{(iv)}$,
$\mathrm{(v)}$,$\mathrm{(vi)}$.
\item  $\mathrm{(a)}\,\,\dim_F (V^{\ast}) \ge 2^{\, \dim_F (V)}$. \quad\quad\quad\quad\quad\quad\quad\quad  $\mathrm{(b)}\,\,\dim_F (V^{\ast}) > \dim_F (V)$.
 \\
$\mathbf{Proof:}$   
  $[\,\dim_F (V^{\ast})\overset{\ref{EKextension} \mathrm{.i}}=| V^{\ast} | = $
  $| F^{\lambda} |=| F |^{\,\lambda}=\kappa^{\,\lambda} \,]\Rightarrow $  $[\,\dim_F (V^{\ast})=\kappa^{\,\lambda} \,]$ (3)
  $\overset{(2)}\Rightarrow$ 
 $[\,\dim_F (V^{\ast})\geq 2^{\,\lambda}=2^{\, \dim_F (V)}\,]\Rightarrow$               
 $[\, \dim_F (V^{\ast}) \ge 2^{\, \dim_F (V)}\,]\, \mathrm{(a)}\Rightarrow$ 
$[\, \dim_F (V^{\ast}) \ge 2^{\, \dim_F (V)} >  \dim_F (V)\,]\Rightarrow$  
$[\, \dim_F (V^{\ast})  >  \dim_F (V)\,]\, \mathrm{(b)}$. 
     
\item $\dim_F (V^{\ast}) = 2^{\, \dim_F (V)}$, when $[\,\dim_F (V) \ge | F |\,]$ (4).  \\
$\mathbf{Proof:}$ $\{(1),(2),(4) \} \Rightarrow$   
 $[\,  2 \le \kappa  \le \lambda\,] \overset{\ref{cardinallemma}\mathrm{.ii}} \Rightarrow [ \, \kappa^{\,\lambda}=2^{\, \lambda} \,]$ 
 $\overset{(3),(2))}\Longrightarrow$
 $[\,\dim_F (V^{\ast}) = 2^{\, \dim_F (V)}\,]$.
 \item $^{\ast}$ $\boxed{\dim_F (V^{\ast}) = 2^{\, \dim_F (V)},\, \text{when}\, | F |  \le  \mathfrak{c}}  $  \\
$\mathbf{Proof:}$   By assumptions $[\, 2\le  \kappa=| F |\le \mathfrak{c}\,]\overset{\ref{cardinallemma}\mathrm{.iii}}\Rightarrow$
$[\, \kappa^{\, \lambda}=2^{\, \lambda}\,]$ $\overset{(3),(2))}\Longrightarrow$
 $[\,\dim_F (V^{\ast}) = 2^{\, \dim_F (V)}\,]$.

 \item  $\dim_F (V^{\ast}) > 2^{\, \dim_F (V)}$, when $[\, \dim_F (V) =\omega, \, | F | > \mathfrak{c}  \,\,]$. \\
$\mathbf{Proof:}$ By assumption we have $[\, \lambda=\omega\,]$ (5)
and $[\,\kappa> \mathfrak{c}\,]$ (6). Hence  
$[\, \dim_F (V^{\ast})\overset{(3)}=$
  $\kappa^{\lambda}\overset{(5)}=\kappa^{\,\omega }$
$\underset{(1)} \ge \kappa \underset{(6)}>\mathfrak{c}=2^{\,\omega }\underset{(5),(2)}=2^{\, \dim_F (V)}\,] \Rightarrow$
$[\, \dim_F (V^{\ast})> 2^{\, \dim_F (V)}\,]$. 

 \item   $ V{\not\cong}_ F V^{\ast\ast}$,   hence $V$ is a
 non-reflexive $F$-module. \\ 
 $\mathbf{Proof:}$ Note that $[\, V^{\ast \ast}\underset F{\cong} (V^{\ast})^{\ast}\,]$ (7). Also 
  $[\, \mathrm{dim}_{F} (V^{\ast}) \geq \mathrm{dim}_{F} (V) \overset{(2)}\geq \omega\,]\Rightarrow$ 
 $[\, \mathrm{dim}_{F} (V^{\ast}) \geq \omega\,]\overset{\mathrm{(ii.b)}}\Rightarrow$
 $[\, \mathrm{dim}_{F} (V^{\ast})^{\ast} > \mathrm{dim}_{F} (V^{\ast})$
 $\geq \mathrm{dim}_F(V)\,]\overset{(7)}\Rightarrow$ 
 $[\, \mathrm{dim}_{F} (V^{\ast\ast}) >  \mathrm{dim}_F(V)\,]\Rightarrow$   $[\, V{\not\cong}_F V^{\ast\ast}\,]$. 
\end{enumerate} 

\begin{remark}  \label{erdosremark}$:$ \\
$\mathrm{(i)}\, $   The full   \emph{Erdős-Kaplansky Theorem} is 
$\ref{EKextension}\mathrm{.i}$ and  has been largely ignored from textbook literature. Some authors refer to it's weaker version $\ref{EKextension}\mathrm{.ii.b}$ instead in order to prove $\ref{EKextension}\mathrm{.vi}$.
This  may be explain why nobody has yet formulated 
$\ref{EKextension}\mathrm{.iv}$, which we claim as a new result.\\    $\mathrm{(ii)}\, $  Note that  $\ref{EKextension}\mathrm{.i}$ is proven  in \cite{jacob}  for a division  ring $F$. More accurately   $[\, | V^{\ast} |=\dim_{\hat{F}} (V^{\ast})$, when   $\dim_F(V) \geq \omega \,]\, (\ref{EKextension}\mathrm{.i}\, ') $ is proven. It follows that \underline{the whole Theorem \ref{EKextension} is valid for a division ring $F$.} (where $\dim_F (V^{\ast})$ is replaced by $\dim_{\hat{F}} (V^{\ast})$. The same holds for Lemma $\ref{isomorphicduals}$ below. \\ $\mathrm{(iii)}\, $ It follows from 
$(\ref{EKextension}\mathrm{.i}')$ that $\ref{vandermonde}$ holds for a division ring $F$ as well since  $ \dim_F (F^{\kappa})= | F |^{\kappa}\ge | F |$. We do not know whether the  \emph{Vandermonde matrix} argument  that we used in the proof of  $\ref{vandermonde}$ can be  modified to work  for a division ring $F$, since determinants  are not applicable to non-commutative rings.\\
 $\mathrm{(iv)}\, $  In \cite{jacob}, Jacobson attributes \ref{vandermonde} to Mackey but does not attribute $\ref{EKextension}\mathrm{.i}$ to anyone. We can't explain why  $\ref{EKextension}\mathrm{.i}$ was named after and Erdős and Kaplansky, since it's earliest reference  is   \cite{jacob} in 1953.  \\
 $\mathrm{(v)}\, $ It follows from the \emph{Löwenheim–Skolem theorem}  that for any infinite cardinal $\kappa$ there is a field $F$ such that $|F|=\kappa$. 
\end{remark}
\begin{theorem} \label{vandermonde3} Let $F$ be a field $($or division ring$)$ and $\lambda \in \mathrm{Card}_{\geq \omega}$. Then $:$ \\
$\mathrm{(i)}$ $ \dim_F (F^{\omega})\ge | F |$. \space
$\mathrm{(ii)}$ $ \dim_F (F^{\omega})= | F \big|$, if $[\, |F|=\lambda_{\cdot n}, \, n\geq 1 \,]$. \space
$\mathrm{(iii)}$ $ \dim_F (F^{\omega}) > |F|$, if $ |F|=\lambda^{\bullet}$. 
\end{theorem}
\begin{proof}
$\ref{vandermonde3}\mathrm{.i}$ was proven in  $\ref{erdosremark}\mathrm{.iii}$,
$\ref{vandermonde3}\mathrm{.ii}$ follows from $\{ \ref{EKextension}\mathrm{.i}, \, \ref{cardinallemma}\mathrm{.vi.c} \}$ and
$\ref{vandermonde3}\mathrm{.iii}$  from $\{ \ref{EKextension}\mathrm{.i}, \, \ref{cardinallemma}\mathrm{.vi.d} \}$.\end{proof}

\subsection{Vector spaces with isomorphic duals are isomorphic ?}
\begin{definition} \label{definitions}
Let $\kappa\in \mathrm{Card}$  and $F$ be a field (or divison ring).
 Consider the statements below :
\end{definition}
\begin{enumerate}[label=(\roman*)]
\item $\mathsf{ICF}_{\kappa}$   
  $\Leftrightarrow[\, \forall \mu \in \mathrm{Card}\,]  [\, \forall \lambda \in \mathrm{Card}\,]$ 
     $[\,(\kappa^{\, \mu}=\kappa^{\, \lambda}) \Rightarrow (\mu={\, \lambda}) \,]$.
\item $\mathsf{ICF}_{\kappa}^{\,'}$  $\Leftrightarrow[\, \forall \mu \in \mathrm{Card}_{\geq \omega} \,]$
$[\, \forall \lambda \in \mathrm{Card}_{\geq \omega} \,]$
  $[\,(\kappa^{\,\mu}=\kappa^{\, \lambda}) \Rightarrow (\mu=\lambda) \,]$. 
\item $\mathsf{ICF}\Leftrightarrow\mathsf{ICF}_{2}$  $\Leftrightarrow[\, \forall \mu \in \mathrm{Card}\,]$
$[\, \forall \lambda \in \mathrm{Card}\,]$
  $[\,(2^{\, \mu}=2^{\, \lambda}) \Rightarrow (\mu=\lambda) \,]$.  
\item $\mathsf{CH}_{\kappa}$  $\Leftrightarrow[\, \forall \lambda \in \mathrm{Card}\,]$
  $\{ \,[\, \kappa \leq \lambda \leq 2^{\, \kappa}\,]   \Rightarrow [\, 
   (\lambda=\kappa) \lor (\lambda=2^{\, \kappa}) \,]\,\}$.
  \item $\mathsf{CH}\Leftrightarrow\mathsf{CH}_{\omega}$
   $\Leftrightarrow[\, \forall \lambda \in \mathrm{Card}\,]$
  $\{ \,[\, \omega \leq \lambda \leq 2^{\, \omega}\,]   \Rightarrow [\, 
   (\lambda=\omega) \lor (\lambda=2^{\, \omega}) \,]\,\}$.
\item $\mathsf{GCH}\Leftrightarrow [\,\mathsf{CH}_{\lambda}\,]\, \forall \, \lambda \in \mathrm{Card}_{\geq \omega}$ .  Note that  $\mathsf{GCH}$ is called  as the \emph{Generalized Continuum Hypothesis}.

  \item  $\mathsf{L}_F \Leftrightarrow$
 $[\, \forall\, V\in F$-$\mathrm{Mod}\,][\,\forall\, W\in F$-$\mathrm{Mod}\,]$
 $[\, ( V^{\ast} \cong_{\hat{F}} W^{\ast}) \Rightarrow ( V \cong_F W) \,]$. Note  $\hat{F}=F$, if $F$ is a field.  
  
\end{enumerate} 

\begin{lemma} \label{isomorphicduals}
Let $V,W$ be $F$-vector spaces over a field $($or division ring$)$ $F$ such that $\big| F \big|=\kappa$.
\end{lemma}
\begin{enumerate}[label=(\roman*)]
\item  Trivially $:$  $\mathrm{(a)}$ $[\,\mathrm{dim}_F(W^{\ast}) \geq \mathrm{dim}_F(W)] $. \quad\quad
 $\mathrm{(b)}$ $[\, \mathrm{dim}_F(V)<\omega\,] \Rightarrow [\,\mathrm{dim}_F(V^{\ast})=\mathrm{dim}_F(V)<\omega\,]$.  
\item  $[\, \mathrm{dim}_F(W^{\ast})<\omega\,]$ (2) $\Rightarrow$
$[\,\mathrm{dim}_F(W^{\ast})=\mathrm{dim}_F(W)<\omega\,] $.    \\
$\mathbf{Proof:}$ Assume that  $[\,\mathrm{dim}_F(W)\geq \omega\,]$ (3). Then $[\,\mathrm{dim}_F(W^{\ast}) \overset{(\mathrm{i.a})}\geq \mathrm{dim}_F(W)\overset{(3)}\geq \omega\,]\Rightarrow$
$[\,\mathrm{dim}_F(W^{\ast}) \geq \omega\,]$, which contradicts (2). Hence
 $[\,\mathrm{dim}_F(W)< \omega\,]\overset{(\mathrm{i.b})}\Rightarrow$
 $[\,\mathrm{dim}_F(W^{\ast})=\mathrm{dim}_F(W)<\omega\,]$. 
\item Assume that  $[\, \mathrm{dim}_F(V)<\omega\,]\lor[\, \mathrm{dim}_F(W)<\omega\,]$. Then 
 $[\, V^{\ast} \cong_F W^{\ast} \,]\Rightarrow [\, V \cong_F W \,]$. \\
$\mathbf{Proof:}$ $\cdot$ Let  $[\, \mathrm{dim}_F(V)<\omega\,]\overset{(\mathrm{i.b)}}\Rightarrow$
$[\,\mathrm{dim}_F(V^{\ast})=\mathrm{dim}_F(V)<\omega\,]$ (5).
 By assumption 
 $[\, V^{\ast} \cong_F W^{\ast} \,]\Rightarrow $
 $[\, \dim_F (V^{\ast})= \dim _F (W^{\ast}) \,(6)$
 $\overset{(5)}\Rightarrow$ 
$[\,  \dim_F  (W^{\ast})< \omega  \,]\overset{(\mathrm{ii})} \Rightarrow$   $[\,\mathrm{dim}_F(W^{\ast})=\mathrm{dim}_F(W)<\omega\,] $ (7). So $\{(5),(6),(7) \} \Rightarrow$    $[\, \dim_F  (V)= \dim_F  (W) \,]\Rightarrow$  $[\, V \cong_F W \,]$. \\
$\cdot$ Let  $[\, \mathrm{dim}_F(W)<\omega\,]$. The proof is similar.
\item  Assume $\{\, \mathsf{ICF}_{\kappa}^{\,'}$,  $[\, \mathrm{dim}_F(V) \geq\omega\,],[\, \mathrm{dim}_F(W)\geq \omega\,]\,\}$. Then 
 $[\, (V^{\ast} \cong_F W^{\ast}) \,]\Rightarrow [\,( V \cong_F W) \,]$.  \\
$\mathbf{Proof:}$ Let $[\, \lambda=\mathrm{dim}_F(V)\geq \omega$, $\mu=\mathrm{dim}_F(W)\geq \omega\,]$ (8). It follows from $\{(8),\, \ref{EKextension}\mathrm{.i} \}$ that $[\, \dim_F (V^{\ast})=\kappa^{\, \lambda}\,]$ (9) and  that $[\, \dim_F (W^{\ast})=\kappa^{\, \mu}\,]$ (10).
 By assumption, we have 
 $[\, V^{\ast} \underset{F}\cong W^{\ast} \,]\Rightarrow $
 $[\, \dim_F  (V^{\ast})= \dim_F  (W^{\ast}) \,]$ (11). Therefore
 $\{(8),(9),(10),(11)  \}\Rightarrow$  $[\, \kappa^{\, \lambda}=\kappa^{\, \mu},\,\lambda \geq\omega,\, \mu\geq \omega\,]\underset{\mathsf{ICF}_{\kappa}^{\, '}}\Longrightarrow$
 $[\,\lambda=\mu \,]\underset{(8)}\Rightarrow$
  $[\, \dim_F  (V)= \dim_F  (W) \,]\Rightarrow$  $[\, V \cong_F W \,]$.
 \item $\mathsf{L}_F \Rightarrow  \mathsf{ICF}_{\kappa}^{\,'}$. \\
$\mathbf{Proof:}$
Let $\lambda,\, \mu$ be cardinals such that  $[\, \kappa^{\, \lambda}=\kappa^{\, \mu},\,\lambda \geq\omega,\, \mu \geq \omega\,]$ (12). We
consider the $F$-vector spaces
 $V=F^{[\lambda]},\, W=F^{[\mu]}$. So $V^{\ast}\cong_F F^{\lambda}$ and 
 $W^{\ast}\cong_F F^{\mu}$. Hence
 $[\, \tmop{dim}_F (V^{\ast})\overset{\ref{EKextension}\mathrm{.i}}=|V^{\ast}  |=|F^{\lambda}|=\kappa^{\, \lambda}\,]\Rightarrow$  
  $[\, \tmop{dim}_F (V^{\ast})=\kappa^{\, \lambda}\,]$ (13). Similarly 
 $[\, \tmop{dim}_F (W^{\ast})=\kappa^{\, \mu}\,]$ (14). Therefore 
 $\{(12),(13),(14) \}\Rightarrow$
 $[\, \dim_F  (V^{\ast})= \dim_F  (W^{\ast}) \,]\Rightarrow$  $[\, V^{\ast} \cong_F W^{\ast} \,]$ 
 $\underset{\mathsf{L}_F}\Longrightarrow$
 $[\, V \cong_F W \,]\Rightarrow$ 
  $[\, \dim_F  (V)= \dim_F  (W)\,] \Rightarrow [\,\lambda=\mu  \,]$.
  Hence $\{\, [\, \kappa^{\, \lambda}=\kappa^{\, \mu},\,\lambda \geq\omega,\, \mu \geq \omega\,]\Rightarrow  [\,\lambda=\mu  \,]\,\}$
$\underset{\ref{definitions}\mathrm{.ii}}\Longrightarrow$
$ \mathsf{ICF}_{\kappa}^{\,'}$.
\item  Assume that $\exists\, \lambda \in \mathrm{Card}:[\,\kappa=2^{\,\lambda},\, \lambda \geq \mathfrak{c}\,]$. Then : $\mathrm{(a)}$ $\mathsf{ICF}_{\kappa}^{\,'}$  fails. \quad \quad
$\mathrm{(b)}$ $\mathsf{L}_F$ fails. \\
$\mathbf{Proof:}$ By assumption $[\, \lambda \geq \mathfrak{c}>\omega  \,]$ (16). So
$[\, \kappa^{\,\omega}=(2^{\,\lambda})^{\,\omega}=2^{\, \lambda \cdot \omega}=2^{\, \mathrm{max} \{\lambda,\, \omega \}}\overset{(16)}=2^{\, \lambda}=\kappa\,]\Rightarrow [\, \kappa^{\,\omega}=\kappa\,]$ (17).
Similarly  $[\, \kappa^{\,\mathfrak{c}}=\kappa\,]$ (18). So
$\{(16),(17),(18), \ref{definitions}\mathrm{.ii} \} \Rightarrow$  
$[\, \mathsf{ICF}_{\kappa}^{\,'}$  fails$\,]\overset{(\mathrm{v})}\Rightarrow$ $[\, \mathsf{L}_F$ fails$\,]$.
\item $\mathrm{(a)} \,\,[\,  \kappa \leq \mathfrak{c}\, ,\, \mathsf{ICF} \,]\Rightarrow 
 \mathsf{ICF}_{\kappa}^{\,'}$. \quad \quad \quad \quad
  $\mathrm{(b)}\,\,  \mathsf{ICF}_{\kappa}^{\,'}\Rightarrow \mathsf{L}_F $.  \quad \quad \quad\quad\quad $\mathrm{(c)}$ $\mathsf{L}_F \Leftrightarrow \mathsf{ICF}_{\kappa}^{\,'}$.\\
$\mathbf{Proof:}$ $\mathrm{(a)}$ Follows from $\ref{cardinallemma}\mathrm{.iii}$.\space\space   $\mathrm{(b)}$ Follows from $\{\mathrm{iii}, \mathrm{iv}\}$. \quad\space  $\mathrm{(c)}$ Follows from $\{\mathrm{(v)}, \mathrm{(vii.b)}\}$.

\item $\mathrm{(a)} \,[\,  \omega\leq \lambda < \kappa\, ,\,\mathsf{CH}_{\kappa},\,  \mathsf{ICF}_{\kappa}^{'} \,]\Rightarrow$
$[\, \kappa^{\lambda}=\kappa\,]$. \quad \quad \quad\quad
$\mathrm{(b)} \,[\,  \kappa >\mathfrak{c}\, ,\,\mathsf{CH}_{\kappa} \,]\Rightarrow$ $[\, \mathsf{ICF}_{\kappa}^{\,'}$ fails$\,]$.\\
$\mathbf{Proof:}$ \\
\underline{Proof of $\mathrm{(a)}$} :         
$[\, \kappa\leq \kappa^{\, \lambda}\leq \kappa^{\, \kappa}\overset{\ref{cardinallemma}\mathrm{.i}}=2^{\, \kappa}\,]\Rightarrow$  
$[\,\kappa \leq \kappa^{\lambda}\leq 2^{\, \kappa}  \,]$
$\overset{\mathsf{CH}_{\kappa}}\Longrightarrow$
$[\,(  \kappa^{\, \lambda}=\kappa) \lor  (\kappa^{\, \lambda}=2^{\, \kappa)}\,]  $ (21). 
Assume that $[\,  \kappa^{\, \lambda}=2^{\, \kappa}\,]$
 then $[\, \kappa^{\, \lambda}=2^{\, \kappa}=\kappa^{\, \kappa}\,]\Rightarrow[\,\kappa^{\, \lambda}=\kappa^{\, \kappa},\,  \omega< \lambda < \kappa\,] $, which contradicts $\mathsf{ICF}_{\kappa}^{'}$. Hence 
 $[\,  \kappa^{\, \lambda} \neq 2^{\, \kappa}\,]\overset{(21)}\Rightarrow$
$[\,  \kappa^{\, \lambda}=\kappa\,]$. \\
\underline{Proof of $\mathrm{(b)}$} : Assume   $[\, \mathsf{ICF}_{\kappa}^{\,'}$ holds$\,]$ and  $[\,  \kappa >\mathfrak{c}\, ,\,\mathsf{CH}_{\kappa} \,]$  . Then 
$[\,  \omega\leq \omega < \kappa\, ,\, \omega\leq \mathfrak{c} < \kappa, \,
\mathsf{CH}_{\kappa}\, ,\,   \mathsf{ICF}_{\kappa}^{'} \,]$
$\overset{\mathrm{(a)}}\Rightarrow$
$[\,\kappa^{\, \omega}=\kappa ,\, \,   \kappa^{\, \mathfrak{c}}=\kappa  \,]\Rightarrow$
$[\, \kappa^{\, \omega}= \kappa^{\, \mathfrak{c}},\, \omega \neq \mathfrak{c}  \,]$, which contradicts $\mathsf{ICF}_{\kappa}^{\,'}$.
Therefore   $[\, \mathsf{ICF}_{\kappa}^{\,'}$ fails$\,]$.
  
\item  $\mathrm{(a)}\,\, \mathsf{GCH}\Rightarrow \mathsf{ICF}$. \quad
Let $\lambda \in \mathrm{Card}_{\geq \omega}$. Trivially : 
    $\mathrm{(b)}\,\,$ $\mathsf{GCH}\Rightarrow \mathsf{CH}_{\lambda}$ . \quad 
  $\mathrm{(c)}\,\,$  $ \mathsf{CH}_{\lambda} \Leftrightarrow
  [\, \lambda^{+}=2^{\, \lambda} \,]$. \\
\underline{Proof of $\mathrm{(a)}$} : Let $\mu, \, \lambda$ be cardinals such that 
 $[\,2^{\, \mu}=2^{\, \lambda}\, ]$. We will show that  $[\, \mu=\lambda\,]$. \\
 $\cdot$ Let $[\, \mu< \omega$ . Then  the claim  follows from 
$\ref{cardinallemma}\mathrm{.v}$, without assuming  $\mathsf{GCH}$.\\
 $\cdot$ Let $[\, \mu\geq \omega,\, \lambda< \mu\,]$. Then
 $[\, \lambda< \mu< 2^{\, \mu}=2^{\lambda}\,]\Rightarrow$
 $[\, \lambda< \mu< 2^{\,\lambda},\,  \mu\ge \omega]$, which contradicts
 $\mathsf{GCH}$. \\$\cdot$ Let $[\, \mu\geq \omega,\, \lambda > \mu\,]$. Then $[\, \mu< \lambda< 2^{\, \lambda}=2^{\, \mu}\,]\Rightarrow$
 $[\, \mu< \lambda< 2^{\, \mu},\,  \lambda > \omega]$, which contradicts
 $\mathsf{GCH}$. \\
 Therefore, it remains the case $[\, \mu\geq \omega,\, \lambda= \mu\,]$. Hence $[\,\mu=\lambda\,]$. 
\end{enumerate}

 \begin{theorem} \label{ICFdual}
 $^{\ast}$ Let $F$ be  a field $($or division ring$)$ such that $\big| F \big|\leq \mathfrak{c}$.  Then  \space $ \mathsf{ICF}\Leftrightarrow \mathsf{L}_F $.
  \end{theorem}
\begin{proof} Let $\kappa=|F|$. Then $[\,2\leq\kappa \leq \mathfrak{c} \,]\,(2)$, by assumption. \\
 $\cdot$ Proof of  $[\, \mathsf{ICF} \Rightarrow \mathsf{L}_F \,]$ :
 It follows from  $\{ \ref{isomorphicduals}\mathrm{.vii.a},\,\ref{isomorphicduals}\mathrm{.vii.b}\} $. \\
 $\cdot$ Proof of  $[\, \mathsf{ICF} \Leftarrow \mathsf{L}_F \,]$ :
 $  \mathsf{L}_F\overset{\ref{isomorphicduals}\mathrm{.v}} \Longrightarrow \mathsf{ICF}_{\kappa}^{\,'}$ 
 $\underset{(2)} {\overset{\ref{cardinallemma}\mathrm{.iii}}  \Longrightarrow} \mathsf{ICF}_{2}^{\,'}\overset{\ref{cardinallemma}\mathrm{.iv}}\Longrightarrow \mathsf{ICF} $.\end{proof}

\begin{theorem} \label{failure}
  Let $F$ be a field $($or division ring$)$ and  $\kappa=\big| F \big|$. Assume that  either   $\mathrm{(i)}$ $\exists\, \lambda \in \mathrm{Card}:[\, \kappa=2^{\,\lambda},\, \lambda \geq \mathfrak{c}\,]$  or 
  $\mathrm{(ii)}$   $[\, \kappa> \mathfrak{c},\, \mathsf{GCH}\,]$. Then
 \ $\mathsf{L}_F $  fails.\end{theorem}
 \begin{proof}
Case $\mathrm{(i)}$ is $\ref{isomorphicduals}\mathrm{.vi.b}$.  Case $\mathrm{(ii)}$ follows from  $\{ \ref{isomorphicduals}\mathrm{.ix.b},\, \ref{isomorphicduals}\mathrm{.viii.b}$, $\ref{isomorphicduals}\mathrm{.vii.c}\}$.  
\end{proof}

 \begin{theorem} \label{ICFdual3}
 $^{\ast}$ Let $F$ be  field $($or division ring$)$. Assume the $\mathsf{GCH}$. Then \quad
 $ \mathsf{L}_F\Leftrightarrow  \big| F \big| \leq \mathfrak{c}$.
  \end{theorem}
  \begin{proof} 
It follows from $\{ \ref{ICFdual},\, \ref{failure},\,\ref{isomorphicduals}\mathrm{.ix.a}  \}$.
 \end{proof} 
 
\begin{theorem}\label{theoremofcohen}
$($Gödel-Cohen$)$. Assume that  
 $ \mathsf{ZFC} $ is consistent. Then $:$ \\  $\mathrm{(i)}$  $\mathsf{ZFC} $ has a model on which  $\mathsf{GCH} $
  holds. \space So $:$  $\mathrm{(a)}$ $ \mathsf{ZFC}\nvdash  \neg \, \mathsf{GCH}$.  $\mathrm{(b)}$ $ \mathsf{ZFC}\nvdash  \neg \, \mathsf{ICF}$.   $\mathrm{(c)}$ $ \mathsf{ZFC}\nvdash  \neg \, \mathsf{CH}$. \\
   $\mathrm{(ii)}$  $\mathsf{ZFC} $ has a model  on which  $[\,2^{\aleph_0}=2^{\aleph_1}  \,] $. So $:$   $\mathrm{(a)}$ $ \mathsf{ZFC}\nvdash \mathsf{GCH}$.\quad $\mathrm{(b)}$ $ \mathsf{ZFC}\nvdash \mathsf{ICF}$. \quad  $\mathrm{(c)}$ $ \mathsf{ZFC}\nvdash \mathsf{CH}$.\\
 $\mathrm{(iii)}$  $ \mathsf{GCH}$  is 
  non-decidable in   $ \mathsf{ZFC} $.
  $\mathrm{(iv)}$  $ \mathsf{CH}$  is 
  non-decidable in   $ \mathsf{ZFC} $.  \quad $\mathrm{(v)}$ \underline{$\mathsf{ICF_{}} $ is non-decidable in $ \mathsf{ZFC} $}.
\end{theorem}
\begin{proof} $\ref{theoremofcohen}\mathrm{.i}$ was proven by Gödel \cite{godel} and $\ref{theoremofcohen}\mathrm{.ii}$ was proven by Cohen in \cite{cohen} (page 109). For the $\mathrm{(a),(b),(c)}$'s we additionally  apply $\ref{isomorphicduals}\mathrm{.ix}$ and $\mathrm{2.iii.4}$ (for $\ref{theoremofcohen}\mathrm{.ii.c}$).
Also $\ref{theoremofcohen}\mathrm{.iii },\, \ref{theoremofcohen}\mathrm{.iv },\, \ref{theoremofcohen}\mathrm{.v}$ follow from the above.
\end{proof}

\begin{remark} \label{LH} $:$  \\   $\mathrm{(i)}$ Theorem  $\ref{theoremofcohen}\mathrm{.v}$ is not popularly known, unlike  Theorem $\ref{theoremofcohen}\mathrm{.iv}$. The
author got aware of $\ref{theoremofcohen}\mathrm{.v}$ after contacting Andreas Blass and Asaf Karagila and exchanging several emails with them. \\
 $\mathrm{(ii)}$ The statement $[\,2^{\aleph_0}=2^{\aleph_1}  \,] $
is called  $ \mathsf{LH}$ (\emph{Luzin's Hypothesis}). It is known that $[\,\mathsf{MA}$ (\emph{Martin's Axiom}), $\neg \, \mathsf{CH}\,]\Rightarrow \mathsf{LH}$ (see \cite{mostowski}). Also $\mathsf{LH}\Rightarrow \neg \, \mathsf{ICF} $. Hence, it follows from  \ref{ICFdual} that \underline{$  \mathsf{L}_F$ fails  if   $[\, \mathsf{MA},\, \neg \, \mathsf{CH},\,  |F| \leq \mathfrak{c}\,] $.}    \\
$\mathrm{(iii)}$. The reader is encouraged to attempt to relax the conditions of  Case $\mathrm{(ii)}$ of Theorem \ref{failure}, i.e try to prove that $[\, \mathsf{L}_F$ fails, when $|F|> \mathfrak{c}\,]$
under weaker conditions than the $\mathsf{GCH}$.
 \\
$\mathrm{(iv)}$ Combining  $\ref{theoremofcohen}\mathrm{.v}$ with
 $\ref{ICFdual}$ we get the  Theorem $\ref{undecidable}$ below. The special  case for $\ref{undecidable}$ and $\ref{LH}\mathrm{.ii } $, where  $F$ is finite has been stated  in 
\cite{exchange} by by  Arturo Magidin, Asaf Karagila and Mark Saving. 
 \end{remark}
  
 \begin{theorem} \label{undecidable} $^{\ast}$
 Consider the statement $\mathsf{L}_F:$ $[\, \forall\, V\in F$-$\mathrm{Mod}\,][\,\forall\, W\in F$-$\mathrm{Mod}\,]$
 $[\, V^{\ast} \underset{\hat{F}}\cong W^{\ast}  \Rightarrow   V \underset{F}\cong W \,]$, 
for a fixed  field $($or division ring$)$ $F$. If $[\,  \mathsf{ZFC} $  is consistent,  $|F| \le \mathfrak{c}$  $]$ then  $[\, \mathsf{L}_F$ is a non-decidable statement in $\mathsf{ZFC}\,]$, i.e neither 
 $\mathsf{L}_F$  or $\neg \, \mathsf{L}_F$  can be proven in $\mathsf{ZFC}$.
\end{theorem} 

\section{Free modules with isomorphic duals, in Artinian rings}

 The next theorem is very indicative on what follows and was the main motivation behind this paper.
\begin{theorem} \label{artinianfreealg}Let $R$ be an affine algebra over a field $F$ and $\lambda$ be a cardinal $($or a set$)$.
 If $R$ is an Artinian ring  then  $R^{\lambda}$ is a free $R$-module i.e $[\, R^ {\lambda}\in \tmop{Free}(R)\,]$. 
\end{theorem}

\begin{proof} Let $R=F[x_1,x_2,...x_n]/J$ be an affine $F$-algebra and let $\phi:F \rightarrow R$ be defined by  $\phi(t)=t+J$.
By application of Theorem 6.54 of \cite{becker} we get that  $[\, \phi$ is f.p $]$(1), since $R$ is Artinian. It follows from $\ref{zornapplication2}$ that $[\, F^{\lambda}\in \tmop{Free}(F)\,]\Rightarrow$
$[\, \exists \,\kappa \in \mathrm{Card}: (F^{\lambda}\cong_F F^{[\kappa]}\,)]$
$ \underset{(1)}{\overset{\ref{rh2}\mathrm{.iii}}\Longrightarrow}$ 
$[\, R^{\lambda} \cong_R R^{[\kappa]}\,] \Rightarrow$
$[\, R^{\lambda}\in \mathrm{Free}(R)\,]$.
\end{proof}

We shall now try to extend the theorem above to  the class of all  Artinian commutative rings and to do so we will need some important theorems of commutative algebra. 
\begin{theorem} \label{collection} Let $R$ be a Noetherian commutative ring and $\lambda$ be a cardinal $($or a set$)$. Then $:$
\end{theorem}
\begin{enumerate}[label=(\roman*)]
\item If $R$ is local then $\tmop{Proj}(R)=\tmop{Free}(R)$. \\ 
$\mathbf{Proof:}$ This is a Theorem of Kaplansky. See   [17] .
\item If $R$  is Artinian then $\tmop{Proj}(R)=\tmop{Flat}(R)$. \\
$\mathbf{Proof:}$ This is Theorem 28.4 of \cite{fuller} which is a Theorem of Bass ( See \cite{bass1}).
\item  $ R^{\lambda} \in \mathrm{Flat}(R)$ \\ $\mathbf{Proof:}  $ Follows from 3.24 of \cite{enochs} since Noetherian rings are coherent.
\item If $R$ is local and Artinian then  $ R^{\lambda} \in \mathrm{Free}(R)$.   \\
$\mathbf{Proof:}$  Follows from the above since Artinian rings are Noetherian.
 \item If $R$ is  Artinian  then $R$ is a finite  product of local Artinian commutative rings . \\  $\mathbf{Proof}:$ See Theorem 27.15 of [11].         
\item If $R$ is  Noetherian  then $R$ is a finite  product of connected Noetherian  rings, i.e  $R=\overset{n}{\underset{k=1}\Pi} R_k $ for some connected Noetherian commutative rings $R_k.$ Also $ R$ is Artinian 
$\Leftrightarrow [\, R_k$ is Artinian $ ]\, \forall \, k\in \mathrm{T}_n$. \\
$\mathbf{Proof}:$ Well known, see at page 37 of  \cite{rao}.
\end{enumerate}
Note that $\ref{collection}\mathrm{.i}$, $\ref{collection}\mathrm{.ii}$, $\ref{collection}\mathrm{.iii}$ are valid even without the commutativity assumption on $R$ \\
The next lemma is very important for this article.

\begin{lemma}\label{strange}
$^{\ast}$ Let $X,\, \Lambda$ be  sets and  $R_1,R_2,\dots,R_n$ be commutative rings and let $R=\overset{n}{\underset{k=1}\Pi} R_k $ such that \\ 
$\mathrm{(i)} $ $[ \,(R_k)^X\in \mathrm{Free}(R_k) \,]$ $\forall \,  k\in  \mathrm{T}_n$  \quad $\mathrm{(ii)}$   $[\, \dim_{R_k} (R_k)^X=\big|\Lambda \big| \, ]$  $\forall k\in  \mathrm{T}_n$. Then $:$ \\
   $$[\,R^X\in \mathrm{Free}(R)\,] \,\, \text{and}\,\,[\, \dim_R (R^X)=\big|\Lambda\big|\,] $$
\end{lemma}
\begin{proof}
Let $ k\in  \mathrm{T}_n$. It follows from the assumptions that there is a  1-1 function $G_k :\Lambda\rightarrow (R_k)^X$ such that $G_k(\Lambda)$ is an  $R_k$-basis of $(R_k)^X$, i.e 
 $[\, G_k(\Lambda)$ is a linearly independent subset of $ (R_k)^X \,]$ (1)   and 
$[\, \langle \, G_k(\Lambda) \, \rangle_{R_k}= (R_k)^X \,]$ (2). For $k\in \mathrm{T}_n$, we consider the r.h $\pi_k :R\rightarrow R_k$ defined in $\mathrm{1.i}$  The function  
 $h:\overset{n}{\underset{k=1}\Pi}  ( R_k)^X \rightarrow $
 $( \overset{n}{\underset{k=1}\Pi}    R_k)^X=R^X$ defined by $ h(f)=(\pi_1\circ f,\pi_2\circ f,\dots, \pi_n\circ f)$ is clearly a bijection. Consider the function
   $G:\Lambda\rightarrow R^X$  be defined by $ G(\lambda)=h(G_1(\lambda),G_2(\lambda),\dots,G_n(\lambda))$. Then $[\,G$ is 1-1 $]$ (3), since $h,\,G_k$'s are so.\\
$\underline{\mathrm{Claim} \, 1_{}} :$ $[\, G(\Lambda)$ is an $R$-linearly independent subset of $R^X\,]$  \\
Proof: Let $r=(r_1,r_2,\dots,r_m)\in R^m$ be such that $ \overset{m}{\underset{k=1}\Sigma} r_k\cdot G(\lambda_k)=0$, where $r_k=(r_k^{i})_{i\in \mathrm{T}_n }\in R$. Then            
  $ [ \overset{m}{\underset{k=1}\Sigma} r^k_i\cdot G_i(\lambda_k)=0 \, ] \, \forall i\in \mathrm{T}_n$ 
$\overset{(1)}\Leftrightarrow$ 
$[(r^k_i)_{i\in \mathrm{T}_m }=0] \, \forall i\in \mathrm{T}_n\Leftrightarrow$   
$ [\,r=0\,] $. Therefore  $[\, G(\Lambda)$ is an $R$-linearly independent subset of $R^X\,]$ .      \\          
$\underline{\mathrm{Claim \, 2_{}} }$: $[\, \langle\, G(\Lambda)\, \rangle_{R}=R^X\,]$ i.e $ G(\Lambda)$ generates  $R^X$ over $R$. \\
Proof : Let  $ f\in R^X $ and  $i\in \mathrm{T}_n$. Let $f_i=\pi_i\circ f\in R_i^X$. It  follows from (2) that   $[\,  \exists  \,(r^k_i)_{ k\in \mathrm{T}_m}\in (R_i)^X :f_i= \overset{m}{\underset{k=1}\Sigma} r^k_i\cdot G_i(\lambda_k)\,]$ (5). For $ k\in \mathrm{T}_m$, we define  $r_k=(r^k_i)_{i\in \mathrm{T}_n}\in R $. Then it follows from (5) that
 $[\,f=\overset{m}{\underset{k=1}\Sigma} r_k \cdot G(\lambda_k)\in \langle\, G(\Lambda)\,\rangle_R \,]$. Hence 
 $[\, \langle\, G(\Lambda)\, \rangle_{R}=R^X\,]$. \\
 Therefore it follows from  $\{ \mathrm{Claim \, 1},\, \mathrm{Claim \,2},\, (3) \}$ that  
   $G(\Lambda)$ is an $R$-basis for $R^X$. Also $\big| G(\Lambda) \big|\overset{(3)}=\big| \Lambda\big|$.  Hence 
  $\{R^X\in \tmop{Free}(R) \,  \text{ and}   \, \dim_R (R^X)=\big|\Lambda \big|\,  \} $. 
\end{proof}

\begin{lemma} \label{noetherian} Let $R$ be a Noetherian commutative ring   and let  $\lambda$ be an infinite cardinal $($or  set$)$. If  
   $[\, R^{\lambda}\in \tmop{Free}(R)\,,\big| R \big | \le \mathfrak{c} \,]$ then  
   $[\,   \dim_R( R^{\lambda})=2^{\lambda} \,]$                         
 \end{lemma}
 \begin{proof} Let $\mathfrak{m}\in \tmop{max}(R)\neq \emptyset$ 
 and let $\phi=\pi^R_\mathfrak{m}:R\rightarrow \dfrac{R}{\mathfrak{m}}$. Then  $[\, \phi$ is a f.p $]$(1), since $R$ is Noetherian. Also  $\dfrac{R}{\mathfrak{m}}=F$ is a field and  $[\, \big| F \big|=\big| \dfrac{R}{\mathfrak{m}} \big|\le \big| R \big| \le \mathfrak{c}\,] $  $
 \Rightarrow [\, \big| F \big| \leq \mathfrak{c}\,]$ (3). It follows from   $\{ \ref{EKextension}\mathrm{.iv},(3) \}$ that
  $[\, \dim_F(F^{\lambda})=2^{\,\lambda}]\,(4)$. By assumption
$[\,  R^{\lambda}\in \tmop{Free}(R)\,]  \Rightarrow [\, \exists \, \kappa \in \mathrm{Card}$ :
 $ (R^{\lambda} \cong_R R^{[\kappa]})\,(5)\,]$ 
 $\underset{(1)}{\overset{\ref{rh2}}\Rightarrow}$ 
 $[\, F^{\lambda}\underset{F}\cong F^{[k]}\,]\, (6)\Rightarrow $  
  $[\,  \kappa=\dim_F (F^{[\kappa]})\overset{(6)}=\dim_F (F^{\lambda})\overset{(4)}=2^{\,\lambda} \,]\Rightarrow [\,\kappa=2^{\lambda} \,]\, (7)$. Similarly it follows from $(5)$  that $[\,\kappa=\dim_R( R^{\lambda})\,] \overset{(7)}\Rightarrow$
 $[\, \dim_R( R^{\lambda})=2^{\, \lambda }\,]$.  \end{proof}

\begin{theorem} \label{artinianfree} $^{\ast}$ Let  $R$ be an Artinian commutative ring such that
 $\big| R \big | \le \mathfrak{c}$  and let  $\lambda$ be  an infinite cardinal $($or a set$)$. Then  $R^{\lambda}\in \tmop{Free}(R) \,$   and   $\, \dim_R (R^{\lambda})=2^{\, \lambda } $.      
\end{theorem}

  \begin{proof} It follows from $\ref{collection} \mathrm{.v}$  that $R=\overset{n}{\underset{k=1}\Pi} R_k $  for some $[$local Artinian commutative  rings $R_k \,]$(1), hence  $[$Noetherian  commutative rings $R_k \,]$(2).  Let  $ k\in \mathrm{T_n}=\{1,2,\dots n \}$. It follows from $\{ (1),\ref{collection}\mathrm{.iii} \}$ that 
$[\,(R_k)^{\lambda}  \in \tmop{Free}(R_k)\,]$(3). 
 Also $[\, \big| R_k \big| \le \big| R \big| \le \mathfrak{c}\,]\Rightarrow$ 
 $[\, \big| R_k \big| \le \mathfrak{c} \,] $ (4). Therefore it follows from 
 $\{(2),(3),(4), \ref{noetherian} \}$ that 
$ [\, \dim_{R_k} (R_k)^{\lambda}=2^{\, \lambda } \,] \, \forall \,  k\in  \mathrm{T}_n\overset{\ref{strange}} \Longrightarrow $ 
 $ [\,R^{\lambda}\in \mathrm{Free}(R) \,\,\,$   and \space    $\, \dim_R (R^{\lambda})=2^{\, \lambda } \,]$.
\end{proof}

\begin{lemma} \label{strange2}
Let $X,\Lambda$ be  sets,  $R_1,R_2,\dots,R_n$ be commutative rings and   $R=\overset{n}{\underset{k=1}\Pi} R_k $.  \\ If
  $[\, R^X\in \mathrm{Free}(R), \,$  $\, \dim_R (R^X)=\big|\Lambda\big| \,]$ then
  $[\, (R_k)^X\in \mathrm{Free}(R_k), \,$  $\, \dim_{R_k} (R_k)^X=\big|\Lambda\big| \,] \, \forall \,  k\in \mathrm{T}_n$.  
\end{lemma}
\begin{proof} For $ k\in \mathrm{T}_n$, we consider the r.h $\pi_k:R\rightarrow R_k$, defined in $1\mathrm{.i}$. Then 
$[\, \phi_k$ is f.p $]$ (1). By  assumption  we have :   
$[\, R^X \underset{R}\cong R^{[\Lambda]}\,] $ $\underset{(1)}{\overset{\ref{rh2}\mathrm{.iii}}\Rightarrow}$
$[\,(R_k)^X\cong(R_k)^{[\Lambda]}\,]\Rightarrow$ 
$[ \,(R_k)^X\in \mathrm{Free}(R_k) \,$   and   $\, \dim_{R_k} (R_k)^X=\big|\Lambda\big| \,]$.
\end{proof}

\begin{lemma} \label{counterexample1} There is an Artinian commutative ring $R$ such that $ [\, R^{\omega} \notin \mathrm{Free}(R)\,] $.
\end{lemma}
\begin{proof}
Let $[\, F$ be a field such that  $|F  | \leq \mathfrak{c}\,]\, \overset{\ref{EKextension}\mathrm{.iv}} \Longrightarrow $
 $[\, \dim_F (F^{\omega}) =|2^{\omega}|=\mathfrak{c} \,\,]$ (1). For example  $F$ could be any finite field. 
Let $K$ be a field such that  $[\, |K| > \mathfrak{c}\,] \overset{\ref{EKextension}\mathrm{.i}} \Longrightarrow[\,  \dim_K (K^{\omega})=
| K^{\omega}| \geq \ |K |> \mathfrak{c} \,]\Rightarrow$  
$[\,  \dim_K (K^{\omega}) > \mathfrak{c}\,]\,(2) $. Note that the  \emph{Löwenheim–Skolem  Theorem} (See  \cite{marker})  guarantees the existence of such a field K. Let now  $R=F\times K$. Assume that $[ \, R^{\omega} \in \mathrm{Free}(R)\,] $  $\overset{\ref{strange2}} \Longrightarrow $
  $[\, \dim_{F} (F^{\omega})=  \dim_K (K^{\omega})\,]\overset{(1),(2)}\Longrightarrow [\, \mathfrak{c}>\mathfrak{c}\,] $, which is a contradiction. Hence $ [ R^{\omega} \notin \mathrm{Free}(R)\,] $.
\end{proof}

 \begin{definition} \label{defdimension}
 Let  $M$ be an $R$-module over be a commutative ring and 
 $\kappa\in \mathrm{Card}$. Then  \\
 $[\, \dim_R(M)=\kappa\,]\Leftrightarrow$
 $[\, M \cong_R R^{[\kappa]}\,]$.
 \end{definition}
\begin{theorem} \label{EKartinian}
$^{\ast}$Let $(R,\mathfrak{m})$ be  local Artinian commutative ring and $\lambda$ be infinite cardinal.  Then \\$[R^{\lambda}\in \tmop{Free}(R) \,  \text{ and}   \, \dim_R (R^{\lambda})=\big|\dfrac{R}{\mathfrak{m}}\big|^{\lambda} \,]$ 
\end{theorem} 
\begin{proof}
We  consider the r.h $\phi=\pi^R_{\mathfrak{m}}:R\rightarrow \dfrac{R}{\mathfrak{m}}=F$. So  $[\, \phi$ is f.p $]$(1), because  $R$ is a Noetherian ring (since $R$ is Artinian) and $[\, F$ is a field $]$ (2).
 By asumption,
 $[\, R$ is local Artinian $]$ 
  $\overset{\ref{collection}\mathrm{.v}}\Rightarrow$
$[\, R^{\lambda}\in \tmop{Free}(R)\,]\Rightarrow$ 
$[ \, \exists \, \kappa \in \mathrm{Card} :$  $\kappa=\dim_R (R^{\lambda})$ and $R^{\lambda}\underset{R}\cong R^{[\kappa]}\, ]\underset{(1)}{\overset {\ref{rh2} \mathrm{.iii}}\Rightarrow}$   
$[\, F^{\lambda}\underset{F}\cong F^{[\kappa]} \,](4) $ 
$\overset{(2)} \Rightarrow$  
$[\, \dim _R (R^{\lambda})=\kappa=\dim _F (F^{[\kappa]})=\dim _F (F^{\lambda})$
$\overset{\ref{EKextension}\mathrm{.i} }=\big| F^{\lambda} \big|$
$=\big| F \big| ^{\lambda} =   \big| \dfrac{R}{\mathfrak{m}} \big|^{\,\lambda}\,] \Rightarrow$
$\boxed{ \dim_R (R^{\lambda}) =\big| \frac{R}{\mathfrak{m}}\big|^{\,\lambda} }$ 
\end{proof}

\begin{remark} The theorem above, is a generalization of the \emph{ Erdős-Kaplansky Theorem} to local Artinian commutative rings, since any field $F$ is a local Artinian ring with maximal ideal 
$\mathfrak{m}=\mathbf{0}$. Therefore $\dim_F(F^{\lambda})\overset{\ref{EKartinian}}= \big|\dfrac{F}{\mathbf{0}} \big|^{\lambda}=$
$\big| F \big|^{\,\lambda}=\big| F^{\lambda} \big|$, when 
$\lambda \geq \omega$.
\end{remark}
\begin{theorem}\label{important}
 Let $R$ be an Artinian commutative ring and $[\, \lambda$ be an infinite cardinal $]\,(1)$. Then $:$  
\end{theorem} 
\begin{enumerate}[label=(\roman*)]
\item There are Artinian local commutative rings $(R_k,\mathfrak{m}_k)$ such that 
 $ R=\overset{n}{\underset{k=1}\Pi} R_k$ \\ $\mathbf{Proof:}$
 This is $\ref{collection}.\mathrm{v}$.
\item$^{\ast}  \mathbf{Freeness\,\, Criterion}:$  \\
$\mathrm{(a)} \,[\,  R^{\lambda}\in \tmop{Free}(R)\,]\Leftrightarrow$
 $ [\,\big| \dfrac{R_k}{\mathfrak{m}_k}\big|^{\, \lambda}=\big| \dfrac{R_1}{\mathfrak{m}_1}\big|^{\,\lambda}]\, \forall k\in \mathrm{T_n}$  
 $ \, \mathrm{(b)} \,[\,  R^{\lambda}\in \tmop{Free}(R)\,]\Rightarrow         [\, \dim_R (R^{\lambda})=\big| \dfrac{R_1}{\mathfrak{m}_1}\big|^{\, \lambda}\,]$
  
$\mathbf{Proof:}$ Follows from $\{$\ref{strange}, $\,$\ref{strange2}, \ref{EKartinian}$\}$ .

\item $(3) \, [\, \big| \dfrac{R_k}{\mathfrak{m}_k}\big|^{\,\lambda} \le 2^{\,\lambda}]\, \forall \, k\in \mathrm{T_n} $
  $\Rightarrow[R^{\lambda}\in \tmop{Free}(R) \,  \text{ and}   \, \dim_R (R^{\lambda})=2^{\, \lambda} \,]$ \\  $\mathbf{Proof:}$ 	Let $ \mu_k=\big| \dfrac{R_k}{\mathfrak{m}_k}\big|$ for $k\in \mathrm{T_n} $.
  Then  
$[\, 2 \le \mu_k \,]\Rightarrow$ 
$[\, 2^{\, \lambda}  \le \mu_k ^{\, \lambda}\underset{(3)}\le (2^{\, \lambda})^{\lambda}=2^{\, \lambda \cdot \lambda}\overset{(1)}=2^{\, \lambda}\,] \Rightarrow $
$[\, 2^{\, \lambda}\le \mu_k^{\lambda}\le 2^{\, \lambda}\,]$   
 $\Rightarrow \{[ 
\,\mu_k^{\, \lambda}=2^{\, \lambda} \,]\, \forall \, k\in\mathrm{T_n} \} $
 (5). So the claim follows from $\{$ (5), $\ref{important}.\mathrm{ii}\}$.
 
 \item  $[\,\big| R \big| \le 2^{\, \lambda}]\,(6) $
  $\Rightarrow[R^{\lambda}\in \tmop{Free}(R) \,  \text{ and}   \, \dim_R (R^{\lambda})=2^{\, \lambda} \,]$  \\
$\mathbf{Proof:}$ Let $k\in \mathrm{T_n} $ then 
$[\, \big| \dfrac{R_k}{\mathfrak{m}_k}\big| \le \big| R_k \big|  \le \big| R \big|\,]\Rightarrow $
$[\, \big| \dfrac{R_k}{\mathfrak{m}_k}\big|^{\, \lambda} \le \big| R \big|^{\, \lambda} \underset{(6)}\le (2^{\, \lambda})^{\lambda}=2^{\, \lambda \cdot \lambda}\overset{(1)}=2^{\, \lambda}\,]\Rightarrow$ \\
$\{ \,  [\, \big| \dfrac{R_k}{\mathfrak{m}_k}\big|^{\,\lambda} \le 2^{\, \lambda}\,]\,\forall k\in \mathrm{T_n \}} \,(7)$. Therefore the claim follows from $\{$(7), \ref{important}$\mathrm{.iii} \}$.
\item Let $\xi$ be a cardinal such that $[\, \xi \geq \lambda \geq \omega\,]$  (8). Then 
$[\,  R^{\lambda}\in \tmop{Free}(R)\,]\Rightarrow [\, R^{\xi}\in \tmop{Free}(R)\,]$  \\
$\mathbf{Proof:}$ Let $\mu_k$ as in the proof of $\mathrm{(iii)}$.
Then  $[\,  R^{\lambda}\in \tmop{Free}(R)\,]\overset{(\mathrm{ii.a})}\Leftrightarrow$
$\{\,[\,\mu_k^{\, \lambda}= \mu_1^{\,\lambda}\,]\, \forall k \in \mathrm{T_n}   \,] \, \} \underset{(8)}     {\overset{\ref{important}\mathrm{.ii} }\Longrightarrow}$
$\{\,[\,\mu_k^{\, \xi}= \mu_1^{\, \xi}\,]\, \forall k \in \mathrm{T_n}   \,] \, \}$     $\overset{(\mathrm{ii.a})}\Leftrightarrow$
$[\,  R^{\xi}\in \tmop{Free}(R)\,]$
\end{enumerate}

 \begin{remark}Note that \ref{artinianfree} follows from 
 $\ref{important}.\mathrm{iv}$ since 
$\{ [\,  \big| R \big| \le \mathfrak{c}\,]\Rightarrow$
 $[\,\big| R \big| \le \mathfrak{c}=2^{\,\omega} \underset{(1)}\le 2^{\, \lambda}\,]\Rightarrow[\, \big| R \big| \le 2^{\, \lambda}\,] \}$. Also note that $\ref{important}\mathrm{.v}$ is known to \cite{oneil}.  \end{remark}
 
\begin{definition} \label{definitions2}
Let $R$ be  commutative ring. Then $:$   \\
 $\mathrm{(i)}\,\mathsf{L}_R$  $\Leftrightarrow[\, \forall \mu \in \mathrm{Card} \,]$
$[\, \forall \lambda \in \mathrm{Card} \,]$
  $[\,(R^{ \mu}\cong_R R^{ \lambda}) \Rightarrow( \mu=\lambda) \,]$.   \\
$\mathrm{(ii)}\,\mathsf{L}_R^w$  $\Leftrightarrow[\, \forall \mu \in 
\mathcal{L}_1' \,]$
$[\, \forall \lambda \in \mathcal{L}_1' \,]$
  $[\,(R^{ \mu}\cong_R R^{ \lambda}) \Rightarrow( \mu=\lambda) \,]$. ($\mathcal{L}_1'$ is the class of non-$\omega$-mesurable cardinals).

 \end{definition}
 \begin{remark}
 Note that \space $\mathsf{L}_R \Leftrightarrow[\,$Any two free $R$-modules with isomorphic duals are themselves isomorphic$\,]$. Definitions $\ref{definitions}\mathrm{.vii}$ and $\ref{definitions2}\mathrm{.i}$ are compatible since all vector spaces are free modules.
 \end{remark}
 \begin{lemma} \label{projections}
 Let $R_1,\ R_2,\, \dots R_n$ be commutative rings, $ R=\overset{n}{\underset{i=1}\Pi} R_i$ and  $\pi_i: R\rightarrow R_i$ as in $\mathrm{1.i}$.
 Let also $M,\, N$ be $R$-modules. Then $[\, M \cong_R N\,]\Leftrightarrow$
 $[\,M _{\pi_i}\cong_{R_i} N_{\pi_i}\,]\,\forall\, i\in \mathrm{T}_n$.   \end{lemma}
 \begin{proof}Let  $i\in \mathrm{T}_n$. 
 Consider the r.h $\xi_i: R\rightarrow R_i$ defined by 
 $\xi_i(x)=x\cdot e_i$ and the r.h  $\phi_i:R \rightarrow R$ defined by  $\phi_i=\xi_i \circ \pi_i$. Let $L \in
 R$-$\mathrm{Mod}$.   It is simple that $[\,  (e_i L)\cap(e_k L)=\mathbf{0},$ if $ i\neq k\,]\,(1)$ 
 and  $[\, L_{\phi_i}\cong_R e_i \, L\,]\,(2)$. Hence
  $[\, L=\overset{n}{\underset{i=1}\Sigma} e_i \, L \overset{(1)} {\underset{R}\cong} $
$\overset{n}{\underset{i=1} \oplus}  e_i \, L $
$\overset{(2)}{\underset{R}\cong}$
$\overset{n}{\underset{i=1} \oplus} L_{\phi_i}=$
$\overset{n}{\underset{i=1} \oplus} L_{\xi_i \circ \pi_i}\underset{R}\cong$ 
$\overset{n}{\underset{i=1} \oplus} (L_{\pi_i})_{\xi_i}\,]
 \Rightarrow$
$[\, L\underset{R}\cong \overset{n}{\underset{i=1} \oplus} (L_{\pi_i})_{\xi_i}\,]$ (3). Therefore the claim follows from (3), since we can recover any module  $L$ from the extensions  $L_{\pi_i}$'s
 
 \end{proof}

 \begin{theorem} \label{dualICFlocal}
 $^{\ast}$ Let $(R,\mathfrak{m})$ be  local Artinian commutative ring, $ \dfrac{R}{\mathfrak{m}}=F $ and $\kappa=\big| F \big|$. Then:
  \end{theorem}
  \begin{enumerate} [label=(\roman*)]
\item $[\, R^{\mu} \cong_R R^{\lambda} \,] \Leftrightarrow$
$[\, F^{\mu} \cong_F F^{\lambda}\,]$ \\
$\mathbf{Proof:}$ $\cdot$ Let $[\,\kappa \geq \omega,\,\lambda \geq \omega  \,]$ (2). Applying $\{ (2),\, \ref{EKartinian}  \}$,  we get that $[\, \dim _R(R^{\mu})=\kappa^{\, \mu}\,]$ (3) and 
$[\, \dim _R(R^{\lambda})=\kappa^{\, \lambda}\,]$ (4). Also applying
 $\{ (2),\, \ref{EKextension}\mathrm{.i} \}$, we get that $[\, \dim _F(F^{\mu})=\kappa^{\, \mu}\,]$ (5) and  $[\, \dim _F(F^{\lambda})=\kappa^{\, \lambda}\,]$ (6).  Therefore  $[\, R^{\mu} \cong_R R^{\lambda} \,] \Leftrightarrow$
  $[\, \dim _R(R^{\mu})= [\, \dim _R(R^{\lambda})\,]$ (7), by application of $\{\ref{rankiswelldefined},\, \ref{EKartinian}  \}$.
  Also  $[\, F^{\mu} \cong_F F^{\lambda} \,] \Leftrightarrow$
  $[\, \dim _F(F^{\mu})= [\, \dim _F(F^{\lambda})\,]$ (8). Therefore  
  by $\{(3),(4),(5),(6),(7),(8)  \} $ it follows that  
 $[\, R^{\mu} \cong_R R^{\lambda} \,] \Leftrightarrow$
$[\, F^{\mu} \cong_F F^{\lambda}\,]$. \\
$\cdot$ Let $[\, \mu <\omega \,]\,(9)$. Then  
$[\, R^{[\mu]}\overset{(9)}= R^{\mu} \cong_R R^{\lambda} \,] \Leftrightarrow$
$[\, R^{\lambda} \cong_R R^{[\mu]} \,]\, (10) \overset{(9),\ref{lost}}\Longrightarrow $
$[\,\lambda <\omega \,] \overset{(10)}\Rightarrow[\,  R^{[\lambda]}=R^{\lambda}\cong_R R^{[\mu]}\,]\Rightarrow$
$[\, R^{[\lambda]} \cong R^{[\mu]}\,]\overset{\ref{rankiswelldefined}}\Rightarrow$
$[\,\lambda=\mu\,]$
$\Rightarrow[\, R^{\mu} \cong_R R^{\lambda}\,]$. Therefore               $[\, R^{\mu} \cong_R R^{\lambda}\,]\Leftrightarrow[\,\lambda=\mu \,]\,(11)$. In a similar way we prove that $[\, F^{\mu} \cong_F F^{\lambda}\,]\Leftrightarrow[\,\lambda=\mu \,]\,(12)$. So the claim follows from  $\{ (11),(12) \}$. \\
$\cdot$ Let $[\, \lambda <\omega \,]$. The proof is similar to the above. We have covered all cases. 

\item $ \mathsf{L}_F   \Leftrightarrow \mathsf{L}_R $ \\
$\mathbf{Proof:}$ Follows from Theorem  $\ref{dualICFlocal}\mathrm{.i}$ and Definition  $\ref{definitions2}\mathrm{.i}$
\item 
If  $\kappa\leq \mathfrak{c}$ then \space
$\mathrm{(a)} \,\,[\,\mathsf{L}_R \Leftrightarrow   \mathsf{ICF}\,] $
\space $\mathrm{(b)}$ If $\mathsf{ZFC}$ is consistent then $\,\mathsf{L}_R $ is non decidable in $\mathsf{ZFC}$ \\
$\mathbf{Proof:}$ $\mathrm{(a)}$ follows from 
$\{\ref{dualICFlocal}\mathrm{.ii},\, \ref{ICFdual}\}$  and $\mathrm{(b)}$ follows from  
  $\{\ref{dualICFlocal}\mathrm{.iii.a},\,\ref{theoremofcohen}\mathrm{.iii}\}$.

\item  Assume that  either   $\mathrm{(i)}$ $\exists\, \lambda \in \mathrm{Card}:[\, \kappa=2^{\,\lambda},\, \lambda \geq \mathfrak{c}\,]$  or  $\mathrm{(ii)}$   $[\, \kappa> \mathfrak{c},\, \mathsf{GCH}\,]$. Then  $\mathsf{L}_R $  fails \\
$\mathbf{Proof:}$ Follows from $\{\ref{dualICFlocal}\mathrm{.ii},\,\ref{failure}\}$ 
\end{enumerate}

\begin{lemma} \label{products}
 Let $R_1,\ R_2,\, \dots R_n$ be commutative rings,  $ R=\overset{n}{\underset{i=1}\Pi} R_i$ and $\mu,\, \lambda$ be cardinals.
\end{lemma}
\begin{enumerate} [label=(\roman*)]
\item $[\, R^{\mu} \cong_R R^{\lambda} \,] \Leftrightarrow$
$[\, R_i^{\mu} \cong_{R_i} R_i^{\lambda}\,] \,\forall\,  i\in \mathrm{T}_n$. \\
$\mathbf{Proof:}$ Let $i\in \mathrm{T}_n$ and $\kappa\in \mathrm{Card}$. Then it follows from  $\ref{rh}\mathrm{.ii}$ that $[\, (R^{\kappa})_{\pi_i} \cong_{R_i} R_i^{\kappa}\,]$
(3), since $\pi_i$ is f.p. Therefore the claim follows from $\{(3),\, \ref{projections} \}$
\item   $ \mathsf{L}_R   \Leftrightarrow [\, \mathsf{L}_{R_i}\,] \,\forall \, i\in \mathrm{T}_n$.   \\
 $\mathbf{Proof:}$  Follows from $\{ \ref{definitions2}, \,  \ref{products}\mathrm{.i}\}$
 \item Assume that $[\, (R_i,\mathfrak{m}_i)$ is local$,\, \big|\dfrac{R_i}  {\mathfrak{m}_i} \big|\leq \mathfrak{c} \, ]\,\forall \, i\in \mathrm{T}_n$. Then  \space $\mathsf{L}_R \Leftrightarrow   \mathsf{ICF} $   \\ 
$\mathbf{Proof:}$   Let $F_i=\dfrac{R}{\mathfrak{m}_i}$.
Then $[\, F_i$ is a field and  $\big| F_i \big|\leq \mathfrak{c}\,]\overset{\ref{ICFdual}}\Rightarrow$
$[\,\mathsf{L}_{F_i} \Leftrightarrow   \mathsf{ICF}\,] $ (2).
Therefore  $\mathsf{L}_R \overset{(\mathrm{ii.a})} \Leftrightarrow 
\{ \, [\,\mathsf{L}_{R_i}\,]\, \forall \, i\in \mathrm{T}_n\,\}\overset{(2)}\Leftrightarrow$
$[\, \mathsf{ICF}\,]\, \forall \, i\in \mathrm{T}_n\Leftrightarrow \mathsf{ICF}$. Hence  $\mathsf{L}_R\Leftrightarrow \mathsf{ICF} $.
\item Assume that $\exists\, i\in \mathrm{T}_n$ such that  $(R_i,\mathfrak{m}_i)$ is local and  either   $\mathrm{(i)}$ $\exists\, \lambda \in \mathrm{Card}:[\,   \big|\dfrac{R_i}  {\mathfrak{m}_i} \big| =2^{\,\lambda},\, \lambda \geq \mathfrak{c}\,]$   or  $\mathrm{(ii)}$   $[\,  \big|\dfrac{R_i}  {\mathfrak{m}_i} \big|> \mathfrak{c},\, \mathsf{GCH}\,]$.
 Then $[\, \mathsf{L}_R$ fails$\,]$ \\
$\mathbf{Proof:}$ Follows from $\{ \ref{products}\mathrm{.ii.a} ,\,  \ref{dualICFlocal}\mathrm{.iv} \}$
\end{enumerate}

\begin{theorem} \label{ICFdual2}
$^{\ast}$Let $R$ be Artinian commutative ring such that  $\big| R \big|\leq\mathfrak{c}$.  Then $:$ \\  $\mathrm{(i)} \,\,\mathsf{L}_R \Leftrightarrow   \mathsf{ICF} $
\quad\quad \quad \quad \quad\quad\ \quad \quad\quad  $\mathrm{(ii)}$
If $\mathsf{ZFC}$ is consistent then  $[\,\mathsf{L}_R $ is non-decidable in $\mathsf{ZFC}\,]$
\end{theorem}
 
\begin{proof}
 It follows from  $\ref{important}\mathrm{.i}$
 that there are Artinian local commutative rings $(R_i,\mathfrak{m}_i)$ such that 
 $ R=\overset{n}{\underset{k=1}\Pi} R_k.$  So $[\, \, \big| \dfrac{R_i}{\mathfrak{m}_i} \big|\leq  \big| R_i \big| \leq  \big| R \big|\leq \mathfrak{c}\,]\Rightarrow$
  $[\, \, \big| \dfrac{R_i}{\mathfrak{m}_i} \leq \mathfrak{c}\,]$  (3). So $\mathrm{(i)}$ follows from $\{(3),\,\ref{products}\mathrm{.iii}\} $. Also $\mathrm{(ii)}$ follows from $\{ \mathrm{(i)},\,\ref{theoremofcohen}\mathrm{.iii}\} $.
\end{proof}
We close this section with an application of projectivity and flatness to elementary commutative algebra. 
\begin{theorem} \label{primitive}
Any primitive polynomial  over any Artinian commutative ring   divides some monic one. \end{theorem}

\begin{proof} 
Let $R$ be a commutative ring and $f=\overset{n}{ \underset{k=0}\Sigma} f_k\cdot t^k\in R[t]$ be such that $f_n\neq 0$, then $[f$ is monic iff $f_n=1]$. The $R$-content of $f$ is 
$c_R(f) \overset{\text{def}}=\overset{n}{ \underset{k=0} \Sigma}   R\cdot f_k $ i.e it is the $R$-ideal generated by the coordinates of the polynomial and $[f$ is $R$-primitive iff $c_R(f)=R]$. We also define the ring  $R^f \overset{\text{def}}=\dfrac{R[t]}{R[t]\cdot f}.$ We consider the ring homomorphism $\psi=\psi^R_f:R\rightarrow R^f$ defined by 
$\psi (r)=r+R[t]\cdot f$, then $R^f$ becomes an $R$-module.  We will use the following:\\
(1): $R^f\in \tmop{Flat}(R)\Leftrightarrow c_R(f)$ is principal idempotent. $\mathrm{Proof:}$ This is the \emph{Nagata Ohi Theorem}. See  \cite{ohi}.
(2): Assume that $R^f\in \tmop{Flat}(R)$.
Then $[R^f\in \tmop{Proj}(R)\Leftrightarrow$  $\psi^R_f$ is integral r.h $]$.  $\mathrm{Proof:}$  See  \cite{vasc}.       \\
Assume now that [\underline {$R_{}$ is Artinian} (3) and \underline{$f$ is primitive} (4) $]\Rightarrow c_R(f)=R\Rightarrow$
$c_R(f)$ is principal idempotent $\overset{(1)}\Rightarrow$
 $R^f\in \tmop{Flat}(R) \underset{(3)} { \overset{\ref{collection}.\mathrm{ii}  }=}\tmop{Proj}(R)\Rightarrow$
 $R^f\in \tmop{Proj}(R)\overset{(2)}\Rightarrow$
$\psi^R_f$ is integral r.h $\overset{\text{def}}\Rightarrow$
$x=t+R[t]\cdot f$ is $R$ integral  $\Rightarrow$
$\exists \,g=\overset{m}{ \underset{k=0}\Sigma} g_k\cdot t^k\in R[t]:$
$[0=g(x)=\overset{m}{ \underset{k=0}\Sigma} g_k\cdot t^k+R[t]\cdot f,\, g_m=1]\Leftrightarrow$ The ideal $R[t]\cdot f$ contains some monic polynomial  $g\Leftrightarrow$ The polynomial $f$ divides some monic polynomial.
\end{proof}
Note that the above Theorem \ref{primitive} is  probably known but the argument used in the proof is  innovative. From the last two sections  we claim Theorems  \ref{secondtheorem2}, $\ref{half}\mathrm{.vi}$, \ref{relief} as unknown.

\section{Free modules with isomorphic duals, in Noetherian rings}

\subsection{Slenderness conditions and freeness of dual}

\underline{$\mathrm{Question \,A }$} : If $[\, \lambda\in \mathrm{Card}_{\geq \omega} ,\,  R^{\lambda}\in  \mathrm{Free}(R),\, R$ is Noetherian commutative ring$\,]$ then  does $R$  have  to be Artinian?
In order to answer  it we will utilize the well known concept of slenderness, stated below and introduce the concept of w-slenderness.

\begin{definition} \label{slenderdef} Let $R$ be a commutative ring, 
$\lambda\in \mathrm{Card}$, $M\in R$-$\mathrm{Mod}$ and $g\in  \mathrm{Hom}_R(R^{\lambda},M)$. Then $:$ \\
$\mathrm{(i)}\,\mathfrak{t}_{\lambda}(g)= \{ i\in \lambda \, | \,g(e_i)\neq 0\}$. \quad\quad\quad $\mathrm{(ii)}\, M$ is a slender $R$-module iff 
$[\,| \mathfrak{t}_{\omega}(f)|< \omega\,]$
$\forall f \in \mathrm{Hom}(R^{\omega},M)$. \\
$\mathrm{(iii)}$ $\mathrm{Slend}(R)$ is the class of slender $R$-modules. \quad\quad\quad\quad  $\mathrm{(iv)}$  $R$ is a slender  ring iff $R\in \mathrm{Slend}(R)$. 
\end{definition}
\begin{lemma} \label{slenderdef2} Let $R$ be a slender commutative ring, $ \lambda\in \mathrm{Card}_{\geq \omega}$ and $ f \in (R^{\lambda})^{\ast}$. Then $| \mathfrak{t}_{\lambda}(f)|< \omega$.  \end{lemma} 
\begin{proof} Assume that $| \mathfrak{t}_{\lambda}(f)| \geq \omega$.
Then $\exists\,  A \in \mathcal{P}(\lambda):\{ \,  |A|=\omega , \, [\, f(e_x)\neq 0\, ]\, \forall \, x\in A\, \}$.  Let $g=f|_{R^A}$, be the restriction of $f$ to $R^A$. Then $[\,g(e_a) \neq 0\,]\, \forall\, a \in A$, which contradicts the slendernerss of $R$. Therefore  $| \mathfrak{t}_{\lambda}(f)| < \omega$.\end{proof}

The first ring that was proven to be slender is $R=\mathbb{Z}$, by Specker (See   \cite{specker}). 
\begin{theorem} \label{proj.slender} Let $\phi :R\rightarrow S$ be r.h,   $L\in S$-$\mathrm{Mod}$ and $M\in R$-$\mathrm{Mod}$. Then $:$ \\
$\mathrm{(i)}$ $[\, N \in \mathrm{Sub}_R M,\, M\in  
\mathrm{Slend}(R)\,]\Rightarrow  [\, N\in  \mathrm{Slend}(R)\,]$. \quad
$\mathrm{(ii)}$ $[\, L\in  \mathrm{Slend}(S)\,]\Rightarrow  [\, \mathrm{K}_{\phi} (L) \in  \mathrm{Slend}(R)\,]$.\\
$\mathrm{(iii)}$ If $R$ is slender ring then $\mathrm{Slend}(R)\subseteq \mathrm{Proj}(R)$. 
\end{theorem} 
\begin{proof} $\mathrm{(i)},\,\mathrm{(ii)} $ follow trivially from Definition $\ref{slenderdef}\mathrm{.ii}$. \\
 \underline{$\mathrm{Proof\, of \ (iii)}$} : It is hard to prove. See Theorem 3 and Corollary 2 of \cite{lady}.\end{proof}

\begin{theorem} \label{slendercollection}
Let $R,S$ be Noetherian commutative rings, $\phi: R\rightarrow S$ be r.h and  $K$ be any field .
\end{theorem}
 \begin{enumerate}  [label=(\roman*)]
 \item If $R$ is Dedekind domain and $\tmop{Spec}(R)$ is infinite then $R$ is slender.   So $\mathbb{Z}$ is slender. \\
  $\mathbf{Proof:}$ See Example 4 of \cite{lady} .
 \item $K[t]$ is slender. \\ $\mathbf{Proof:}$   If  $K $ i s infinite then $\Lambda=\{x-a \,|\, a\in K \} \subseteq \tmop{Spec}(K[t])$ is infinite hence  $\tmop{Spec}(K[t])$ is infinite. If  $K $ is finite then there are   irreducible polynomials of any positive degree  hence  $\tmop{Spec}(K[t])$ is infinite . Hence the claim follows from $\mathrm{(i)}$ since  $K[t]$ is a Dedekind domain as a PID.
 \item  $[\, R$ is slender, $S \in \tmop{Proj}(R) \,]\Rightarrow$ 
 $[\, S$ is slender$\,]$. \\ $\mathbf{Proof:}$ It is $\ref{proj.slender}\mathrm{.ii}$.
 \item $[\, R$  is slender$\,]\Rightarrow[\,   R[t]$  is 
 slender$\,]$  \\
 $\mathbf{Proof:}$ 
 Consider the r.h \space  $\phi:R\rightarrow S=R[t]$, defined by $\phi(r)=r$. The set $\Delta=\{t^k\,|\,k\in \mathbb{N}  \}$ is clearly an $R$-basis for $S$, hence  $S\in \tmop{Free}(R)\subseteq \tmop{Proj} (R)\Rightarrow$
 $S\in \tmop{Proj}(R)\overset{\mathrm{(iii)}}\Rightarrow$
 $S=R[t]$ is a slender ring.
 \item If  $R$ is slender then $T=\dfrac{R[t]}{t^2 \cdot R[t]}$  is non-reduced (hence non-domain) and slender.\\
$\mathbf{Proof:}$ Consider the r.h  \space $\phi:R\rightarrow T$ defined by 
$\phi(r)=r+J$, where $J=t^2 \cdot R[t]$. The set $\Lambda=\{1+J,\,t+J \}$ is clearly an $R$-basis of $T\underset{R}\cong R_{\phi}$. Hence 
  $T\in \tmop{Free}(R)\subseteq \tmop{Proj} (R)\Rightarrow$
 $T\in \tmop{Proj}(R)\overset{\mathrm{iii}}\Rightarrow$
 $T$ is a slender ring. Let $x=t+J \in T$. Then $x^2=0_T$ (since $t^2 \in J$) and $x \neq  0_T$ ( since $t\notin J$). Therefore $R$ is not reduced.
 \item The polynomial ring in $n$ variables, $K[t_1,t_2,\dots t_n]$ is slender. \\ $\mathbf{Proof:}$ It  follows from
  \{$\mathrm{(ii)},\, \mathrm{(iv)}\}$.
 \item The ring $K[[t]]$  of formal power series over $K$ is not slender.  \\ $\mathbf{Proof:}$ See  \cite{kremer} .
  \item Let   $\phi$ be 1-1 and $S$ be a domain. Then  $[\, R$  is slender $\Leftrightarrow$  $S$  is slender $]$ \\  $\mathbf{Proof:}$ See Proposition 2.1 of  $\cite{jensen}$.
 \item If $R$ is a countable domain but not a field then $R$ is slender. \\ $\mathbf{Proof:}$  See Corollary 2.4  of  \cite{eklof}.
 \end{enumerate}

\begin{lemma}\label{reflexiveslender1} Let  $R$ is a slender commutative ring. Then  both  $R^{ \omega} $ and $ R^{[\omega]}$ are reflexive
 $R$-modules. Hence  
$(R^{\omega})^{\ast} \cong_R R^{[\omega]}$.
\end{lemma} 
\begin{proof}See Corollary $3.8$ of  \cite{eklof}.\end{proof}

\begin{lemma} \label{refnofree}
Let $R$ be a  Noetherian commutative ring and $\lambda$ be an infinite cardinal.\\ If $(R^{\lambda})^{\ast} \cong_R R^{[\lambda]}$   then 
$R^{\lambda}\notin \mathrm{Free}(R)$. 
\end{lemma}
\begin{proof} $[$Assume that $R^{\lambda}\in \tmop{Free}(R)\,]$ 
$\Rightarrow  \exists \, \kappa \in \mathrm{Card} :$
 $[\, R^{\lambda} {\cong}_R R^{[\kappa]}\,]$(4). Let $\mathfrak{m}\in \tmop{max}(R)\neq \emptyset$ 
 and  $\phi=\pi^R_\mathfrak{m}:R\rightarrow (R/ \mathfrak{m})$. Then $(R / \mathfrak{m})=F$ is a field, the r.h $[\,\phi$ is a f.p $]$(1), since $R$ is Noetherian.  Also  $[\,(R^{\lambda})^{\ast} {\cong}_R R^{[\lambda]}\,]$ (2) by assumption .
 So                         
 $(4)\Leftrightarrow [\, R^{\lambda} {\cong}_R R^{[\kappa]}\,]$
 $\underset{(1)}{\overset{\ref{rh2}.\mathrm{iii}}\Longrightarrow}$
$[\,F^{\lambda} \cong_F F^{[\kappa]}\,]\,(5)$. Also  
$(4)\Leftrightarrow [\, R^{\lambda} {\cong}_R R^{[\kappa]}\,]\Rightarrow$
$[\,(R^{\lambda})^{\ast} {\cong}_R (R^{[\kappa]})^{\ast}\,]$ $\overset{(2)}\Rightarrow$
$[\, R^{[\lambda]} {\cong}_R R^{\kappa}\,] \underset{(1)} {\overset{\ref{rh2}.\mathrm{iii}}\Longrightarrow}$    
$[\, F^{[\lambda]}\cong_F F^{\kappa}\,]\,(6)$. We  consider the $F$-vector space $V=F^{[\lambda]}$ and we get: $[\, V^{\ast}=(F^{[\lambda]})^{\ast}\cong_F F^{\lambda} {\overset{(5)}\cong}_F F^{[\kappa]}\,] $   
$\Rightarrow$
$[\, V^{\ast\ast }\cong_F (F^{[\kappa]})^{\ast} \cong_F F^{\kappa} {\overset{(6)}\cong}_F F^{[\lambda]}=V\,] \Rightarrow $ 
$ [\,V^{\ast\ast} \cong_F V \,]$ (7),
 which contradicts  Theorem $\ref{EKextension}\mathrm{.vi}$ since $\lambda=\tmop{dim}_K (V)$ is infinite by assumption. Therefore $[\, R^{\lambda} \notin \tmop{Free}(R)\,]$.
\end{proof}

\begin{lemma} \label{slendernofree} 	If $R$ is a slender Noetherian commutative ring then $R^{\omega}$ is not a free $R$-module.
\end{lemma}
\begin{proof}
By assumption $[\,R$  is slender, \underline{$R_{}$ is Noetherian} (2)  $]$ $\overset{\ref{reflexiveslender1}} \Rightarrow  [\,(R^{\omega})^{\ast} \cong_R R^{ [\omega]}\,]$     $\underset{(2)}{\overset{\ref{refnofree}}\Rightarrow}$
$[\, R^{\omega}\notin  \tmop{Free}(R) \,]$ .
\end{proof}

\begin{theorem}\label{slender10}
Let $R$ be a countable Noetherian commutative ring. If $R^{\omega} \in \tmop{Free}(R)$  then $R$ is Artinian.
\end{theorem}
\begin{proof} If  $R$ is non-Artinian then $[\, \exists \, \mathfrak{p} \in (\mbox{Spec}(R)\smallsetminus\mbox{max}(R))\neq \emptyset\,]$.
Let  $\phi=\pi^R_\mathfrak{p}:R\rightarrow \dfrac{R}{\mathfrak{p}}=S$. Then \\  $[\, \phi$ is f.p $]$(1)  since $R$ is Noetherian. By assumption  $[\, R^{\omega} \in \mathrm{Free}(R) \,]\Rightarrow$
 $[\, R^{\omega}  {\cong}_R R^{[\lambda]}$, for some 
 $\lambda \in \mathrm{Card}\,]$. So
$[\, R^{\omega} {\cong}_R R^{[\lambda]}\,] \underset{(1)}{\overset{\ref{rh2}.\mathrm{iii}} \Longrightarrow }$
$[\,  S^{\omega} \cong_S S^{[\lambda]}\,]$
  $\Rightarrow$   
 $[\,S^{\omega} \in \mathrm{Free}(S)\,]$ (2). Note that $S$ is a domain since  $\mathfrak{p} \in\tmop{Spec}(R)$,  $S$ is not a field since 
$ \mathfrak{p} \notin \tmop{max}(R)$ and $S$ is countable since $R$ is countable. Therefore $[\, S$ is a countable domain that is not a field
$]$ $\overset{\ref{slendercollection}\mathrm{.ix}} \Longrightarrow$
 $[\, S$  is slender$]$(3) $\overset{\ref{slendernofree}}\Rightarrow$  $[\,  S^{\omega} \notin \mathrm{Free}(S)\,]$, which contradicts (2). Therefore $R$ is Artinian.\end{proof}
The proof above motivates the definition below.
\begin{definition} \label{w-slenderdef} A  commutative ring $R$ is w-slender (weakly slender) iff there is a  r.h  $\phi:R\rightarrow S $ such that the ring  $S$  is commutative and  slender  and $\phi$
is f.p (i.e $S$ is a finitely presented $R$-module).
\end{definition}
We  have actually proven the next two theorems in  the above proof of 
\ref{slender10}
\begin{theorem} \label{countableslender2}
Any countable Noetherian commutative ring $R$  that is not  Artinian is  $\mathrm{w}$-slender.
\end{theorem}
\begin{proof}
We use the notation of the the proof of \ref{slender10} where $R_{\phi}$ is slender commutative ring and a f.p $R$-module by $\{ (3),(1) \}$, 
hence $R$ should be w-slender commutative ring.
\end{proof}

\begin{theorem}\label{slender12} Let  $R$ be a $\mathrm{w}$-slender commutative ring. Then $ R^{\omega} \notin \mathrm{Free}(R)$.
\end{theorem}
\begin{proof}
$[\, R$ be a w-slender$\,]$   $\overset{\ref{w-slenderdef}}\Rightarrow \{$ There is a r.h \space $\phi:R\rightarrow S$ such that  $[\,\phi$  is f.p $]$(1) and  
$[\, S$ is slender $]\,(3)$.
Assume that $[\, R^{\omega} \in \mathrm{Free}(R)\,]\Rightarrow$ 
$\{ \exists \,  \lambda \in \mathrm{Card} :[\,$
$R^{\omega} {\cong}_R R^{[\lambda]}\,]\} \underset{(1)}{\overset{\ref{rh2}.\mathrm{iii}} \Rightarrow }$
   $[\,S^{\omega} \cong_S S^{[\lambda]}\,]$   
 $\Rightarrow[\, S^{\omega} \in \mathrm{Free}(S)\,]$ (2).
But  $[\, S$  is slender $]\,(3)$ $\underset{\ref{slendernofree}}\Rightarrow$  $[\, S^{\omega} \notin \mathrm{Free}(S)\,]$, which contradicts (2). Therefore $[\,  R^{\omega} \notin \mathrm{Free}(R)\,]$.
\end{proof}

\begin{theorem} \label{wslenderaffinedomains2}Let $R$ be an affine algebra   over a field  $K$.\end{theorem}
\begin{enumerate} [label=(\roman*)] 
\item If $R$ is not Artinian  then $R$  is w-slender. \\
  $\mathbf{Proof:}$ Let $n=\mbox{dim}(R)\ge 0$. Then  
there is a r.h \space  $\phi:S=K[t_1,t_2,...,t_n]\rightarrow R$  that is 1-1 and finite by application of the \emph{ Noether Normalization Theorem}
(See Theorem 8.19 of \cite{kemper}). So $R=$ is a
f.p $S$-module since  $S$ is Noetherian (by Hilbert's basis Theorem) hence the r.h $[\, \phi$ is f.p $]$ (1). Also  by assumption   
  $[\, R$ is not Artinian ]  $\Rightarrow[\,  n=\tmop{dim}(R)\ge 1\,] \overset{\ref{slendercollection}\mathrm{.vi}}   \Rightarrow$ 
$[\, S=K[t_1,t_2,...,t_n]$ is slender]  $\overset{(1)}\Rightarrow \, [\,R$ is w-slender] .
 
\item If $R$ is a domain that is not a field then $R$  is slender. \\
$\mathbf{Proof:}$ Since any Artinian domain is a field and $R$ is not a field it follows that $R$ not Artinian hence by applying the above  construction the claim follows from $\ref{slendercollection}\mathrm{.vi},\, \ref{slendercollection}\mathrm{.ix}$
\end{enumerate}
 So the class of w-slender commutative rings is huge since it contains all non Artinian affine algebras (over fields) and all non Artinian countable Noetherian commutative rings by $\{\ref{wslenderaffinedomains2},\, \ref{countableslender2}  \}$.

\begin{theorem} \label{freeartinian}
$^{\ast}$Let $R$ be an affine algebra  over a field  $K.$ If $[\, R^{\omega}\in \mathrm{Free}(R)\,](1)$   then $R$ is Artinian. 
\end{theorem}
\begin{proof}
Assume that [$R$ is not  Artinian] $\overset{\ref{wslenderaffinedomains2}\mathrm{.i}}\Longrightarrow [\,  R$ is w-slender $]\overset{\ref{slender12}}\Longrightarrow $ 
 $[\, R^{\omega} \notin \mathrm{Free}(R)\,]$, which contradicts the assumption $(1)$. Therefore $R$ is an Artinian Ring.
\end{proof}

\begin{lemma} \label{noartinianisslender}Let $R$ be an Artinian  commutative ring:
\end{lemma}
\begin{enumerate} [label=(\roman*)]
\item If $[ \,R$ is field or  $R$ is local or $\mathrm{card}(R)\le \mathfrak{c} \,] $  then  $R^{\omega} \in  \mathrm{Free}(R)$. \\
$\mathbf{Proof:}$ Follows from  $ \ref{zornapplication2},\ref{collection}.\mathrm{iv},\ref{artinianfree}$ .  
\item If $ R^{\omega} \in  \mathrm{Free}(R)$ then $R$ is not w-slender. \\ $\mathbf{Proof:}$ This is logically equivallent to Theorem \ref{slender12} .
\item  If $[ \,R$ is field or  $R$ is local or $\tmop{card}(R)\le \mathfrak{c} \,] $  then   $R$ is not  w-slender. \\ $\mathbf{Proof:}$ 
Follows from the above.
\item $\mathrm{(a)}$ R is not slender and  $\mathrm{(b)}$ $R$ is not w-slender. So \underline {Artinian commutative rings are not w-slender}.\\
 $\mathbf{Proof:}$ $\mathrm{(a)}$:
Follows from Example 5.4 of \cite{jensentwo} since Artinian rings are perfect.\\
 $\mathrm{(b)}$: Assume that $R$ is w-slender then there is a r.h \space  $\phi : R\rightarrow S$  such that $S$ is slender and $\phi$ is f.g ($R$ is Noetherian) hence $S$ is Artinian as a finite extension of an Artinian commutative ring. So $S$ is slender and Artinian which contradicts 
 $\ref{noartinianisslender}.\mathrm{.iv.a}$. Hence $R$ is not w-slender.
\end{enumerate}
So far we have been able to establish the non-freeness of $R^{\omega}$ for various non-Artinian commutative rings but what about the non freeness of $R^X$
for an uncountable set $X$? 
In order to tackle this case we will use the concept of $\omega$-measurable cardinal that was introduced in $\ref{sets}$.

\begin{lemma} \label{measurablefacts}
Let $R$ be a slender commutative ring such that $\mathrm{card}(R)$ is not $\omega$-measurable and let 
$ \mu$ be a cardinal. Let also $\zeta=\mathrm{min}(\mathcal{L}_1)$ be the first non-$\omega$-measurable cardinal defined in $\ref{sets}$. Then $:$
\end{lemma}
\begin{enumerate}[label=(\roman*)]
\item $(R^{\mu})^{\ast}$ is a free $R$-module.
\item If $\mu$ is  not $ \omega$-measurable then 
$(R^{\mu})^{\ast} \underset{R} \cong R^{[\mu]}$  and both $R^{\mu},R^{[\mu]}$
are reflexive $R$-modules.
\item $\mathrm{(a)}$ If $\mu \geq \zeta $ then $\mu \in \mathcal{L}_1$  \space $\mathrm{(b)}$ \underline{ $(\mu < \zeta)\Leftrightarrow (\mu \in \mathcal{L}_1')$ } \space $\mathrm{(c)}$ The class $\mathcal{L}_1'$ downward closed.   
 \item  If $\mu$ is not $\omega$-measurable then   then $\mu^{\bullet}$ is not $\omega$-measurable.
\end{enumerate}

\begin{proof}
 $\mathrm{(i)}$ is  Corollary 3.7 of \cite{eklof},  $\mathrm{(ii)}$ is  Theorem 2.11  of \cite{eklof} .\\
 \underline{Proof of $\mathrm{(iii)}:$} See Theorem 2.11 of chapter II of $\cite{eklof}$ for a proof of $\mathrm{(a)}$. Also  $(b), \mathrm{(c)}$ follow  from $\mathrm{(a)}$. \\ 
 \underline{Proof of $\mathrm{(iv)}:$} It follows from
 $\{ \mathrm{(iii.b)},\, \ref{simple}   \}$  \end{proof}
\begin{lemma} \label{fpcardinal}
Let $M$ be a f.p  $R$-module over a ring
$R$ and $\phi: R\rightarrow S$ be a  r.h that is f.p. Then $:$ \\
$\mathrm{(i)}\, \mathrm{(a)}\,  M\cong_R (R^{[\kappa]} /L)$, for some $\kappa \in \mathrm{Card}_{\, < \omega}$ and $\, L\in \mathrm{Sub}_R R^n $. \quad\quad $\mathrm{(b)}\, |M|=| R^n /L |$\\
 $\mathrm{(ii)}\,\mathrm{(a)}$ If $ |R|< \omega $ then  $  |M| < \omega $  \quad\quad\quad\quad\quad\quad $\mathrm{(b)}$ If $ |R|\geq \omega  $ then  $  |M| \leq |R| $. \\
$\mathrm{(iii)}$ If  $|R|  \in \mathcal{L}_1'$ then   $[\, |M| \in \mathcal{L}_1'$ and  $ |S| \in \mathcal{L}_1' \,]$.\end{lemma} 
\begin{proof} $\mathrm{(i)},\,\mathrm{(ii.a)} $ are trivial and 
$\mathrm{(iii)} $ follows from $\{ \mathrm{(ii)},\,\ref{measurablefacts}\mathrm{.iii.c} \}$. \\
\underline{Proof of $\mathrm{(ii.b)}:$} 
 $[\, |R| \geq \omega \,]\,(1)\Rightarrow [\,  |M| \underset{\mathrm{(i.b)}}\le |R^n|\underset{(1)}=|R|\,]\Rightarrow$
$[\, |M| \leq |R|\,]$. \end{proof}

\begin{theorem}\label{measurable} 
Let $R$ be a $\mathrm{w}$-slender Noetherian commutative ring and $\lambda$ be an infinite cardinal such that  $[\, \lambda$ is countable or  both  $\tmop{card}(R),\,  \lambda$ are non $\omega$-measurable $]$. Then
$R^{\lambda} \notin  \tmop{Free}(R)$.
\end{theorem}
\begin{proof}  $R$ is w-slender $\overset{\ref{w-slenderdef}}\Rightarrow$ There is a r.h \space  $\phi:R\rightarrow S$ such that 
$[\,S$ is slender $]$(1) and  $[\, \phi$ is f.p $]$(2) .  
\underline {Assume that $R^{\lambda}\in \tmop{Free}(R) $}  $\Rightarrow$  $\{ \exists \, \kappa\in \mathrm{card}:$ 
  $[ \, R^{\lambda} {\cong}_R R^{[\kappa]} \, ]\} \underset{(2)}{\overset{\ref{rh2}.\mathrm{iii}}\Longrightarrow}$
$[\, S^{\lambda}\cong_S S^{[\kappa]}\,]\Rightarrow$  
$  [\,S^{\lambda}\in\tmop{Free}(S)\,]$ (3). It also follows by $\{ (2), \ref{fpcardinal}\mathrm{.iv} \}$  that $[\, \tmop{card}(S)$ is not $\omega$-measurable $]$(5). Therefore 
 $\{\ref{measurablefacts}.\mathrm{ii},  \ref{slendernofree}, (1),(5) \}\Rightarrow$  
$[\,(S^{\lambda})^{\ast}\cong_S (S^{[\lambda] }\,]$ (4).  So $\{(1), (4) \}\overset{ \ref{refnofree}} \Rightarrow$  $  [\,S^{\lambda} \notin \mathrm{Free}(S)\,]$, which contradicts (3).
Therefore $R^{\lambda} \notin \tmop{Free}(R)$.
\end{proof}

 Note that in  view of Theorem \ref{relief} and Lemma \ref{noartinianisslender} the restriction on cardinalities on the  above theorem  are redunant.  

\subsection{Free modules with isomorphic duals in weakly slender rings}

 \begin{theorem}  \label{secondtheorem} $^{\ast}$ 
 Let $R$ be a $\mathrm{w}$-slender  Noetherian commutative ring, 
  $\{\kappa,\,\lambda\}\subseteq \mathrm{Card}$,  $M \cong_R R^{[\kappa]}$ and $N \cong_R R^{[\lambda]}$. If both $\mathrm{card}(R), \, \kappa$ are not $\omega$-measurable then 
$[\,   M^{\ast} \cong_R N^{\ast}\,] \Rightarrow [\, M \cong_R N \,]$  
\end{theorem}
\begin{proof}  
 By assumption  $M \cong_R R^{[\kappa]}$ and  $N \cong_R R^{[\lambda]}$ hence $M^{\ast}\cong_R R^{\kappa}$ and $N^{\ast}\cong_R R^{\lambda}$.  There is   a r.h  \space
  $\psi :R\rightarrow S$ such that $[\,S$ is slender $]$(1) and 
  $[\, \psi$ is f.p $]$(2), since $[\, R$ is w-slender $\,]$(12). Also  $[ \, \tmop{card}(R)$ is not  $\omega$-measurable $]$(11) and $[\, \kappa$ is not $\omega$-measurable $]$(13), by assumption.  It also follows by $\{ (2), (11), \ref{fpcardinal}\mathrm{.iv} \}$  that $[\, \tmop{card}(S)$ is not $\omega$-measurable $]$(3).  
 By assumption, we have that $M^{\ast}\cong_R N^{\ast}$. 
 So \\ 
$[\, M^{\ast}\cong_R N^{\ast}\,] \Leftrightarrow$ 
$[\, R^{\kappa}\cong_R R^{\lambda}\,] \overset{\ref{rh2}\mathrm{.iii}}{\underset{(2)} \Longrightarrow} $
$[\, S^{\kappa} \cong_S S^{\lambda}\,] \Rightarrow$
$[\, (S^{\kappa} )^{\ast}\cong_S( S^{\lambda})^{\ast} \,]\,(4)$. 
 From $ \{  \ref{measurablefacts}\mathrm{.ii},\, (1),\,(3),\, (13)   \}$  that 
$[\, (S^{\kappa})^{\ast}\cong_S S^{[\kappa]}\,]$(5) and it follows from 
 $\{ (1),(3),\ref{measurablefacts}\mathrm{.i} \} $  that $(S^{\lambda})^{\ast }\in \mathrm{Free}(S)\Rightarrow$   $[\,(S^{\lambda})^{\ast }\cong_S S^{[\mu]}$,  for some  $\mu\in \mathrm{Card}\,]\,(6)$. So $\{(4),(5),(6)  \} \Rightarrow$          
$[\, S^{[\kappa]}\cong_S S^{[\mu]}\,] $  
$\overset{\ref{rankiswelldefined}}\Rightarrow [\,\kappa=\mu  \,]\,(7)$.
It also follows from  $\{ (6),\ref{dualfree} \}$ that  $[\,  \mu \ge \lambda\,] \,(8) $ $\overset{(7)}\Rightarrow$ 
$[\, \kappa \ge \lambda\,] \underset{(13)} {\overset{\ref{measurablefacts}\mathrm{.iii}}\Longrightarrow} [\, \lambda$ is  non-$\omega$-measurable  $]\,(9)$. It follows from   $ \{  \ref{measurablefacts}\mathrm{.ii},(1), (9)   \}$  that 
$[\, (S^{\lambda})^{\ast}\cong_S S^{[\lambda]}\,]$ (10). Therefore  
$\{ (4), (5), (10)\}\Rightarrow$
$ [\, S^{[\kappa]}\cong_S S^{[\lambda]}\,] \overset{\ref{rankiswelldefined}}\Rightarrow $
$[\, \kappa=\lambda\,] \Rightarrow $
$[\, R^{[\kappa]}\cong_R R^{[\lambda]}\,] \Rightarrow[\,  M \cong_R N \,]$. So we  have shown that $[\ (M^{\ast} \cong_R N^{\ast})\Rightarrow  (M \cong_R N) \,] $. 
\end{proof}

\begin{remark}Let $R$ be a Noetherian commutative ring. Then $:$\\
 $\mathrm{(i)}$ It follows from $\{\ref{secondtheorem},\, \ref{definitions2}\}$ that $\mathcal{L}^w_R$ if $R$ is w-slender.\\
$\mathrm{(ii)}$  If $R$ is a non-Artinian ring that is either countable or Hilbert then $R$ is w-slender hence $\ref{secondtheorem}$ is applicable to $R$.      \\
$\mathrm{(iii)}$ The cardinality restrictions placed on $\ref{secondtheorem}$ are very mild and will be later dropped in  Theorem $\ref{w-slender4}$.  \\
$\mathrm{(iv)}$ Combining the $\ref{secondtheorem}$ with  $\ref{ICFdual}$, we get the next theorem, which extends 
both.\end{remark}

\begin{theorem}  \label{secondtheorem2} 
Let $M\cong_R R^{[\kappa]}, N\cong_R R^{[\lambda]} $ be  free modules over a Noetherian commutative ring $R$ such that $\mathrm{(i)\,}[\,  \exists \,\mathfrak{m} \in \tmop{max(R)}: \big| \dfrac{R}{\mathfrak{m}} \big| \le \mathfrak{c}, \, \mathsf{ICF}\,]  $ \space  or \space
 $\mathrm{(ii)\,}[ \, |R|,\lambda $  are not $\omega$-measurable and  $R$ is $\mathrm{w}$-slender$\,]$.  Then $[\, M^{\ast} \cong_R N^{\ast}\,]\Rightarrow[\,  M \cong_R N\,] $
 \begin{proof}
 Case $\mathrm{(ii)}$ is Theorem \ref{secondtheorem}. \underline{Proof of case $\mathrm{(i)}$}  
 Let  $\phi=\pi^R_\mathfrak{m}:R\rightarrow \dfrac{R}{\mathfrak{m}}$.  
Then $\dfrac{R}{\mathfrak{m}}=F$ is a field  such that $[\, |F |\le \mathfrak{c} \,]\,(1)$  and  such that  the r.h $[\, \phi$ is a f.p $]$ (2) since $R$ is Noetherian. We consider the $F$-vector spaces $V=F^{[\kappa]},W=F^{[\lambda]}$. Note that $M^{\ast}\cong_R R^{\kappa}$,
$N^{\ast}\cong_R R^{\lambda}$ and $V^{\ast}\cong_F F^{\kappa}$, 
$W^{\ast}\cong_F F^{\lambda}$.  So \\
$[\, M^{\ast} \cong_R N^{\ast} \,]\Leftrightarrow  $
$[\, R^{\kappa} \cong_R R^{\lambda}\,] \underset{(2)}{\overset{\ref{rh2}.\mathrm{ii}}\Rightarrow}$
$[\, F^{\kappa} \cong_F F^{\lambda}\,]\Leftrightarrow$
$[\, V^{\ast}\cong_F W^{\ast}\,]\underset{(1)} {\overset{\ref{ICFdual}} \Leftrightarrow}$
$[\,  V\cong_F W\,] \Leftrightarrow$
$[\, F^{[\kappa]} \cong_F F^{[\lambda]}\,] \overset{\ref{rh}  } \Leftrightarrow$
$[\, \kappa=\lambda \, ]\Rightarrow$ 
$[\, R^{[\kappa]} \cong_R R^{[\lambda]}\,]\Leftrightarrow$
$[\, M\cong_R N\,]$. Hence
 $[\,(M^{\ast}\cong_R N^{\ast})\Rightarrow (M\cong_R N)\,]$.  
 \end{proof}
 
\end{theorem}
\begin{lemma} \label{dualfree}
Let  $R$ be any  commutative ring and  $\kappa,\lambda$ be any cardinals. If
$(R^{\kappa} )^{\ast}\cong_R R^{[\lambda]}$ then $\kappa \le \lambda $.
\end{lemma}
 \begin{proof}
The $R$-module $L=R^{[\kappa]}$  is torsionless  $\Rightarrow$[the map 
$\sigma_L:L\rightarrow L^{\ast\ast}$ is 1-1$]$ 
  $\Rightarrow [\,  L \underset{R}\rightarrowtail L^{\ast\ast}\,]$ (1). Also
  $[\, L^{\ast\ast}\underset{R}\cong (R^{\kappa} )^{\ast}\underset{R}\cong R^{[\lambda]} \,]$  $\overset{(1)}\Rightarrow$ 
  $[\, R^{[\kappa]}\underset{R}\rightarrowtail R^{[\lambda]}\,]$
$\overset{\ref{innocuous}}\Rightarrow [\, \kappa\le \lambda\,]$ .
\end{proof}

\begin{lemma} \label{innocuous}
Let $R$ be a commutative ring and $\kappa,\lambda$ be cardinals. If $R^{[\kappa]}\underset{R}\rightarrowtail  R^{[\lambda]} $ then $\kappa  \le  \lambda$ .
\end{lemma} 
 \begin{proof}
 This innocuous looking statement is actually  hard to prove. We will  need to invoke  Theorem \ref{Lazarus} of Lazarus  to do so. Fortunately the proof of \ref{Lazarus}  will lead us to Theorem \ref{unexpected} .\\
 \underline{Case 1: $R_{}$ is a domain.} Consider the r.h  $\phi=\mathrm{i}_R: R\rightarrow \mathrm{Q}(R)=F$. So $[\, F$ is a field$\,]\,(1)$ and $[\, \phi$ is flat$\,]\, (2)$. By assumption 
  $[\, R^{[\kappa}] \underset{R} \rightarrowtail  R^{[\lambda]}\,] \overset{\ref{rh2}\mathrm{.v.b}}{\underset{(2)}\Longrightarrow} $     
$[\, F^{[\kappa]} \underset{F}\rightarrowtail F^{[\lambda]}\,]\overset{(1)}\Rightarrow$     
$[\, \dim_F (F^{[\kappa]}) \le \dim_F (F^{[\lambda]})\,] \Rightarrow$  
$[\, \kappa \le \lambda \,]$. \\
\underline{Case 2:  $\lambda_{}$ is infinite}. By assumption  $R^{[\kappa]}\underset{R}\rightarrowtail  R^{[\lambda]} \Rightarrow$
 $[\, \exists\, f \in \mathrm{Hom}_R( R^{[\kappa]}, R^{[\lambda]}) $ such that $f$ is  1-1$\,]\,(3)$. Therefore the set 
  $S=\{ f(e_x)\,|\,x\in \kappa \}$ is an $R$-linearly independent subset of $R^{[\lambda]}$ of cardinality $\big|S \big|=\kappa$         that is contained in a maximal one  $T$ (hence $\kappa=\big| S\big| \le \big| T\big|\,(1)$) by application of \ref{zornapplication1}. The set   $W=\{ f(e_x)|x\in \lambda \}$ is an $R$-linearly independent subset of $R^{[\lambda]}$  of cardinality $ \big|W \big|=\lambda$ that is a maximal one. It follows from the  Case 2 of the  Theorem of Lazarus \ref{Lazarus} that
   $\big|T \big|=\big| W \big|=\lambda \overset{(1)} \Rightarrow$
   $\kappa \le \lambda$. \\
\underline{Case 3:  $\lambda_{}$ is finite.} 
It is trivial that the set  $S=\{ f(e_x)|x\in \kappa \}$ is an
 $R$-linearly independent subset of $R^{[\lambda]}$ of cardinality $\big|S \big|=\kappa$ hence     
$\kappa=\big|S \big|\le \lambda$  by application of  of $\ref{Lazarus}\mathrm{.i}$. Therefore  $\kappa\le \lambda$.  
 \end{proof}
 The next lemma is much easier to prove.
\begin{lemma} \label{noninnocuous}
Let $R$ be a commutative ring and $\kappa,\lambda$ be cardinals. If $R^{[\kappa]}\underset{R} \twoheadrightarrow  R^{[\lambda]} $ then $\kappa  \geq  \lambda$ \end{lemma} 
\begin{proof} Let $\mathfrak{m}\in \tmop{max}(R)\neq \emptyset$ 
 and  $\phi=\pi^R_\mathfrak{m}:R\rightarrow \dfrac{R}{\mathfrak{m}}=F$  
. Then $[\, F$ is a field$\,]\,(1)$. By assumption, we have   $[\, R^{[\kappa]}\underset{R} \twoheadrightarrow  R^{[\lambda]}\,] \overset{\ref{rh2}\mathrm{.v.a}}\Longrightarrow$
 $[\, F^{[\kappa]}\underset{F} \twoheadrightarrow  
 F^{[\lambda]}\,]\overset{(1)}\Rightarrow$  
 $[\, \dim_F (F^{[\kappa]}) \geq \dim_F (F^{[\lambda]})\,] \Rightarrow[\, \kappa \geq \lambda \,]$.                
 \end{proof} 
 
\begin{lemma} \label{zorn}
$( \emph{Zorn's \, Lemma}) $ Let $(P,\le)$ be a poset such that any chain $($i.e totally ordered subset $)$ of  $P$ has an upper bound  and let $a\in P$. Then $\mathrm{(i)}\, \mathrm{Max}(P )\neq \emptyset$  \quad $\mathrm{(ii)}\, \exists \, b \in \tmop{Max}(P) :b\ge a$. 
\end{lemma}

 \begin{lemma} \label{zornapplication1}Any $R$-linearly independent subset  $S$ of an $R$-module  $M$ 
over a commutative ring  $R$ is contained in a maximal one 
$(R$-linearly independent subset of $M)$
\end{lemma} 
 \begin{proof}
 Let $P=\tmop{Ind}_R(M)\overset{\text{def}}=\{T\in \mathcal{P}(M): T \text{ \, is}\,$   linearly independent over $R \}$. Consider the order relation $\le$ on $P$ defined by
 $C\le D\overset{\text{def}}\Leftrightarrow C\subseteq D.$  Let $[\, \Lambda\subseteq P$ be totally ordered subset of $P\,]$(1) and let 
  $[\, \hat{\Lambda}=\cup \Lambda\,]$(2).  It is trivial that  $\hat{\Lambda}$ is an upper bound of $\Lambda$ in $\mathcal{P}(M),$ so it remains to show that   $\hat{\Lambda}\in P=\tmop{Ind}_R(M)$. \underline{Assume that $ \hat{\Lambda}\notin P_{} $ } $\Rightarrow [\, \exists \, A\in \mathcal{P}_{\text{fin}}(\hat{\Lambda}) : \underline{A\notin P} \,\, (4) \, ] $, hence $[\, A$ is a finite set $]$ (5). But $A\subseteq  \hat{\Lambda}=\cup \Lambda\overset{(5)}\Rightarrow$         
 $[\, A\subseteq \overset{n} {\underset{k=1} \cup} A_k \,(6),$ for  $A_k\in \Lambda\, ].$ Also $(1)\Rightarrow$ $\exists\,  m\in \mathrm{T_n}:A_m=\max \{A_1,A_2\dots A_n \}$ $\Rightarrow$   
$[\,  \overset{n} {\underset{k=1} \cup} A_k=A_m \,] \,(7)$                                                      
 $\overset{(6)}\Rightarrow[\, A\subseteq A_m\in \Lambda\subseteq P \, ]\Rightarrow$
 $[\,  A\in P \,](8)$ (since subsets of linearly independent sets are linearly  independent). So we reached a contradiction by $\{(4),(8)  \}$. Therefore $[$ every chain of $P$ has an upper bound (in $P$)$]$(9). Let now $S\in P$  be $R$-linear independent subset of $M$.
 Then it follows by $\{(9),\ref{zorn} \} $ that $S$ is contained in a maximal $R$-linear independent subset of $M.$ 
\end{proof}

The next lemma is important on what follows, is very simple to prove and  it can be found in \cite {lazarus}.
\begin{lemma} \label{torsion}
Let $R$ be  a commutative ring $R$, $M\in$R-$\mathrm{Mod}$, $S\in \tmop{Ind}_R(M)$  and $L=<S>_{R}$. Then  \\ $S\in \tmop{Max}(\tmop{Ind}_R(M))$ $\Leftrightarrow \{ \, \forall \, m\in M \, [ \,\exists \,  r\in (R\smallsetminus \mathbf{0}) : r\cdot m\in L \, ]\, \}\Leftrightarrow $ $\dfrac{M}{L}$ is a torsion 
$R$-module.  \end{lemma}

\begin{lemma} \label{zornapplication2}
Any $F$-linearly independent subset  $S$ of a vector space  $V$ 
over a field  $F$ can be extended to an $F$ basis of $V$. In particular any vector space has a basis.\end{lemma} 

\begin{proof}  Follows from \ref{zornapplication1}, \ref{torsion} since
any torsion $F$ vector space $V=\dfrac{M}{L}$  is zero. 
\end{proof} 
 For the next theorem, whose proof can be found in  \cite{blocki}, $\mathcal{I}^R_k(A)$ is the $k$-th determinantal ideal of $A$ i.e  the ideal generated by all $k\times k$ minors of $A$.
 \begin{theorem} \label{mccoy} $($McCoy$)$
   Let  $\mathfrak{T}^R_{A} : R^m \rightarrow R^n$ 
  defined by $\mathfrak{T}^R_{A} (m) = A \cdot m$, for $A\in R^{n\times m }$, where $R$ is commutative ring . Then  \space
  $\mathfrak{T}^R_{A}$ is 1-1  $\Leftrightarrow $
  $\{ m \le n, \,\tmop{ann}(\mathcal{I}^R_m(A))=\mathbf{0 } \}$ 
   \end{theorem}

\begin{theorem} \label{Lazarus}
 $(  \emph{Lazarus} )$  Let $M\cong R^{[X]}$ be a free 
  $R$-module over a commutative ring $R$  and $S\subseteq M$.
 
\end{theorem}

\begin{enumerate}[label=(\roman*)]
\item $\mathbf{Claim \,1:}$ If $X$ finite and  $S\in \tmop{Ind}_R(M)$
  then $ \big| S \big| \le \big| X \big| $.(Hence $S$ is finite)\\
$\mathbf{Proof:}$ Let $\big| X \big|=n$. We may assume, without loss of generality, that $M=R^n$ and 
$S\subseteq R^n$. Let     $T \in \mathcal{P}_{\text{fin}}(S)$    then  $[\,T\in \mathrm{Ind}_R(M)\,]\,(21)$. Let $m=\big| T \big|$ ,   $A=[s_1|s_2|\dots|s_m]\in R^{n\times m}$ and  $\mathfrak{T}^R_{A} : R^m \rightarrow R^n$
  defined by $\mathfrak{T}^R_{A} (m) = A \cdot m.$  
  Let      $x=(x_1,x_2\dots,x_m)^{\top}\in R^m$    then  $[\, A\cdot x=0\,]\Leftrightarrow$
$[\,  \underset{k=1}{ \overset{m}  \Sigma} x_k \cdot s_k =0\,]$
 $\overset{(21)}\Leftrightarrow [\, x=0 \,]$.
Hence $\mathfrak{T}^R_{A}$ is 1-1 
$\overset{\ref{mccoy}}\Leftrightarrow$
$\{m \le n, \,\mathrm{ann}(\mathcal{I}^R_m(A)=\mathbf{0}  \}$ 
$\Rightarrow[\,  m\le n\,]\Leftrightarrow $
$\big| T \big| \le \big| X \big|.$ Hence $\{[\, \big| T \big| \le \big| X \big| \, ] \, \forall \, T  \in \mathcal{P}_{\text{fin}}(S) \} \Rightarrow$
$\big| S \big| \le \big| X \big|.$
 \item $\mathbf{Claim \,2:}$ If $X$ is infinite and    $S\in \tmop{Max}(\tmop{Ind}_R(M))$ then $ \big| S \big| = \big| X \big| $.     \\
 $\mathbf{Proof:}$
\underline{Assume that  $\big| S \big| >  \big|X \big|$} (2).  Let   $H:S\rightarrow \mathcal{P}_{ \text{fin}}(X)$ defined by 
$H(f)=\tmop{Supp}(f)$. We  will show that the above map  has an uncountable fiber. If not then 
\underline{$[\, \big|  H^{-1}(\{ a \}) \big|\le \big|\mathbb{N} \big|\,]\, \forall a\in B \,$} (3), where $B=\mathcal{P}_{ \text{fin}}(X)$.  Note $\big| \mathcal{P}_{ \text{fin}}(T) \big|=\big| X \big|\,(4)$, since $X$ is infinite. \\
$S=\underset{a\in B}\cup H^{-1}(\{ a \})\Rightarrow$
$\big| S\big| \le\underset{a\in B} \sum \big|  H^{-1}(\{ a \}) \big|$
$\underset{(3) }\le \big| B \big|\cdot \big|\mathbb{N} \big|=$
$\big| \mathcal{P}_{ \text{fin}}(X) \big| \cdot \big|\mathbb{N} \big|\overset{(4)}= \big| X \big|\cdot\big|\mathbb{N} \big|=\max \{\big| X \big|,\big|\mathbb{N} \big| \} $
$=\big| X \big|\Rightarrow$
$[\, \big| S\big| \le \big| X \big| \,]$, which contradicts (2).
Therefore $H$ has an infinite  (as uncountable) fiber, i.e $\{\exists \, A \in  \mathcal{P}_{ \text{fin}}(X) :[\, T=H^{-1}(A)\,  \text{is infinite} \,] \}$. So  $[\, T\subseteq R^A \,]$ (5), $[\, T$ is infinite $]$(6), and $[\, A$ is finite $]$(7). Also     
  $[\,T\in \tmop{Ind}_R(M)\,]$(8) since $[\, S\in \tmop{Ind}_R(M) \,] (9)$   and $T\subseteq S$. 
Therefore we reached a contradiction  from  $\{(5), (6),(7),(8), \text{Claim 1}  \}$
Hence $\boxed{\big| S \big| \le  \big|X \big|}$ (1). Note that we have not used the maximality of $S$ yet. It is trivial that $[\,\underset{f\in S}\bigcup \tmop{Supp}(f)\subseteq X \,]$ (11). If equality  does not hold  in (11) then $\{ \exists \, x  \in X :[x\notin \tmop{Supp}(f)\,] \, \forall f \in X   \}$ (12) $\Rightarrow e_x\notin X$  $ \Rightarrow [\, (S\cup \{e_x \}) \notin \tmop{Ind}_R(M)$ $] \,$(13) (as it follows from  the maximality of $S$ ) 
$\Rightarrow\{ \exists \, r=(r_0,r_1,r_2\cdots,r_n)\in (R^{n+1}\smallsetminus \{0\}):$
$ [\, r_0\cdot e_x + \underset{k=1}{ \overset{n} 
\Sigma} r_k\cdot s_k =0 \,] \}$ (14)
  $ \overset{(13)}\Rightarrow  $ 
   $[\,  r_0=0\, ] $      $\overset{(14)}\Rightarrow$
$[\, \underset{k=1}{ \overset{n}  \Sigma} r_k \cdot s_k =0 \,]$
$\Rightarrow $
$(r_1,r_2\cdots,r_n)=0 \Rightarrow [\, r=0 \,]$, which is contradiction since $r\neq 0$, by construction. Hence
 $[\, \underset{f\in S}\bigcup \tmop{Supp}(f)= X \,] $ (15). For $f\in S, $ $\tmop{Supp}(f)$ is finite hence  $\big| \tmop{Supp}(f)\big|\le \big| \mathbb{N}\big| $ (16). If $S$ was finite then it would follow from $\{(15),(16)\}$ that $X$ is finite, which would contradict the assumption. Hence  $\big| S \big|\ge \big| \mathbb{N}\big| $ (17). Finally we get :\\ 
 $(15)\Leftrightarrow$  $[\, \underset{f\in S}\bigcup \tmop{supp}(f)= X \,] $
 $\Rightarrow$
 $\big|  X\big| = \big| \underset{f\in S}\bigcup \tmop{Supp}(f)  \big|$
 $\le \underset{f\in S}\Sigma \big| \tmop{Supp}(f)\big| $  
$  {\overset{(16)} \le }  \underset{f\in S}\Sigma \big| \mathbb{N}\big|$
 $\le  \big| S\big| \cdot  \big| \mathbb{N}\big|$
 $=\max\{  \big| S\big|,\big| \mathbb{N}\big| \}\overset{(17)}= \big| S\big| \Rightarrow \boxed{ \big|  X\big| \le \big|  S \big|}$ (18). It follows from $\{(1),(18)\}$ that  $\boxed{ \big|  X\big| = \big|  S \big|}$

\item $\mathbf{Claim \,3:}$ If $R$ is Noetherian,  $X$ is finite and  $S\in \tmop{Max}(\tmop{Ind}_R(M))$  then $ \big| S \big| = \big| X \big| $. \\  
 $\mathbf{Proof:}$  This is Proposition 3.3 of \cite{lazarus} and is not going to be used here, except in claim 4.
 
\item $\mathbf{Claim \,4:}$  If $R$ is Noetherian, $S\in \tmop{Max}(\tmop{Ind}_R(M))$  and $X$ is any set then   $ \big| S \big| = \big| X \big|$. \\ $\mathbf{Proof:}$ Follows from $\mathrm{ii,iii}$.
\end{enumerate}

\begin{theorem} \label{unexpected} 
Let $M$  be a free $R$-module over  a \underline{local Artinian and self injective} commutative ring $R$. Then  any  $R$-linearly  independent subset of $M$  can be  extended to an $R$-basis of $ M$.
\end{theorem} 
\begin{proof} Let $S\in \tmop{Ind}_R(M)$  
$\overset{\ref{zornapplication1}}\Rightarrow$
$[\, \exists \, T \in \tmop{Max}( \tmop{Ind}_R(M))  :$
 $S\subseteq T \subseteq M\,]$.  In order to  show that $T$ is an
  $R$-basis  of $M$ it remains to show that $<T>_R=M$. Consider the 
   $R$-linear map $f:R^{T}\rightarrow M$  defined by 
$f(\overset{n}{\underset{k=1}\Sigma}r_k\cdot e_{x_k})=\overset{n}{\underset{k=1}\Sigma}r_k\cdot x_k$. The map $f$ is 1-1 since $T\in \tmop{Ind}_R(M) $ and $<T>_R=\tmop{Im}(f)$. \underline{Assume that $T_{}$ is not an $R$-basis of $M_{} $} then $f$ is not onto  hence $ \tmop{Cok}(f)\neq \mathbf{0}\,(1).$ The $R$-module $R^{[\lambda]}=\underset{x\in \lambda}\Sigma R\cdot e_x\underset{R}\cong\underset{x\in \lambda} \Sigma R $ is $R$-injective,  as a direct sum of injective modules, since $R$ is self injective over a Noetherian ring $R$, by application of Proposition 18.13 of \cite{fuller}. The short exact sequence 
 $\mathbf{0}\rightarrow R^{[\lambda]} \overset{f}\rightarrow M\rightarrow \tmop{Cok}(f)\rightarrow \mathbf{0}$ splits since 
 $R^{[\lambda]}\in \mathrm{Inj}(R)$ hence $[\, M\underset{R}\cong  R^{[\lambda]} \oplus  \tmop{Cok}(f)\,] (3)\Rightarrow $ 
$\tmop{Cok}(f)\in \tmop{Proj}(R)\overset{\ref{collection}.\mathrm{i} } = \tmop{Free}(R)\Rightarrow $ 
$[\, \mathrm{Cok}(f)\in \mathrm{Free}(R)\,](4)$. Also it follows by \ref{torsion}
that  $[\, \mathrm{Cok}(f)$ is $R$-torsion ] (5), since
$T \in \mathrm{Max}( \mathrm{Ind}_R(M))$. Therefore we have reached a contradiction by $\{ (4),(5) \}$, since non-trivial free modules are never torsion. Hence $T$ is an $R$-basis of $M$. Therefore any $R$-linearly independent subset $S$ of $M$ can be extended (i.e it is contained) to an $R$-basis $T$ of $M$.
\end{proof}

\begin{definition}\label{def.steinitz.bad} A commutative ring $R$  is Steinitz iff any $R$-linear independent subset of a free $R$-module $M$ extends to an  
$R$-basis of $M$.
\end{definition}
\begin{theorem} Let  $R$ be a local commutative ring. If $R$ is Artinian the $R$  is Steinitz.
\end{theorem}
\begin{proof}It Theorem 1.1  of \cite{nichols}.\end{proof}
\begin{lemma} \label{steinitz7}
Let $R$ be a   ring that is Artinian and local, $\{ P, \, Q  \} \subseteq \mathrm{Free}(R)$ and  $f\in \mathrm{Hom}_R(P,Q)$. Then 
$\mathrm{(i)\, }\, [\,f$ is $1$-$1\,]\,(1) \Leftrightarrow$  $[\, f$ is $R$-left invertible$\,]\,(2)$. \quad \quad\quad   $\mathrm{(ii)\, }\, [\,f$ is onto$\,]\,(3)  \Leftrightarrow$  $[\, f$ is $R$-right invertible$\,]\,(4)$.\end{lemma}
\begin{proof}  \ \underline{$\mathrm{Proof\,of}$  $\mathrm{(i)}$}: The direction $(2)\Rightarrow (1)$ holds trivially.  
\quad\quad   
$\cdot \,\mathrm{Proof\,of} \, (1)\Rightarrow(2) :$
 It follows from  Theorem \ref{steinitz7} that $[\, R$ is Steinitz$]\,(5)$. Let $[\, B_1 \in \mathrm{Bas}_R(P)\neq \emptyset\,]\,(7)$ and $[\,S=f(B_1)\subseteq Q\,]\,(8)$. Then $\{(1),\, (7),\, (8)  \}\Rightarrow$  $[\, S\in \mathrm{Ind}_R(Q)\,]\overset{(5)}\Rightarrow$ 
$[\, S\subseteq B_2$, for some $ B_2 \in \mathrm{Bas}_R(Q)\,]\,(9)$. There is a factorization $[\, f=h\circ f'\,]\,(10)$, where $f':P\rightarrow \langle S \rangle_R$  and
$h:\langle S \rangle_R\rightarrow \langle B_2 \rangle_R=Q$ is  the inclusion homomorphism. Note that $[\, h$ is $R$-left invertible$\,]\,(11)$ since
$[\,\{S,\, B_2\} \subseteq \mathrm{Ind}_R(Q),\, S\subseteq B_2\,] $ 
and trivially 
 $[\, f'$ is $R$-isomorphism $\,]\,(12)$ since $f$ is 
$1$-$1$. It follows from $\{ (10),\,(11),\,(12)$ that $f$ is $R$-left invertible.\\
 \underline{$\mathrm{Proof\,of}\, \mathrm{(ii)}$}:
 $[\,(4)\Rightarrow (3)\,]$ holds trivially and   
$[\, (3)\Rightarrow(4)\,] $ follows from the projectivity of $Q$.

\end{proof}

\begin{remark}  $\mathrm{(i) }$ Theorem $\ref{unexpected}$ is an extension of  Theorem \ref{zornapplication2} and a special case of of Theorem $ \ref{steinitz7}\mathrm{.i}$.\\ $\mathrm{(ii) }$   
Also note that proving $R^{\lambda}\in \tmop{Proj}(R)$, when $R$ is QF, is simpler than \ref{collection}.$\mathrm{ii}$ since 
$\tmop{Inj}(R)=\tmop{Proj}(R)$( by Theorem 15.9 of \cite{lam}) and  $R^{\lambda}\in \tmop{Inj}(R)$
(product of injectives is always injective) 
\end{remark}

\begin{theorem} \label{Bass}
Let $R$ be a connected Noetherian commutative ring and $M$ be a \underline{non f.g} $R$-module.\\ If 
 $M\in \tmop{Proj}(R)$ then  $M\in \tmop{Free}(R)$. 
\end{theorem}
\begin{proof}
This is a  famous Theorem of Bass. It is Corollary 4.5 of \cite{ bass2}. Note that, in this article , it will be applied to $M=R^{\lambda}$, for an infinite cardinal $\lambda$, which is non f.g by \ref{lost}.
\end{proof}

\begin{lemma} \label{remainprojective} Let $R,S$ be commutative rings and $\phi:R\rightarrow S$ be r.h. Then  \\
 $[\,P\in \tmop{Proj}(R)\,]\Rightarrow [\, P_{\phi}\in \tmop{Proj}(S)\,]$ i.e   extensions of projective modules remain projective.\end{lemma}
\begin{proof}
$[\, P\in \mathrm{Proj}(R)\,] \Rightarrow \{ [\,\exists \,  \lambda \in \mathrm{Card}\,]\,[ \exists \, Q \in \mathrm{Proj}(R)\, ]:[P \oplus Q  \cong_R R^{[\lambda]}\, ] \}\Rightarrow$
$[\,(P\oplus Q)_{\phi} \cong_{R_{\phi}} (R^{[\lambda]})_{\phi}\,] \overset{\ref{rh}\mathrm{.i}} \Longrightarrow  $ 
$[\, (P_{\phi}\oplus Q_{\phi})  \cong_{R_{\phi}} (R_{\phi})^{[\lambda]}\,]$
$\Rightarrow [\,  P_{\phi}\in \tmop{Proj}(R_{\phi})\,]$
$\Leftrightarrow  [\,  P_{\phi}\in \tmop{Proj}(S)\,]$ \end{proof}

\subsection{Projectivity of the dual and half-slender rings}

\begin{definition} \label{hilbertdef}
Let $R$ be a commutative . Then $R$ is Hilbert  $\overset{\text{def}}\Leftrightarrow$ Any prime ideal of $R$ is the intersection of a set of maximal ideals of $R\Leftrightarrow$
$\forall\, \mathfrak{p}\in \mathrm{Spec}(R)\, \exists\,  \Delta \in \mathcal{P}(\mathrm{max}(R)):[\,  \mathfrak{p} =\cap \Delta\, ]$.
\end{definition}
 On next theorem, we will extend the class of w-slender commutative rings to include all non-Artinian Hilbert commutative rings. 

\begin{theorem} \label{hilbert} Let $R$ be Noetherian commutative 
ring. Then:
\end{theorem}
\begin{enumerate} [label=(\roman*)]
\item  If $R$ Artinian then  $R$ is Hilbert. \\
 $\mathbf{Proof:}$ It is trivial since $\tmop{Spec}(R)=\tmop{max}(R)$.
\item If $R$ is a domain such that  $\tmop{dim}(R)=1$ then 
$[\, R$ is Hilbert $\Leftrightarrow \mathfrak{J}(R)=\mathbf{0} \,]$. So $\mathbb{Z}$ is  Hilbert. \\ $\mathbf{Proof:}$ 
$[\, R$ is a domain,  $\tmop{dim}(R)=1\,]\Rightarrow$
$[\, \tmop{Spec}(R)=\tmop{max}(R)\cup \mathbf{0} \,]$ (9). \\ 
Assume that $[\, \mathbf{0}=\mathfrak{J}(R)\overset{\text{def}}=\cap \tmop{max}(R)\,] $ (10)  Let $\mathfrak{p}\in \tmop{Spec}(R)$. If $\mathfrak{p}\in \tmop{max}(R)$ then we may take $\Delta=\{\mathfrak{p} \}\subseteq \tmop{max}(R)$. If $\mathfrak{p}\notin \tmop{max}(R)\overset{(9)}\Rightarrow \mathfrak{p}=\mathbf{0} \overset{(10)} =\cap \tmop{max}(R)$, so   we may choose $\Delta=\tmop{max}(R)$. Hence we proved that 
$[\, R$ is Hilbert $\Leftarrow \mathfrak{J}(R)=\mathbf{0} \,]$. For the other direction  $[\, R$ is Hilbert $\Rightarrow \mathfrak{J}(R)=\mathbf{0} \,]$: \\
$[\,R$ is a domain$\,] \Rightarrow[\,  \mathbf{0} \in \tmop{Spec}(R)\,] 
\overset{R \text{\, is \, Hilbert} } \Rightarrow$  $\{\,[\, \exists 
\Delta \subseteq \tmop{max}(R) \,] \,(11):[\mathbf{0}  \overset{\ref{hilbertdef}}= \cap \Delta $
$ \underset{(11)}\supseteq  \cap  \tmop{max}(R) \overset{\text{def}} =   \mathfrak{J}(R) \,]\, \}\Rightarrow$
$[\, \mathbf{0}\supseteq \mathfrak{J}(R)\,] \Rightarrow [\,  \mathfrak{J}(R)=\mathbf{0}\,]$.
 \item $R$ is Hilbert  $\Rightarrow$   $\dfrac{R}{I}$
 is Hilbert, for any proper ideal $I$ of  $R$.\\ $\mathbf{Proof:}$
 Trivial. See \cite{kaplansky} .
 \item $R$ is Hilbert  $\Leftrightarrow$  The polynomial ring $R[t]$ is Hilbert. \\ $\mathbf{Proof:}$ Theorem 31 of \cite{kaplansky}.
  \item $R$ is Hilbert $\Leftrightarrow$  The polynomial ring $R[t_1,t_2\cdots,t_n]$ is Hilbert. \\ $\mathbf{Proof:}$ Folllows from
  $\mathrm{(iv)}$.
  \item $R$ is Hilbert $\Rightarrow$ Any affine $R$-algebra $S=\dfrac{R[t_1,t_2\cdots,t_n]}{I}$ is Hilbert. \\ $\mathbf{Proof:}$ Folllows from $\{\mathrm{iii},\mathrm{v}   \}$.
 
  \item Assume that $R$ is local. Then $[\, R$ is Hilbert $\Leftrightarrow$
   $R$ is Artinian $]$. \\ $\mathbf{Proof:}$ Simple.
   \item  $[\, R$  is domain, $R$ is non-local, $\tmop{dim}(R)=1 \,] \Rightarrow [\, R$  is slender$\,]$ \\  $\mathbf{Proof:}$ This is Lemma 5.3 of  \cite{jensen}
   \item $[\, R$ is Hilbert, $R$ is non-Artinian  $ \,]\Rightarrow$
   $[\,R$ is w-slender$]$. \\    $\mathbf{Proof:}$
    $[\, R$ is non-Artinian$\,]$ $\Rightarrow [\, \exists\,  \mathfrak{p}\in \tmop{Spec}(R):$
    $\tmop{dim}(\, \dfrac{R}{\mathfrak{p}}\,)=1\,]$. Let  $\phi=\pi^R_\mathfrak{p}:R\rightarrow \dfrac{R}{\mathfrak{p}}$.  
Then 
  $[\, \phi$ is f.p ] (1)  since $R$ is Noetherian and $[\, \dfrac{R}{\mathfrak{p}}=$ $S$ is a domain] (2)   such  that  
$[\, \tmop{dim } (S)=1 \,]$ (3).
Also $[\, S$ is Hilbert ](4)  by $\mathrm{(iii)}$ since $R$ is so. Hence $\{ (3),(4)  \} \overset{\mathrm{(vii)}}\Rightarrow$ $[\, S$ is non-local ] (5) . Therefore 
 $\{ (2),(3),(5)  \} \overset{\mathrm{(viii)}}\Rightarrow $  
  $[\, S$ is slender ] $\overset{(1)}\Rightarrow[\,  R$  is
   w-slender]. 
   \item Assume that \underline{$R_{}$ is Hilbert}. 
    If $R^\mathbb{N} \in \tmop{Free}(R)$  then $R$ is Artinian ring.\\
    $\mathbf{Proof:}$ Assume that $[\, R$ is non-Artinian ]$\overset{\mathrm{(ix)}}\Rightarrow [\, R$ is w-slender ] $\overset{\ref{slender12} } \Rightarrow$
    $[\, R^\mathbb{N} \notin \tmop{Free}(R)\,]$, which  contradicts the assumption. Therefore $R$ is Artinian.
   \item $^{\ast}$ Assume that \underline{$R_{}$ is Hilbert}. Then 
   $\boxed{R$ is non
    Artinian $\Leftrightarrow  R$ is w-slender $}$ \\
    $\mathbf{Proof:}$ Folllows from $\{ \ref{hilbert}\mathrm{x},\ref{noartinianisslender}\mathrm{.iv.b}\}$. 
  
\end{enumerate}

\begin{definition}
A commutative ring $R$ has bounded type $\Leftrightarrow$ either 
$\mathrm{(i)}$ $|R| < \mathfrak{c}$ or $\mathrm{(ii)}$ $\{ \,R$ is an algebra over a field $K$ and $ \mathrm{dim}_KR<|K|^{\omega}\, \} $ or 
$\mathrm{(iii)}$  $\{ \, R$ is an algebra over an Artinian commutative ring $T$  such that   $[ \, {\tmop{dim}}_{\frac{T}{\mathfrak{m}}}( \, \dfrac{R}{\mathfrak{m}R})<\big| \dfrac{T}{\mathfrak{m}} \big|^{\, \omega} \,\,]$
$\forall \, \mathfrak{m}\in \tmop{max}(T) \, \}$. Note that this is  the Definition 3.1 of of \cite{jensentwo}.
\end{definition}
Next, we introduce and apply  the well-known concept of half-slenderness.

\begin{definition} \label{def.half.slender} Let $R$ be a  ring, 
$\kappa\in\mathrm{Card}$ and $\Phi^{\kappa}_{R}:R^{[\kappa]}\rightarrow R^{\kappa}$ be the  natural inclusion. \\
$\mathrm{(i)}\, R^{\cdot k}=\mathrm{Cok}(\Phi^{\kappa}_{R})=R^{\kappa}/R^{[\kappa]}$.\space $\mathrm{(ii)}\,R $ is $\kappa$-half-slender $\Leftrightarrow  (R^{\cdot k})^{\ast}=0$.
\quad $\mathrm{(iii)}\,R $ is half-slender iff $ (R^{\cdot \omega})^{\ast}=0$.   \\
$\mathrm{(iv)}\,R $ is  globally  half-slender iff $[\, (R^{\cdot k})^{\ast}=0\,]\, \forall \, \kappa \in \mathcal{L}_1'$. ($\mathcal{L}_1'$ is the class of non-$\omega$-measurable cardinals).    \end{definition}
\begin{remark} A ring  $R$ is is $\kappa$-half-slender iff   $ \forall F\in \mathrm{Hom}_R(R^{\kappa},R)$  
$[ \,  ( F(e_i)=0) \, \forall\, i \in \kappa \,]\Rightarrow (F=0)\,]$
\end{remark}

\begin{theorem} \label{half}
Let $R$ be a Noetherian  commutative ring. \end{theorem}
\begin{enumerate} [label=(\roman*)]
\item $[\, R$ is slender$\,]\, \Rightarrow [\, R$ is half-slender$\,]$. \\
$\mathbf{Proof:}$  See Corollary 1.3 of  \cite{eklof} .
\item If $R$ is an affine algebra  over an Artinian commutative ring $T$ then   $R$ has bounded type. \\ 
$\mathbf{Proof:}$  $R$ is a countably generated $T$-algebra hence $R$
has bounded type, according to  the comments    succeeding  Definition 3.1 of \cite{jensentwo}.
\item  Let \underline{$R_{}$ be of bounded type}. Then $[\,(R$  is slender$)\Leftrightarrow( R$  is half-slender$)\,]$  \\ $\mathbf{Proof:}$ See Corollary 3.3 of \cite{jensentwo}.

\item Let there be  proper ideals  $I,J$ of $R$  such that 
$\mathrm{(a)}\,[\,I^k+J^k=R \,]\, \forall k \in \mathbb{N}$ \space
$\mathrm{(b)}\,\,[\,  \mathbf{0}=\overset{\infty}{\underset{n=1}\cap}I^n=\overset{\infty}{\underset{n=1}\cap}J^n\,]$, \\ 
$\mathrm{(c)}\,[\,$ $I,J$ are principal ideals and $R$ is a domain$\,]$.
Then $R$ is half-slender. \\ $\mathbf{Proof:}$ See the proof of Theorem D.0.24 of page 287 of \cite{rao}, where the  authors actually prove that $R$ is half-slender and they falsely assume that this implies 
the reflexivity of $R^{\mathbb{N}},R^{\mathbb{[N]}}$ and probably the slenderness of $R$. We elaborate the proof below:  \\
Given any  proper ideal $L\neq \mathbf{0}$ of $R$ such that 
$[\, \mathbf{0}=\overset{\infty}{\underset{n=1}\cap}L^n \,]$ (1), we define the  set $B_{y, L}= \{n\in\mathbb{N} \, | \, y\in L^n \} $ and the    function 
 $\mathrm{ord}_{L}:R\rightarrow \bar{\mathbb{N}}=\mathbb{N} \,\cup \{+ \infty \} $ by the formula 
 $\mathrm{ord}_{L}(y)=\mathrm{sup}\, B_{y, L}$.  It follows from (1) that $[\, \mathrm{ord}_{L}(y)=+ \infty \Leftrightarrow y=0\,]$. It is simple to see that $[\, \mathrm{ord}_{L}(x\cdot y) \geq \mathrm{ord}_{L}(x)+ \mathrm{ord}_{L}(y)]$ and  that $[\, \mathrm{ord}_{L}(x+ y) \geq \mathrm{min} \{ \mathrm{ord}_{L}(x), \mathrm{ord}_{L}(y
 )\}\,]$. 
  We define the set $\mathcal{H}_L \subseteq R^{\mathbb{N}}$ by 
 $\{\, (x_n)_{n\in \mathbb{N}}=x\in \mathcal{H}_L \Leftrightarrow     \underset{n\rightarrow +\infty}  {\mathrm{lim}} \mathrm{ord}_{L}(x_n)=+\infty \}$.    \\
 \underline{$\mathrm{Claim\, 1_{}}$}:  $[\, \mathcal{H}_L \underset{R} \leq  R^{\mathbb{N}}\,]$ (3). \\
 Proof: Let $(a_n)_{n\in \mathbb{N}}=a\in \mathcal{H}_L $, $(b_n)_{n\in \mathbb{N}}=b\in\mathcal{H}_L  $,  
 $0 \neq r\in R$, $c=a+b$ and $d=r\cdot a$. Then \\
  $\underset{n\rightarrow +\infty}  {\mathrm{lim}} \mathrm{ord}_{L}(c_n)=+\infty $  $=\underset{n\rightarrow +\infty}  {\mathrm{lim}} \mathrm{ord}_{L}(a_n+b_n)=\underset{n\rightarrow +\infty}  {\mathrm{lim}} \mathrm{ord}_{L}(a_n)+\underset{n\rightarrow +\infty}  {\mathrm{lim}} \mathrm{ord}_{L}(b_n)=+\infty+\infty=+\infty$. Hence
  $[\,c=a+b\in  \mathcal{H}_L\,]$. Also  $\underset{n\rightarrow +\infty}  {\mathrm{lim}} \mathrm{ord}_{L}(d_n)=  $
  $\underset{n\rightarrow +\infty}  {\mathrm{lim}} \mathrm{ord}_{L}(r\cdot a_n)\geq  $  
  $\underset{n\rightarrow +\infty}  {\mathrm{lim}} \mathrm{ord}_{L}(a_n)=+\infty\Rightarrow$
   $\underset{n\rightarrow +\infty}  {\mathrm{lim}} \mathrm{ord}_{L}(d_n)=+\infty\Rightarrow$ 
   $[\,d=r\cdot a\in  \mathcal{H}_L\,]$. It follows from the above that
    $[\, \mathcal{H}_L \underset{R} \leq  R^{\mathbb{N}}\,]$.            \\
 \underline{$\mathrm{Claim\, 2_{}}$}:  $[\,R^{[N]}+L \cdot \mathcal{H}_L =\mathcal{H}_L \,] $ (4), when   $L$ is principal ideal and $R$ is a domain.  \\
Proof: $L$ is principal hence $\exists\, t \in R :[\,L=R\cdot t\, ]$. It follows from  \emph{ Krull's Intersection Theorem} (See Corollary 2.5.6 of \cite{balcerzyk}) that  $ \mathbf{0}=\overset{\infty}{\underset{n=1}\cap}L^n$. Let $a\in R^{[N]}$. Then $\{ \, \exists \,  n_0 \in \mathbb{N} : [(\, a_n=0\, )\, \forall n \geq n_0\,]\, \}$  hence $ \underset{n\rightarrow +\infty}  {\mathrm{lim}} \mathrm{ord}_{L}(a_n)=  \underset{n\rightarrow +\infty}  {\mathrm{lim}} \mathrm{ord}_{L}(0)=$
$\underset{n\rightarrow +\infty}  {\mathrm{lim}} (+\infty) =+\infty \Rightarrow[\,  a\in  \mathcal{H}_L \,]$. Hence $R^{[N]} \subseteq  \mathcal{H}_L \Rightarrow$
$[\,R^{[N]}+L \cdot \mathcal{H}_L \underset{R} \leq \mathcal{H}_L \,] $ (6). For the inverse inclusion, let $a\in \mathcal{H}_L$. Then 
 $ \underset{n\rightarrow +\infty}  {\mathrm{lim}} \mathrm{ord}_{L}(a_n)= +\infty \Rightarrow$
  $\{ \, \exists \,  n_1 \in \mathbb{N} : [(\, \mathrm{ord}_{L}(a_n)\geq 1   \, )\, \forall n \geq n_1\,] \} $. Hence  (7) $[\, a \underset{R}| a_n\,] \, \forall n\geq n_1$. Let $b\in R^{\mathbb{N}}$
  defined by $ \{\,[\,  b_n=a_n $ if $n\leq n_1\,], [\,  b_n=0 $ if $n> n_1\,]\, \}$ and $c=a-b$. It follows from (7) that  $ \{\,[\,  c_n=0 $ if $n\leq n_1\,], [\,  t \underset{R}|c_n=a_n $ if $n> n_1\,]\, \}$. Let $d=\dfrac{a-b}{t}$. Then $d\in R^{\mathbb{N}}$  and 
  $ \underset{n\rightarrow +\infty}  {\mathrm{lim}} \mathrm{ord}_{L}(d_n)= \underset{n\rightarrow +\infty}  {\mathrm{lim}} \mathrm{ord}_{L}(\dfrac{a_n}{t}\,)= \underset{n\rightarrow +\infty}  {\mathrm{lim}} [\, \mathrm{ord}_{L}(a_n)-1\,]=+\infty -1=+\infty \Rightarrow$
  $d=\dfrac{a-b}{t}\in \mathcal{H}_L$
  $\Rightarrow a=b+t\cdot d\in R^{[\mathbb{N}]}+L \cdot \mathcal{H}_L\Rightarrow$
  $[\,R^{[\mathbb{N}]}+L \cdot \mathcal{H}_L \supseteq \mathcal{H}_L \,]$ (9). Finally it follows from  $\{(6),(9) \}$ that $[\,R^{[\mathbb{N}]}+L \cdot \mathcal{H}_L =\mathcal{H}_L \,] $. \\                          
 \underline{$\mathrm{Claim\, 3_{}}$}: $[\, R^{\mathbb{N}}= \mathcal{H}_I
  + \mathcal{H}_J\,]$ (12) \\  Proof: Let  $(x_n)_{n\in \mathbb{N}}=x\in R^{\mathbb{N}}\overset{\mathrm{(a)}}\Rightarrow $
  $[\,( \exists\, y_n\in I^n )( \exists\, z_n\in J^n ):(x_n=y_n+z_n)\,]\, \forall n\in \mathbb{N}$. Let $y=(y_n)_{n\in \mathbb{N}}$ and 
 $z=(z_n)_{n\in \mathbb{N}}$ then $[\,x=y+z \,]$ (13). Also $y_n\in I^n\Rightarrow \mathrm{ord}_J(y_n)\geq n\Rightarrow$    
 $\underset{n\rightarrow +\infty}  {\mathrm{{lim}}} \mathrm{ord}_{L}(y_n) \ge \underset{n\rightarrow +\infty}  {\mathrm{{lim}}} n=+\infty \Rightarrow$
 $\underset{n\rightarrow +\infty}  {\mathrm{{lim}}} \mathrm{ord}_{L}(y_n)=+\infty\Rightarrow [\, y\in \mathcal{H}_I \,] $(14). Similarly, we show that  $[\, z\in \mathcal{H}_J \,] $(15). It follows from $\{(13),(14),(15) \}$ that  $R^{\mathbb{N}}= \mathcal{H}_I
  + \mathcal{H}_J$. So (12) has been proven.

Let now   $[\, F\in \tmop{Hom}_R(R^{\mathbb{N}},R)\,]$(16) be such that 
 (17) $[\,F(e_n)=0\,]\, \forall n\in \mathbb{N}$. We will show 
 that  $F=0$. Let $[\,M=F(\mathcal{H}_I)\,] $ (18). It follows from $\{ (3),(13) \}$ that $[\,M\underset{ R} {\leq} R\,]$ (19). By application of (4) for 
 $L=I$ we receive that  $[\,R^{[N]}+I \cdot \mathcal{H}_I =\mathcal{H}_I \,] $ (20). It follows from $\{(16),(17),(18),(19),(20) \}$ that
  $[\, M \underset{R} {\leq} J\cdot M\,]$ (21)
   $ \Rightarrow M \underset{R} {\leq} J\cdot M \underset{R,(21)} {\leq} J\cdot (J\cdot M)=J^2\cdot M\Rightarrow$ 
   $[\,M \underset{R} {\leq} J^2\cdot M \,]$ (22). Similarly, using induction we  show   
 $[\,M \underset{R} {\leq} J^n \,]\, \forall n\in \mathbb{N}\Rightarrow$
 $[\, M \underset{R} {\leq}\underset{n\in \mathbb{N}} \cap( J^n)\overset{\mathrm{(b)}}=\mathbf{0}\,]\Rightarrow[\, M=\mathbf{0}\,]\Rightarrow$
 $[\,F(\mathcal{H}_I)=\mathbf{0}\,]$ (24). Similarly $[\,F(\mathcal{H}_J)=\mathbf{0}\,]$ (25). It follows from $\{(12),(24),(25)  \}$ that $F(R^{\mathbb{N}})=\mathbf{0}\Rightarrow [\, F=0\,]$. Hence $R$ is a half-slender ring. Note that according to \cite{rao}  
  the assumption $\mathrm{(c)}$ is redundant. Here it was used only in the proof of $\mathrm{Claim\,2}$. We encourage the reader to try  to  prove  $\mathrm{Claim\,2}$  without using assumption $\mathrm{(c)}$ .

\item If $R$ is  non-local domain, then there are principal proper ideals  $I,J$ of $R$  such that \\
  $\mathrm{(a)}\,[I^k+J^k=R \,]\, \forall k \in \mathbb{N}$ \space\space
$\mathrm{(b)}\,[\,  \mathbf{0}=\overset{\infty}{\underset{n=1}\cap}I^n=\overset{\infty}{\underset{n=1}\cap}J^n\,]$            \\
$\mathbf{Proof:}$ We know that $\{ R$ is locall $\Leftrightarrow$
$\forall x\in R \, [\, x\in \tmop U(R) \text{ \, or  } (1-x) \in \tmop U(R)\,] \, \}$. See 15.15 of \cite{fuller}.  \\ 
But  $R$ is non local by asumption, hence $\{ \exists\, x\in R \, [\, x\notin \tmop U(R) \text{ \, and  \, } (1-x) \notin \tmop U(R)\,]\, \}  $ (2).
Consider the  ideals $I=R\cdot x $  , $ \, J=R\cdot(1-x)$ of $R$. They are  both proper and principal satisfying the relation $[\, I+J=R \,]$ (3). For $k\in \mathbb{N}$, we consider the polynomials $f_k(t),g_k(t)\in \mathbb{Z}[t]$ of \ref{polynomials} such that $[\, t^k\cdot f_k(t)+(1-t)^{k}\cdot g_k(t)=1\,]$ (4). Then $x^k \cdot f_k(x)+(1-x)^k\cdot g_k(x)=1  $  $\Rightarrow [\, I^k+J^k=R\,]$ (6). Let $K=\overset{\infty}{\underset{n=1}\cap}I^n$, then by applying  \emph{ Krull's Intersection Theorem} (See Corollary 2.5.6 of \cite{balcerzyk}) we get that
 $\mathbf{0}=\overset{\infty}{\underset{n=1}\cap}I^n$, since $R$ is a Noetherian domain. In a similar way we get that $\mathbf{0}=\overset{\infty}{\underset{n=1}\cap}J^n$. Therefore both $\mathrm{(a),(b)}$ are satisfied.

\item   $^{\ast}\boxed{$ If $R$ is a non-local domain, then $R$ is half-slender $}$  \\ 
$\mathbf{Proof:}$ It follows from $\mathrm{iv}$ and $\mathrm{v}$. Note that it is  proven in Remark 4.6 of  \cite{bashir} that any non-semilocal domain is half-slender.
\item  If $R$ is non-local domain  of bounded type, then $R$ is slender. \\ $\mathbf{Proof:}$ Follows from $\{\mathrm{(iii), (vi)}\}$. Note that it could also follow from Proposition 3.6 of \cite{jensentwo}.

\end{enumerate}
\begin{theorem} \label{affinealgebras2} Let $R$ be an affine algebra  over an Artinian commutative ring $T$. If $R$  is a non-Artinian domain $($i.e an integral domain that is not a field $)$  then $R$  is slender.                     
\end{theorem}
\begin{proof} $[\, T$ is Artinian$\,] \Rightarrow [\, T$ is Noetherian$\,]\Rightarrow [\, R$ is Noetherian $]$(1), as an affine algebra over a Noetherian  commutative ring. 
It follows by Lemma $\ref{hilbert}\mathrm{.i} $ that $[\,T$ is 
Hilbert$\,]$ $\overset{\ref{hilbert}\mathrm{.vi} }\Longrightarrow$ $[\,R$ is Hilbert$\,]$. Hence $[\, R$ is Hilbert and  R is non-Artinian $]$   $\overset{\ref{hilbert}\mathrm{.vii}}\Rightarrow$
$[\, R$ is non-local $]$ (2). Also $[\, R$ has bounded
 type$]$ (3)
 by $\ref{half}\mathrm{.ii}$ and $[ R$ is domain $]$(4)  by assumption. So $\{(1), (2), (3),(4)  \}  \overset {\ref{half}\mathrm{.vii}}  \Longrightarrow$  [$R$  is slender].
\end{proof}
\begin{remark} Note that Theorem \ref{affinealgebras2} above, which extends $\ref{wslenderaffinedomains2}\mathrm{.ii}$, could also follow from Proposition 3.6 of \cite{jensentwo}, although not stated there . So we do not claim \ref{affinealgebras2} as unknown but we do  claim  Theorem $\ref{half}\mathrm{.vi}$ as an unknown result.
\end{remark} 
\begin{lemma} \label{polynomials}  There are polynomials   $f_k(t),g_k(t)\in \mathbb{Z}[t]$ such that $t^k\cdot f_k(t)+(1-t)^{k}\cdot g_k(t)=1$, $\forall k \in \mathbb{N}$ .
\end{lemma} 
\begin{proof}
Consider the r.h \space $\pi:R=Z[t]\rightarrow S=\dfrac{R}{R\cdot t^k}$
defined by $\pi (h(t))=h(t)+R\cdot t^k$.   \\
$(t-1)^k-t^k=(-1)^k+(\overset{k-1}{\underset{i=1}\Sigma}$
${k}\choose {i} $
  $(-1)^{k-i}\cdot t^i) +t^k-t^k=(-1)^k (1+t\cdot g(t))$, where
  $g(t)=\overset{k-1}{\underset{i=1}\Sigma}(-1)^{i}$
  ${k}\choose {i} $ 
  $t^i$ \\
  $\Rightarrow \pi((t-1)^k)=(-1)^k \cdot \pi (1+t\cdot g(t))$ (2).
  Also $(\pi(t))^k=\pi(t^k)=t^k+ R\cdot t^k=0\Rightarrow$
$[\, \pi(t)\in  \mathfrak{N}(S)\,]$ (3) $\Rightarrow$  \\
$\pi(t)\cdot \pi(g(t)) \in  \mathfrak{N}(S)$
$\overset{\ref{units}}\Rightarrow$
$1_{S} +\pi(t)\cdot \pi(g(t)) \in  \tmop{U}(S)\Rightarrow$  
 $\pi (1+t\cdot g(t))\in \tmop{U}(S)\overset{(2)}\Rightarrow$ 
 $\pi((t-1)^k)  \in \tmop{U}(S)\Rightarrow$ \\
 $\pi( (-1)^k \cdot(1-t)^k)  \in \tmop{U}(S)\Rightarrow$ 
 $\pi((1-t)^k)  \in \tmop{U}(S)\Rightarrow$
$\exists \,g_k(t)\in \mathbb{Z}[t]:[\, \pi((1-t)^k) \cdot  \pi(g_k(t) )=1_{S} \,]\Leftrightarrow $ \\
$1-(1-t)^k \cdot \,g_k(t) \in Z[t]\cdot t^k \Leftrightarrow $
$\exists \,  f_k(t) \in Z[t]:[\, t^k\cdot f_k(t)+(1-t)^{k}\cdot g_k(t)=1        \,] $. Note that the argument is valid for $k\ge 2$ and that the claim is trivial for $k\in \{0,1  \}$.

\end{proof}

\begin{lemma} \label{units} Let $S$ be any commutative ring and $x\in S$. Then $[ \,x\in \mathfrak{N}(S)\,]\Rightarrow [\, (1+x)\in \tmop{U}(S) \,]$
\end{lemma} 
\begin{proof} By assumption $ x\in \mathfrak{N}(S)\Rightarrow[\,  \exists \, n\in \mathbb{N}: x^n=0 \,] \Rightarrow$
$(1+x)\cdot \overset{n-1}{\underset{k=1}\Sigma} \,(-1)^k\cdot x^k=1-x^n=1-0=1\Rightarrow$  $
[\,(1+x)\in \tmop{U}(S)\,]$ and 
$(1+x)^{-1}=\overset{n-1}{\underset{k=1}\Sigma} \,(-1)^k\cdot x^k$. Note that this is well known and  it is an exercise in many textbooks.
\end{proof}

\begin{lemma}\label{lost} Let $R$ be a commutative ring and $\mu,\lambda$ be cardinals such that $R^{\lambda} \underset{R}\cong R^{[\mu]}$. 
Then $:$ \\
 $\mathrm{i) }$     $\mu \ge 2^{\lambda}$  \quad   $\mathrm{ii) }$ If $\mu$ is finite then $\lambda$ is finite.
\end{lemma}
\begin{proof}
 Let $\mathfrak{m}\in \tmop{max}(R)\neq \emptyset$ 
 and let $\phi=\pi^R_\mathfrak{m}:R\rightarrow \dfrac{R}{\mathfrak{m}}=K$. Then $[\, \dfrac{R}{\mathfrak{m}}=K$ is a field$\,]$ (3). Note that $\phi$ may not be f.p (it would make the proof easier) since $R$ may not be Noetherian but  $[\,\phi$  is   f.g $]$ (2)  since $ \dfrac{R}{\mathfrak{m}}$ is  a f.g $R$-module. 
 By assumption we have 
$[\,R^{[\mu]} \underset{R}\cong R^{\lambda}\,]\overset{\ref{rh2}\mathrm{.iv} } \Rightarrow$
$[\, K^{[\mu]}\underset{K} \twoheadrightarrow K^{\lambda}\,]\overset{(3)}\Rightarrow$
$[\, \mu=\tmop{dim}_K  ( K^{[\mu]})\ge \tmop{dim}_K  ( K^{\lambda})\overset{\ref{EKextension}\mathrm{.ii}} \ge 2^{\lambda}\,] \Rightarrow $
$[\,\mu \ge 2^{\lambda}\,]$ (6).
Let now $\mu$ be finite at i.e that $[\, \mu <\omega\,]$ (5). Then it follows from  $ \{(5),(6), \ref{cardinallemma}\mathrm{.v}   \}$ that  $\lambda$ is finite. Note that for  \ref{oneil} we only need 
$\ref{lost}.\mathrm{ii}$ for Noetherian commutative rings which is  simpler to prove since $\phi$ is f.p but for the proof of $\ref{dualICFlocal}\mathrm{.i}$ we need the full statement.
\end{proof}
\begin{definition} \label{semiprimary}
A ring   R is semiprimary  $\Leftrightarrow$
$[\, (R / \mathfrak{J} (R))$ is semisimple ring, $\mathfrak{J} (R) $ is nilpotent ideal $]$
\end{definition}
\begin{remark} In   \cite{koh},  a seemingly different than \ref{semiprimary}  definition of semiprimaryness is given, where the condition that the ring   $ (R/ \mathfrak{J} (R))$ is semisimple is replaced by the condition that it is Artinian. But both definitions are equivallent in view of Proposition 15.17 of \cite{fuller}. The definition of semiprimaryness that is provided above is from page 175 of \cite{fuller}, on which   O'Neil, in \cite{oneil}, refers to. Also for the definitions of the concepts of  semisimple ring, nilpotent ideal and the Jacobson radical in  non-commutative rings the reader may look at \cite{fuller}. 
\end{remark}
Note that for Definition \ref{semiprimary} above and for Lemma \ref{oneil}
below the ring $R$ does not have to be commutative.  
\begin{lemma} \label{oneil} Let $R$ be a ring, $\mathfrak{J} (R)$ be it's \emph{Jacobson radical} and  $\lambda$ be an infinite cardinal. 
\end{lemma} 
 \begin{enumerate} [label=(\roman*)]
 
\item If $R^{\lambda} $ is non-f.g and  free left $R$-module then $R$ is semiprimary .  \\ $\mathbf{Proof:}$ It is Lemma 5.1 of \cite{oneil}.

\item  R is left Artinian   $\Leftrightarrow$ 
$[\,  (R/ \mathfrak{J} (R))$ is semisimple ring, $\mathfrak{J} (R) $ is nilpotent ideal, R is left Noetherian$]$ \\
$\mathbf{Proof:}$ It is  Proposition 2.3.26 of \cite{enochs}.
\item Let R be Noetherian commutative ring. If $ R^{\lambda}   \in \tmop{Free}(R) $  then $R$  is Artinian . \\  $\mathbf{Proof:}$ By assumption $[\,R$ is Noetherian$\,]$ (1). Also by assumpion $[\, R^{\lambda}   \in \tmop{Free}(R)\,] \overset{\ref{lost}\mathrm{.ii}}\Rightarrow $
 $[\, R^{\lambda}$ is a non f.g free left $R$-module$\,]$
$\overset{\ref{oneil}\mathrm{.i}}\Rightarrow$
$[\, R$ is a semiprimary ring  $]$(2). Therefore it follows from 
$\{(1), (2), \ref{semiprimary},\ref{oneil}\mathrm{.ii}   \}$ that the ring $R$ is Artinian.
 
\end{enumerate}
We are now ready to answer  Question 1  in the case where $R$ is a Noetherian commutative ring.  
\begin{theorem}\label{relief} $^{\ast}$  
Let $R$ be a Noetherian commutative ring and  $\lambda$ is an infinite cardinal. Then \\
$$\boxed{R^{\lambda}\in \tmop{Proj}(R) \Leftrightarrow R$ is Artinian  
$}$$

 \end{theorem} 
$ \bullet$ Proof of $[\, R^{\lambda}\in \tmop{Proj}(R) \Leftarrow R$
 is Artinian $]$:  Follows from \ref{collection}$\mathrm{.ii}$,
  \ref{collection}$\mathrm{.iii}$ \\
 . \quad  $\bullet$ Proof of $[\, R^{\lambda}\in \tmop{Proj}(R) \Rightarrow R$
 is Artinian $]$: By applying     \ref{collection}$\mathrm{.vi}$ we get that 
 $R=\overset{n}{\underset{k=1}\Pi} R_k $   (3), for some $\{$[ connected Noetherian commutative rings  $R_k  ]\,\forall\, k\in \mathrm{T_n} \}$ (4). Let  $k\in \mathrm{T_n}$  and $\phi: R\rightarrow R_k$ like in \ref{strange2} then $[$ $\phi$ is f.p $]$(5). So  
 $[\,  R^{\lambda}\in \tmop{Proj}(R)\,]\underset{(5)} {\overset{\ref{remainprojective}}\Rightarrow}$
$[\, (R^{\lambda})_{\phi}\underset{R_k}\cong (R_k)^{\lambda}\in \tmop{Proj}(R_k)\,]$  $\underset{(4)} {\overset{\ref{Bass}}\Rightarrow } [\, (R_k)^{\lambda}\in \tmop{Free}(R_k)\,]$ 
$\overset{\ref{oneil}\mathrm{.iii}}\Rightarrow$
$[\,R_k$ is Artinian$\,]$. So 
$[\, R_k$ is Artinian $]\, \forall\, k\in \mathrm{T_n} \underset{\ref{collection}} {\overset{(3)} \Rightarrow}$
$\overset{n}{\underset{k=1}\Pi} R_k=$ 
\underline{$R_{}$ is Artinian},
as a product of Artinian rings.

\begin{remark}
Note that $\ref{oneil}\mathrm{.iii}$ is valid even if $R$ is non-commutative but  the argument used in the proof \ref{relief} does not apply for non commutative rings.  Also note that in page 287 of   \cite{rao}, it  is  claimed  that  $ \bigstar$[If $R$ is an uncountable-principal-complete-$\text{d.v.r}$ then $R^{\mathbb{N}}\in \mathrm{Free}(R)]$, which is not compatible to \ref{relief} since   $\text{d.v.r}$'s are not Artinian and free modules are always projective. Hence $\bigstar$ should be invalid.
\end{remark}

\section{Lazarus and Goldie dimensions} \label{lazarus-section}

\begin{lemma}\label{shortEKrings}
Let $R$ be a commutative ring and $M\in \mathrm{Free}(R)$. If
$[\, \dim_R(M)\geq \omega$ or $| R| \geq \omega\,]$ then 
$\boxed{\big| M \big|= \max \{ \big| R\big|, \dim _R (M) \}}$, i.e
$   | R^{[\lambda]}|= \max \{ | R |, \,\lambda \} $, when $[\, \lambda\geq \omega
$ or $|R|\geq \omega\,]$.
\end{lemma} 

\begin{proof} By assumption   $M\underset{R}\cong R^{[\lambda]}$ for  $\lambda=\dim_R(M)$. So the claim follows from  \ref{cardinallemma1}.
\end{proof} 
Let $R$ be a commutative ring and $S\subseteq R$ be multiplicative (i.e $1\in S$, $0\notin S$ and S is closed under multiplication). The r.h  $\psi=\mathrm{i}_R^{S} :R\rightarrow S^{-1}(R)=R_{\centerdot S}$ is defined by $\psi(r)=[\,\frac{r}{s}\,]$. If $S=R \smallsetminus \mathrm{Z}(R)$ then $S^{-1}(R)=\mathrm{Q}(R)$ is the total ring of fractions of $R$ and $\mathrm{i}_{R}=\mathrm{i}_R^{S}: R \rightarrow \mathrm{Q}(R)$.
Let  also $M\in R$-$\mathrm{Mod}$ then $ S^{-1}(M)\overset{\text{def}}=$ $\{\, [\, \frac{m}{s}\,]\, \big| \, m\in M,\,s\in S    \}$
\begin{lemma} \label{localization}
Let $R$ be a commutative ring, $\kappa \in \mathrm{Card}$ and $M\in R$-$\mathrm{Mod}$. Let also $S\subseteq R$ be multiplicative,  $\psi=\mathrm{i}_R^{S} :R\rightarrow R_{\centerdot S}$ and $\phi=\mathrm{i}_{R}: R \rightarrow T=Q(R)$. Assume that $[\, S\subseteq R\smallsetminus \mathrm{Z}(R)\,]\,(1)$. Then :
\end{lemma}
\begin{enumerate} [label=(\roman*)] 
\item It is well known that : $\mathrm{(a)}\,\, \psi $ is a flat r.h \quad $\mathrm{(b)}\,\,  M_{\psi}\underset {R_{\psi}}\cong S^{-1}(M)$
\item
 $ [\, M$ is $R$-torsion$\,] \Rightarrow$ 
$[\, M_{\psi}$ is $R_{\psi}$-torsion$\,]$ \\               
$\mathbf{Proof:}$  Let $\xi=[\,\frac{m}{s} \,]\in $
$  S^{-1}(M) \underset {R_{\psi}}\cong M_{\psi}$ such that $m\in M,\, s \in S$. Then  $\exists \, r\in (R\smallsetminus \mathbf{0}):r\cdot m=0$, since $M$ is $R$-torsion. Note that $r\cdot s\neq 0$, since
$s\in S\subseteq  R\smallsetminus \mathrm{Z}(R)$ and $(r\cdot s) \cdot \xi= r\cdot s \cdot [\,\dfrac{m}{s} \,]=r\cdot m=0$.
\item 
   $\mathrm{(a)}$ $\psi$ is 1-1  \quad  $\mathrm{(b)}$ $\big|R \big|\leq \big| R_{\psi} \big|$ \quad\quad
   $\mathrm{(c)}$ 
 $ \big| R \big| \leq  \big| R_{\psi} \big|\leq  \big| R \big|^{\, 2}$  \\               $\mathbf{Proof:}$ Note that the r.h 
 $\psi :R\rightarrow R_{\centerdot S}\cong R_{\psi}$ is  defined by $g(r)=[\,\frac{r}{1}\,]$. The assumption (1) implies that $[\,$g is 1-1$\,] \,\mathrm{(a)}\Rightarrow [\,  | R \big|\leq \big| R_{\psi} |\,]\, \mathrm{(b)} $. Consider the function $h:R \times S\rightarrow  R_{\centerdot S}\cong R_{\psi}$, defined by $h((r,s))=[\, \frac{r}{s}\,]$. Trivially 
$[\, h$ is onto$\,]\Rightarrow \big| R_{\psi} \big|\leq \big| R\times S \big|=\big| R \big| \cdot \big|  S \big|\leq$
$\big| R \big| \cdot \big|  R \big|= \big|  R \big|^2\Rightarrow$
$[\,\, \big| R_{\psi} \big|\leq  \big| R \big|^2\,\,]$ (4).
It follows from $\{ \mathrm{(b)},\, (4) \} $ that 
 $ \big| R \big| \leq  \big| R_{\psi} \big|\leq  \big| R \big|^{\, 2}\, \mathrm{(c)}$ 
\item Assume that $ M$ is a torsionless $R$-module. Then  
$ \big|  M \big| \leq \big|  M_{\psi} \big|$.\\ $\mathbf{Proof:}$
The function $g:M\rightarrow  S^{-1}(M)$ defined by $g(m)=[\,\frac{m}{1} \,]$ is trivially $R$-linear. There is $\lambda\in \mathrm{Card}$ and $h\in \mathrm{Hom}_R(M,R^{\lambda})$ such that $[\, h$ is 1-1$\,]$ (7), since $M$ is torsionless. Let $m\in \mathrm{Ker}(g)$. Then $g(m)=[\,\frac{m}{1} \,]=0\Rightarrow[ \,  \exists\,  s\in S: s\cdot m=0\,]\Rightarrow$
 $[\, s\cdot h(m)=0\,]$ (8). So $[\, s\in  R\smallsetminus \mathrm{Z}(R),\, h(m)\in R^{\lambda},\,(8)\,]\Rightarrow [\, h(m)=0\,]\overset{(7)}\Rightarrow [\, m=0\,]$. Hence $[\, \mathrm{Ker}(g)=\mathbf{0}\,]\Rightarrow[\, g$ is 1-1$\,]\Rightarrow$
 $[\,  \big|  M \big| \leq \big| S^{-1}(M) \big|\overset{\mathrm{(i.b)}}=
 \big| M_{\psi} \big|\,]\Rightarrow$
  $[\,  \big|  M \big| \leq
 \big| M_{\psi} \big|\,]$

\item Assume that $[\,|R|\leq \omega\, $  or  $\, \kappa \geq \omega\,]$ (11). Then  \quad   
  $ \big|R^{[\kappa]} \big|=\big|(R_{\psi})^{[\kappa]} \big|$ \\
 $\mathbf{Proof:}$  $\centerdot$ Let
 $[\,|R| <\omega\,  , \, \kappa \geq \omega\,]$. It follows from 
 $\mathrm{(iii.c)}$ that $[\,|R| <\omega\, , \,|R_{\psi}| <\omega\,]
 $ (12) and $[\, \kappa \geq \omega\,]$ (13). So 
 $[\,  | R^{[\kappa]}|\overset{\ref{shortEKrings}}= \max \{ | R |, \,\kappa \}= \kappa \overset{\ref{shortEKrings}}=$
 $ \max \{ | R_{\psi} |, \,\kappa \}= | (R_{\psi})^{[\kappa]}| \,]\Rightarrow$
 $ [\, |R^{[\kappa]} |=|(R_{\psi})^{[\kappa]} |\,]$.\\
  $\centerdot$ Let$[\,|R| \geq\omega\,  , \, \kappa < \omega\,]$.
   It follows from  $\mathrm{(iii.c)}$ that $[\, | R|=| R_{\psi}|\,]$
   (15). So \\
 $[\,  | R^{[\kappa]}|\overset{\ref{shortEKrings}}= \max \{ | R |, \,\kappa \} \overset{(15)}=$
 $ \max \{ | R_{\psi} |, \,\kappa \}= | (R_{\psi})^{[\kappa]}| \,]\Rightarrow$
 $ [\, |R^{[\kappa]} |=|(R_{\psi})^{[\kappa]} |\,]$.\\
 It follows from the assumption (11) that we have covered all cases.  
\item Assume that $[\,|R|\leq \omega\, $  or  $\, \kappa \geq \omega\,]$ (11). Then    $ \big|R^{\kappa} \big|=\big|(R_{\psi})^{\kappa} \big|$ \\ $\mathbf{Proof:}$ It follows from   $\mathrm{(iii.c)}$ that
$ [\, | R |^{\, \kappa} \leq  | R_{\psi} \big|^{\,\kappa}\leq  (| R |^2)^{ \kappa}= | R |^{\, 2 \kappa}\,]\,(17) $ \\ 
 $\centerdot$ Let
 $[\,|R| <\omega\,  , \, \kappa \geq \omega\,]$. Then $[\, \kappa \geq \omega\,]\Rightarrow [\, 2 \kappa=\kappa\,]\, (18)$. Hence the claim follows from $\{(17),\, (18)  \}$\\
$\centerdot$ Let$[\,|R| \geq\omega\,  , \, \kappa < \omega\,]$.
Then $[\,|R| \geq \omega\, ]\Rightarrow [\,|R|^{\,2}= |R|\,  ]$ (19).
   The claim  follows from $\{(17),\, (19)  \}$.

\item   Assume that $[\,|R|\leq \omega\, $  or  $\, \kappa \geq \omega\,]$. Then : \quad   $\mathrm{(a)}$ 
  $ \big|R^{[\kappa]} \big|=(\mathrm{Q}(R))^{[\kappa]} \big|$
        \quad  $\mathrm{(b)}$ 
   $ \big|R^{\kappa} \big|=\big|(\mathrm{Q}(R))^{\kappa} \big|$\\
$\mathbf{Proof:}$  It follows from the above $\{ \mathrm{(v)},\, \mathrm{(vi)}  \}$  applied for $S= R\smallsetminus \mathrm{Z}(R)$.
  \
\end{enumerate}

So far, we have defined $\dim _R (M)$, only when 
$M\in \mathrm{Free}(R)$. Next we define the \emph{Lazarus dimension} 
 $ \mathrm{L.dim}_R(M)$, in a way that will extend  both $\ref{shortEK},\, \ref{EKextension}\mathrm{.i}$ to larger classes of   commutative rings and  modules.
 
\begin{definition} \label{defldim}
Let $R$ be a commutative ring,  
$M \in R$-$\mathrm{Mod}$, $\kappa \in \mathrm{Card}$ and $\phi=\mathrm{i}_{R}: R \rightarrow T=Q(R)$ \\
$\mathrm{(i)\,\,}\mathrm{L}_R(M)=  \{ \kappa \in \mathrm{Card} \,|\, \exists \, S \in \tmop{Max}(\tmop{Ind}_R(M)) :[\, \kappa=\mathrm{Card}(S)\,] \}$. Note that it follows from \ref{torsion}
that $\{ \kappa\in \mathrm{L}_R(M)\Leftrightarrow \, \exists f \in \mathrm{Hom}_R(R^{[\kappa]},M) :[\, f \text{\, is \,}1-1, \, \mathrm{Cok}(f)\text{\, is \,} R$-torsion$] \, \}$.\\
 $\mathrm{(ii)\,\,} \mathrm{L.dim}_R(M) =\kappa$  $   \Leftrightarrow  \mathrm{L}_R(M)=\{ \kappa \} $.
 Then we say   that $M$ has $\mathrm{L.dimension}$ or that $\mathrm{L.dim}_R(M)$ exists. \\ 
 $\mathrm{(iii)\,\,}\mathrm{dim}_R(M)=\kappa$ $ \Leftrightarrow$
$ M \underset{R}\cong R^{\, [\kappa]}$. This is repetition of 
\ref{defdimension}. In this case we say   that $M$ has 
$\mathrm{dimension}$.\\ 
 $\mathrm{(iv)}\,\,\mathrm{rank}_R(M)=\kappa$  $\Leftrightarrow$
$ M_{\phi} \underset{T}\cong T^{\, [\kappa]}$. In this case we say that $M$ has $\mathrm{rank}$ (Well established terminilogy)\\
$\mathrm{(v)}\, \, R$ is an $\mathrm{L}$-ring $ \Leftrightarrow$
 $[\, \mathrm{L.dim}_R(N)$ exists$\,]\, \forall N\in R$-$\mathrm{Mod}$.

 \end{definition}
 Let  $S\subseteq R$ be multiplicative,  $\psi=\mathrm{i}_R^{S} :R\rightarrow S^{-1}(R)$ and $\phi=\mathrm{i}_{R}: R \rightarrow T=Q(R)$. 
\begin{lemma} \label{lazarus2}
Let $R$ be a commutative ring, $\kappa \in \mathrm{Card}$ and $M\in R$-$\mathrm{Mod}$. Let also $S\subseteq R$ be multiplicative,  $\psi=\mathrm{i}_R^{S} :R\rightarrow R_{\centerdot S}$ and $\phi=\mathrm{i}_{R}: R \rightarrow T=Q(R)$. Assume that $[\, S\subseteq R\smallsetminus \mathrm{Z}(R)\,]\,(1)$. Then :
\end{lemma} 
 
  \begin{enumerate} [label=(\roman*)] 
   \item $\{\, [\,\mathrm{dim}_R(M)=\kappa\,\,]$, $[\,R$ is Noetherian or $\kappa \geq \omega\,]\, \}$
 $\Rightarrow [\, \mathrm{L.dim}_R(M)=\kappa\,]$ \\
 $\mathbf{Proof:}$ Follows from $\{ \ref{Lazarus},\, \ref{defldim}\mathrm{.ii},\, \ref{defldim}\mathrm{.iii}\}$  
  \item $[\,\mathrm{dim}_R(M)=\kappa\, \,]$ $\Rightarrow [\, \mathrm{rank}_R(M)=\kappa\,]$\\
  $\mathbf{Proof:}$ Follows from $\{ \ref{rh}\mathrm{.i},\, \ref{defldim}\mathrm{.iii},\,\ref{defldim}\mathrm{.iv} \}$

 \item  $ [\, \dim_{R_{\psi}}(M_{\psi})=\kappa$,
 $R_{\psi}$ is Noetherian $]\Rightarrow$   
 $[\, \mathrm{L.dim}_R(M)=\kappa \,]$. \\
$\mathbf{Proof:}$ 
 Let $\lambda \in \mathrm{L}_R(M)$. Then 
$ \exists f \in \mathrm{Hom}_R(R^{[\lambda]},M) :[\, f \text{\, is \,}1-1, \,]\,(2), [\, \mathrm{Cok}(f)\text{\, is \,} R$-torsion$]\,(3) $. It follows from $\{\ref{localization}\mathrm{.i.a} ,(2)\}$ that $[\, f_{\psi}$ is 1-1$\,]$ (4).
Also $\mathrm{Cok}(f_{\psi})\underset{R_{\psi}}\cong (\mathrm{Cok}(f))_{\psi}\overset{(3),\, \ref{localization}\mathrm{.ii} } \Longrightarrow$  $[\, \mathrm{Cok}(f_{\psi})$ is $R_{\psi}$-torsion$\,]$ (5). Note that $f_{\psi}:(R^{[\lambda]})_{\psi}\rightarrow M_{\psi}$,  
$(R^{[\lambda]})_{\psi}\underset{R_{\psi}}\cong (R_{\psi})^{[\lambda]}$ and 
$M_{\psi}\underset{R_{\psi}}\cong (R_{\psi})^{[\kappa]}$.
Therefore it follows from  $\{(4),(5),\,\ref{torsion},\ \ref{Lazarus}\mathrm{.iv}  \}$ that $[\, \lambda=\kappa \,]\Rightarrow $
$[\, \mathrm{L}_R(M)=\{ \kappa \}\,]\overset{\ref{defldim}\mathrm{.ii}}\Rightarrow$ 
 $[\, \mathrm{L.dim}_R(M)=\kappa \,]$      
\item If $ R_{\psi}$ is a field, then:  $\mathrm{(a)\,}\, \mathrm{L.dim}_R(M)=\mathrm{rank}_{R}(M)$
$=\dim_{R_{\psi}} (M_{\psi})\,]$ \space  $\mathrm{(b)}\,R$ is $\mathrm{L}$-ring. \\
$\mathbf{Proof:}$ Follows from $\mathrm{(iii)}$ and it's proof, since all $R_{\phi}$-modules are free  and fields are Noetherian.
\item If $R$ is  a domain then   
  $\mathrm{(a)\,}\, \mathrm{L.dim}_R(M)=\mathrm{rank}_R(M)$
$=\dim_T (T\underset{R} \otimes M)\,]$ \quad   $\mathrm{(b)}\,R$ is an $\mathrm{L}$-ring.\\
 $\mathbf{Proof:}$ Follows from $\mathrm{(iv)}$ applied for $S=R\smallsetminus \mathbf{0}$, since $T=\mathrm{Q}(R)$ is a field.
\item  If $[\, \kappa \geq \omega,\, R$ is a domain$\,]$ then $\mathrm{L.dim}_R(R^{\kappa}) \geq |R |$  \\
$\mathbf{Proof:}$ Follows from $\{ \ref{vandermonde2},\, \ref{lazarus2}\mathrm{.iv.b} \}$

\end{enumerate} 
 
\begin{theorem} \label{shortEKlazarus}
Let $M$ be a tonsionless module over a domain $R$. If  $[\, \big| R \big| \geq \omega$  or  $ \mathrm{L.dim} _R (M)\geq \omega          \,]\,(1)$ then  \space
$\big| M \big|= \max \{ \big| R\big|, \mathrm{L.dim} _R (M) \}$
\end{theorem}
\begin{proof}
Let $F=\mathrm{Q}(R)$ and $\phi=\mathrm{i}_{R}: R \rightarrow F=Q(R)$ and  $[\, \kappa=\mathrm{L.dim}_R(M)\,]$ (2). Then $[\, R^{[\kappa]} \underset{R} \rightarrow M\,]$ (3), $[\, \kappa=\dim_F M_{\phi}\,]$ (4) and 
$[\,  M_{\phi} \underset{F}\cong F^{[\kappa]}\,]$ (5). It follows from $\ref{localization}\mathrm{.iv}$ that $\big| M \big| \leq \big| M_{\phi} \big|$ (6) and from $\{\ref{localization}\mathrm{.iv},\, (1)\}$ that 
$[\, |R^{[\kappa]}|=|F^{[\kappa]}|\,]$ (7). Therefore 
$\big| R^{[\kappa]} \big|  \underset{(3)}\leq  \big| M \big| \underset{(6)} \leq \big| M_{\phi} \big|\underset{(5)}= \big| F^{[\kappa]} \big| \underset{(7)}\Rightarrow$
$\big| M \big|= \big| R^{[\kappa]} \big| \underset{\ref{shortEKrings}}=$
$\max\{\big| R \big| ,\kappa \} \underset{(2)}=\max \{ \big| R\big|, \mathrm{L.dim} _R (M) \}\Rightarrow$
$\big| M \big|= \max \{ \big| R\big|, \mathrm{L.dim} _R (M) \}$
\end{proof}

\begin{remark}
Theorem \ref{shortEKlazarus} above  generalizes Theorem \ref{shortEK} and Theorem \ref{EKlazarus} below generalizes Theorem $\ref{EKextension}\mathrm{.i}$ to domains replacing $\mathrm{dim}_R(M)$ by $\mathrm{L.dim}_R(M)$. 
\end{remark}

\begin{theorem} \label{EKlazarus}
Let $R$ be a domain and $\lambda$ be a cardinal. If $\lambda\geq \omega$ then  
$\boxed{ \mathrm{L.dim} _R (R^{\lambda})= \big| R^{\lambda}  \big| }$
\end{theorem}
\begin{proof}
By assumption $[\, \lambda \geq \omega\,]\Rightarrow$
 $[\, R^{[\omega]} \underset{R} \rightarrowtail R^{[\lambda]}  \underset{R} \rightarrowtail R^{\lambda}\,]\Rightarrow$
  $[\, R^{[\omega]}   \underset{R} \rightarrowtail  R^{\lambda}\,]$ (1).
 It follows from  $\{(1),\, \ref{lazarus2}\mathrm{.v}  \}$ that $[\, \mathrm{L.dim}_R(R^{\lambda})\geq \omega\,]$ (2).
  The module $M=R^{\lambda}$ is $R$-torsionless. Hence it follows from $\ref{shortEKlazarus},\, (2)$ that $ \big| R^{\lambda} \big| =\max \{ \big| R\big|, \mathrm{L.dim} _R (R^{\lambda}) \}$ (3). Also it follows from \ref{vandermonde2} that
 $[\,  \mathrm{L.dim} _R (R^{\lambda})\geq \big| R \big|\,] $ (4).
 Hence it follows from $\{ (3),(4) \}$ that $ \mathrm{L.dim} _R (R^{\lambda})= \big| R^{\lambda}  \big| $ \end{proof}

\begin{remark} \label{lazarus-remark}  Let $R$ be a Noetherian commutative ring, $M\in R$-$\mathrm{Mod}$ and $\lambda \in \mathrm{Card}_{\geq \omega}$.  \\
$\mathrm{(i)}$ In \cite{lazarus}, it is proven that  $[\,M 
\underset{R} \rightarrowtail F,\, F\in \mathrm{Flat}(R)\,]\Rightarrow[\,  \mathrm{L.dim}_R(M)$ exists$\,]$ \\
$\mathrm{(ii)}$ Assume that   $M$  is a torsionless $R$-module. Then 
$\exists \,  \lambda \in \mathrm{Card}:[\, M\underset{R} \rightarrowtail R^{\lambda}\,]$
$ \overset{\ref{collection} \mathrm{.iii}} { \underset{\mathrm{(i)}} \Longrightarrow }$                                  
$[\,  \mathrm{L.dim}_R(M)$ exists$\,]$ \\ 
\underline{$\mathrm{Question\,\, B}$} : Is Theorem \ref{shortEKlazarus} valid under weaker assumption than 
$R$  being is a domain ?             \\                    
$\mathrm{(iii)}$ It follows from $\mathrm{(ii)}$ that 
$\mathrm{L.dim}_R(R^{\lambda})$ exists. \underline{$\mathrm{Question\,\, C}$}:  $\mathrm{L.dim} _R (R^{\lambda})=\, ? $ \\
$\mathrm{(iv)}$ In \cite{lazarus}, it is proven that if $R$ is reduced then $R$ is an $\mathrm{L}$-ring.  \\
$\mathrm{(v)}$ In \cite{lazarus}, it is  proven that if
$[\, R^{[\omega]} \underset{R} \rightarrowtail M\,]$ then $[\,  \mathrm{L.dim}_R(M)$ exists$\,]$. 
\end{remark}

\begin{definition}\label{defgoldie}
Let $R$ be a ring (with identity) and $M\in R$-$\mathrm{Mod}$.
 Then : \\
$\mathrm{(i)}\,\mathrm{G}_R(M)=\{\kappa \in \mathrm{Card}\, \big|\,    [\,\forall\, i \in \kappa\,][\, \exists\,  M_i \in R $-$\mathrm{Mod}\,]:[\, \underset{i\in \kappa}  \oplus M_i \underset{R}  {\overset{_l}\rightarrowtail}  M \,]\,\}$ \quad
$\mathrm{(ii)}\,\mathrm{G.dim}_R(M)=\mathrm{sup}\, \mathrm{G}_R(M)$
\end{definition}

\begin{theorem} \label{dimensions1}Let $R$ be a ring, $\{M,N\}\subseteq R$-$\mathrm{Mod}$, $\{\kappa,\lambda\}\subseteq \mathrm{Card}$ and $\mu=\mathrm{G.dim}_R(M)$. Then
\end{theorem}
\begin{enumerate} [label=(\roman*)] 
\item Trivially  $[\,  \kappa \in \mathrm{G.dim}_R(M)\,]\Leftrightarrow[\,  M$ contains an internal direct sum of $\kappa$ submodules$\,]$
\item Trivially $\mathrm{G}_R(M)$ is downward closed.
\item If $\mu=\omega$ then $ \mu \in \mathrm{G}_R(M)$. \\
$\mathbf{Proof}$ : It is Lemma 5 of \cite{fuchs}.  Note that, in Theorem 6 of  \cite{fuchs}, it is  even proven that $ \mu \in \mathrm{G}_R(M)$ when $\mu$ is any  strongly inaccessible cardinal. 
\item If $[\, \kappa \geq \omega\, ,\, \kappa \in \mathrm{G.dim}_R(M)  \,]$  then $M$ is not a  left Noetherian  $R$-module.  \\
$\mathbf{Proof}$ : It follows from the assumption that $[\, \omega \in \mathrm{G.dim}_R(M)\,]$, since  $\mathrm{G.dim}_R(M)$ is downward closed. Hence $[\, \underset{i\in \omega}  \oplus M_i  \subseteq  M \,]\,(1)$, for some  $[\, 0 \neq M_i \in R $-$\mathrm{Mod}\,]\,(2)$. We define $[\, L_n= \overset{n}{\underset{i=0}\Sigma} M_i\,]\, (3)$,
for $n\in \omega$. It follows from $\{(1),(2),(3)\}$ that the sequence $(L_n)_{n\in \omega}$ is strictly increasing sequence of left $R$-submodules of $M$. Therefore $M$ is not a left Noetherian $R$-module. 
\item If $M$ is a  Noetherian left $R$-module then $ \mu= \mathrm{G.dim}_R(M) < \omega$. \\
$\mathbf{Proof}$: It follows from $\{ \mathrm{(iii)}, \mathrm{(iv)}  \}$.
\item If $M$ is a left Noetherian $R$-module. Then
$[\, M\oplus N \underset{R}{ \rightarrowtail} M\,]\Rightarrow[\,  N=\mathbf{0}\,]$   \\
$\mathbf{Proof}$: Assume that $[\, N \neq \mathbf{0}\,]\,(4)$. We assume
$[\, M\oplus N \underset{R}{\rightarrowtail} M\,]\,(5)\Rightarrow[\, M\oplus N \oplus N \underset{R}{\rightarrowtail}   M\oplus N \,]\Rightarrow$
$[\, M\oplus N \oplus N \underset{R}{\rightarrowtail}   M\oplus N \,]\overset{(5)}\Rightarrow$
$[\, M\oplus N^2 \underset{R}{\rightarrowtail}   M \,]$. In a similar way, by induction, we can show that $[\, M\oplus N^n \underset{R}{\rightarrowtail}   M \,]\, \forall\, n\in \omega$. Hence $[\, \overset{n}{\underset{i=1}\oplus} N \underset{R}{\rightarrowtail}   M \,]\, \forall\, n\in \omega$  $\overset{(4)}\Rightarrow [\, \omega \subseteq   \mathrm{G}_R(M)\,]\Rightarrow$    $ [\, \omega \leq   \mathrm{G.dim}_R(M)\,]\overset{\mathrm{(v)}}\Rightarrow$
$M$ is a not a left   Noetherian $R$-module, which contradicts the assumpion. Therefore
$N=\mathbf{0}$.  
\item If $R$ is left Noetherian and $[\ \kappa<\omega,\, \lambda< \omega \,]$  then  $[\, R^{[\kappa]} \underset{R}{\rightarrowtail} R^{[\lambda]} \,]\Rightarrow [\,\kappa \leq \lambda \,]$ \\
$\mathbf{Proof}$ : Assume that $ \kappa > \lambda\ $. Then $[\, \lambda=\kappa+ \mu\,]\,(6)$  for some cardinal $\mu $ such that $ 1 \leq \mu < \omega $. So, $[\, R^{[\kappa]}$ is  left Noetherian$\,]\,(7)$  since $[\,R$ is left Noetherian,  $\kappa< \omega  \,]$ and $[\,R^{[\mu]}\neq \mathbf{0}\,]\,(8)$ since $\mu \geq 1$. So
 $[\, R^{[\kappa]} \underset{R}{\rightarrowtail} R^{[\lambda]} \,]\overset{(6)}\Rightarrow [\, R^{[\kappa]} \underset{R}{\rightarrowtail} (R^{[\kappa]}\oplus R^{[\mu]}) \,]\overset{\mathrm{(iv)},(7)}\Longrightarrow$
 $[\, R^{[\mu]}= \mathbf{0}\,]$, which contradicts (8). Hence $[\, \kappa \leq \lambda\,]$.
\item If $R$ is left Noetherian  then  $[\, R^{[\kappa]} \underset{R}{\cong} R^{[\lambda]} , \,  \kappa < \omega, \, \lambda < \omega\,]\Rightarrow [\,\kappa = \lambda  \,]$ \\
$\mathbf{Proof}$ :  $[\, R^{[\kappa]} \underset{R}{\cong} R^{[\lambda]} , \,  \kappa < \omega, \, \lambda < \omega\,]\Rightarrow [\, R^{[\kappa]} \underset{R}{ \rightarrowtail} R^{[\lambda]}     ,\,  R^{\kappa} \underset{R}{ \rightarrowtail} R^{[\lambda]},\,   \kappa < \omega, \, \lambda < \omega  \,]\overset{\ref{dimensions1}\mathrm{.v}} \Longrightarrow$
$[\, \kappa \leq \lambda\, ,\,  \lambda \leq \kappa    \,]\Rightarrow [\, \kappa=\lambda\,]$. 
\item Let $R$ be left Noetherian and $a,b$ be elements of $R$. The  equation  $x\cdot a= y \cdot b$  has non-zero solutions.  \\
$\mathbf{Proof}$ : $\centerdot $ Let $(a,b)=(0,0)$. Then the claim is trivial \\$\centerdot $ Let $(a,b)\neq(0,0)$. Consider the function  $h:R^2 \rightarrow R$ be defined by $h(x,y)=x\cdot a- y \cdot b$. Trivially $h$ is left $R$-linear. If the equation $x\cdot a= y \cdot b$   had no nontrivial solutions then $h$ would be 1-1. Therefore 
 $ R^{[2]} \underset{R}{\overset{_l} \rightarrowtail} R^{[1]}$, which contradicts  $\ref{dimensions1}\mathrm{.viii}$. Therefore the  equation has non-trivial solutions.
\end{enumerate}

\begin{theorem} \label{dimensions2}
$^{\ast}$ Let $R$ be a commutative  domain, $M$ be a torsionless $R$-module and $\lambda \in \mathrm{Card}$.
Then $:$ \\ $\mathrm{\,(i)\,G.dim}_R(M)= \mathrm{L.dim}_R(M)$. \quad\quad \quad\quad \quad\quad \quad\quad
$\mathrm{(ii)} $ If  $\lambda \geq \omega$ then   
$ \mathrm{G.dim}_R(R^{\lambda})=|R^{\lambda}|$.
\end{theorem}
\begin{proof} \underline{$\mathrm{Proof \, of \, (i)}$} :
Let $\lambda=\mathrm{L.dim}_R(M)$. Then $[\, \underset{i\in \lambda}\oplus R \underset{R} \cong R^{[\lambda]}\underset{R}\rightarrowtail M\,]\Rightarrow [\, \lambda \in \mathrm{G}_R(M)\,]\Rightarrow$
$[\, \lambda \leq \mathrm{sup}\,\mathrm{G}_R(M)= \mathrm{G.dim}_R(M)\,] \Rightarrow$
$[\, \mathrm{L.dim}_R(M)\leq \mathrm{G.dim}_R(M)  \,]\,(1)$. Note that (1) holds in general. Let $[\, \kappa\in \mathrm{G}_R(M)\,]\,(2)$. Then 
$[\,  \underset{i\in \kappa} \oplus M_i \subseteq M\,]\,(3)$, for some submodules $[\, \mathbf{0} \neq M_i\,]\,(4)$ of $M$. 
Pick $(5)\,[\, 0\neq m_i \in M_i\,]\, \forall \, i \in \kappa$ and let 
$S=\{m_i\,|\, i\in \kappa  \}$. It follows from $\{(3),(5)   \}$  that 
$[\,|S|=\kappa$, $S\in \mathrm{Ind}_R(M)\,]\,(6)$ and from $\{ \ref{zornapplication1},\, \ref{lazarus-remark}\mathrm{.ii}, \, (6) \}$
that $[\, \kappa\leq  \mathrm{L.dim}_R(M)\,]\overset{(2)}\Rightarrow $
$[\, \mathrm{sup}\,\mathrm{G}_R(M)\,]\leq \mathrm{L.dim}_R(M)\,]\Rightarrow[\,\mathrm{G.dim}_R(M)  \leq   \mathrm{L.dim}_R(M)  \,]\,(8)$. Hence, it follows from $\{ (1),\, (8) \}$ that  
$[\, \mathrm{\,G.dim}_R(M)= \mathrm{L.dim}_R(M)\,]$.\\
\underline{$\mathrm{Proof \, of \, (ii)}$} : It follows from 
$\{ \mathrm{\ref{dimensions2}\mathrm{.i}}\, , \ref{EKlazarus} \}$
\end{proof}

\begin{remark} Let $R$ be a ring (with identity) and $M\in R$-$\mathrm{Mod}$ and $\{\kappa, \lambda  \}\subseteq \mathrm{Card}$. \\ $\mathrm{(i)}$ In \ref{defgoldie}  the ring may be non-commutative and  $\mathrm{G}_R(M)$ is downward closed. Also $\mathrm{G.dim}_R(M)$ is called Goldie of uniform dimension, and is well defined because $\mathrm{G}_R(M)$ is bounded above, since $[\,\kappa\leq |M |\,]\, \forall \, \kappa \in \mathrm{G}_R(M)$. Also $[\, \kappa \in \mathrm{G.dim}_R(M)$ iff $M$ contains an internal direct sum of $\kappa$ submodules$\,]\,(2)$.\\
$\mathrm{(ii)}$  The condition $[\, R^{[\kappa]} \underset{R}{ \rightarrowtail} R^{[\lambda]},  \, \kappa< \omega, \, \lambda <\omega]\Rightarrow [\,\kappa \leq \lambda  \,]$, holds if $R$  is any commutative ring (see \ref{mccoy}) but does not hold in general
(see Remark 1.32 of \cite{lam}).  \\
$\mathrm{(iii)}$ If $[\, R^{[\kappa]} \underset{R}{ \cong} R^{[\lambda]},  \, \kappa< \omega, \, \lambda <\omega]\Rightarrow [\,\kappa = \lambda  \,]$, then we say that $R$ is $\mathsf{IBN}$
(Invarianvce basis number). So, in $\ref{dimensions1}\mathrm{.vi}$ we proved that left Noetherian rings are $\mathsf{IBN}$.\\
$\mathrm{(iv)}$  The proofs of $\{  \ref{dimensions1}\mathrm{.iv}, \, \ref{dimensions1}\mathrm{.v}, \, \ref{dimensions1}\mathrm{.vi}, \ref{dimensions1}\mathrm{.vii}\} $ have been adjusted from \cite{lam}.
But the argument used for the proof of $ \ref{dimensions1}\mathrm{.vi}$ in \cite{lam} (which is Lemma 1.36 of \cite{lam}) is different since it does not invoke Lemma 5 of  \cite{fuchs}. \\
$\mathrm{(v)}$ We are not aware of a simple proof of
$\ref{dimensions1}\mathrm{.ix}$. Note that the claim is trivial when
 $R$  is commutative, since  $(x,y)=(b,a)$ is an obvious solution.

\end{remark}

\section{Projective modules with isomorphic duals} \label{projectivemodules}

\begin{definition}  Let $R$ be a commutative ring, $M\in R$-$\mathrm{Mod}$ and $\kappa,\, \lambda$ be cardinals.
 Then $:$ \\   $\mathrm{(i)}$ Let $f\in \mathrm{Hom}_R(M,N)$.                     Then $f$  is $R$-dense (or a weak $R$-epimorphism) iff $(\mathrm{Cok}(f))^{\ast}=\bold{0}$. \\
 $\mathrm{(ii)}$  $(\, M  \overset{_w}{\underset{R} \twoheadrightarrow}  N\,) \Leftrightarrow$                            $(\, \exists\,   f\in \mathrm{Hom}_R(M,N) :f $ is $R$-dense$\,)$.\\
 $\mathrm{(iii)} \,$   $\mathrm{E}_R(M)=\{\mu \in \mathrm{Card}\, \big|\, R^{[\mu]}  \overset{_w}{\underset{R} \twoheadrightarrow} M  \} $.
\quad \quad\quad \quad 
$\mathrm{(iv)} \,\, \mathrm{E.dim}_R(M)=\mathrm{min}
(\mathrm{E}_R(M))$                         \\
 $\mathrm{(v)} $     $(\, M \overset{_s} {\underset{R} \rightarrowtail}  N\,) \Leftrightarrow$                            $(\, \exists\,   f\in \mathrm{Hom}_R(M,N) :f $ is $R$-left invertible$\,)$.\\          
 $\mathrm{(vi)} $  $\mathrm{K}_R(M)=\{\mu \in \mathrm{Card}\, \big|\,$  $  R^{[\mu]} \overset{_s} {\underset{R} \rightarrowtail }  M\}$.\quad \quad\quad \quad   $\mathrm{(vii)} \,\, \mathrm{K.dim}_R(M)=\mathrm{sup}
(\mathrm{K}_R(M))$  \quad \quad\quad \quad\quad \quad 
\end{definition}

\begin{definition}\label{strangedef}  Let $[\, R_k$ be commutative ring$\,]\, \forall \, k\in  \mathrm{T}_n$, $R=\overset{n}{\underset{k=1}\Pi} R_k$ and $[\, M_k\in R_k$-$\mathrm{Mod}\,]\, \forall \, k\in \mathrm{T}_n$. The abelian group
$\overset{n}{\underset{k=1} \oplus}  M_k $ becomes $R$-module by
$r\ast \overset{n}{\underset{k=1} \oplus}  m_k=\overset{n}{\underset{k=1} \oplus} (  \pi _k(r)\cdot m_k )$, which we denote by 
$\overset{n}{\underset{i=1} \boxplus}  M_k $. Let also   $[\, N_k\in R_k$-$\mathrm{Mod}\,]\, \forall \, k\in \mathrm{T}_n $, $[\,f_k\in \mathrm{Hom}_{R_k}(M_k,N_k)\,]\, \forall \, k\in \mathrm{T}_n $, $M=\overset{n}{\underset{k=1} \boxplus}  M_k$ and  $N=\overset{n}{\underset{k=1} \boxplus}  N_k$. The map $f=\overset{n}{\underset{k=1} \times}  f_k:M\rightarrow N $ defined by $f(\overset{n}{\underset{k=1} \oplus} m_k)=\overset{n}{\underset{k=1} \oplus}f_k(m_k)$ is $R$-linear
and will be denoted by $\overset{n}{\underset{k=1} \boxplus}  f_k $. Hence  $\overset{n}{\underset{k=1} \boxplus}  f_k\in \mathrm{Hom}_R(\overset{n}{\underset{k=1} \boxplus}  M_k, \overset{n}{\underset{k=1} \boxplus}  N_k)$.  
\end{definition}

\begin{lemma} \label{strangesum}  Under the above notation, let also $g\in \mathrm{Hom}_R(M,N)$ and $g_k=g_{\pi_k}:M_{\pi_k}\rightarrow  N_{\pi_k}$. Then $:$  \end{lemma} .\\
$\mathrm{(i)}$  $\mathrm{(a)}\,\, M_{\pi_k} \underset{R_k}\cong M_k $ \quad\quad
$\mathrm{(b)}\,\, f_{{\pi_k}}\underset{R_k}\sim f_k $ \quad\quad
$ \mathrm{(c)} \,\, g \underset{R}\sim   \overset{n}{\underset{k=1} \boxplus}  g_k $. \quad\quad$ \mathrm{(d)}\,\,[\,M=\bold{0}\,]\Leftrightarrow [\, M_k=\bold{0}\,]\, \forall \, k\in
 \mathrm{T}_n$. \\
$\mathrm{(ii)}$  $\mathrm{(a)}\,  \mathrm{Ker}(g)\underset{R}\cong  \overset{n}{\underset{k=1} \boxplus}   \mathrm{Ker}(g_k)$. \quad\quad\quad\quad\quad\quad\quad\quad
$\mathrm{(b)}\,[\, g$  is 1-1$\,]\Leftrightarrow [\, g_k$ is 1-1$\,]\, \forall \, k \in \mathrm{T}_n$. \\   
$\mathrm{(iii)}$  $\, \mathrm{(a)}\,\,\mathrm{Cok}(g)\underset{R}\cong  \overset{n}{\underset{k=1} \boxplus}   \mathrm{Cok}(g_k)$
 \quad\quad\quad\quad\quad\quad\quad\quad $\mathrm{(b)}\,[\, g$  is onto$\,]\Leftrightarrow [\, g_k$ is onto$\,]\, \forall \, k \in \mathrm{T}_n$. \\
$\mathrm{(iv)}$  $\mathrm{(a)}\,[\, M \underset{R}\rightarrowtail N \,] \Leftrightarrow[\, M_k \underset{R_k}\rightarrowtail N_k \,]\, \forall \, k\in \mathrm{T}_n$. \quad\quad
$\mathrm{(b)}\,[\, M \underset{R} \twoheadrightarrow N \,] \Leftrightarrow[\, M_k \underset{R_k} \twoheadrightarrow N_k \,]\, \forall \, k\in \mathrm{T}_n$.\\
$\mathrm{(v)}\,\,\mathrm{(a)}\, [\,g$  is bijection$\,]\Leftrightarrow [\, g_k$ is bijection$\,]\, \forall \, k \in \mathrm{T}_n$. \space
$\, \mathrm{(b)}\,\, [\, M \underset{R}\cong N \,] \Leftrightarrow[\, M_k \underset{R_k}\cong N_k \,]\, \forall \, k\in \mathrm{T}_n$.  \\ 
$\mathrm{(vi)}\, \mathrm{(a)}$ $[\, g$  is  
$R$-left invertible$\,]\Leftrightarrow [\, g_k$ is so$\,]\, \forall \, k \in \mathrm{T}_n$.
$\, \mathrm{(b)}\,[\, g$  is  
$R$-right invertible$\,]\Leftrightarrow [\, g_k$ is so$\,]\, \forall \, k \in \mathrm{T}_n$.\\
$\mathrm{(vii)}$ $[\, g$  is  
 $R$-dense $\,]\Leftrightarrow [\, g_k$ is $R$-dense $\,]\, \forall \, k \in \mathrm{T}_n$.\\
$\mathrm{(viii)}\,$  $\mathrm{(a)}\,\, \mathrm{Hom}_R(\overset{n}{\underset{k=1} \boxplus}  M_k, \overset{n}{\underset{k=1} \boxplus}  N_k)\underset{R}\cong \overset{n}{\underset{k=1} \boxplus} \mathrm{Hom}_{R_k}(M_k,N_k)$.
\quad\quad\quad  $\mathrm{(b)}\,$  $[\,(M_{\pi_k})^{\ast} \underset{R_k}\cong (M^{\ast})_{\pi_k} \,$.  \\
$\mathrm{(ix.a)}\,\,[\, M^{\ast} \underset{R} \rightarrowtail N^{\ast} \,]\Leftrightarrow[\, (M_k)^{\ast} \underset{R_k} \rightarrowtail (N_k)^{\ast} \,]\, \forall \, k\in \mathrm{T}_n $.  \space
$\mathrm{(ix.b)}\,\,[\, M^{\ast} \underset{R} \twoheadrightarrow N^{\ast} \,]\Leftrightarrow[\, (M_k)^{\ast} \underset{R_k} \twoheadrightarrow (N_k)^{\ast} \,]\, \forall \, k\in \mathrm{T}_n $. \\
$\mathrm{(ix.c)}\,\,[\, M^{\ast} \underset{R} \cong N^{\ast} \,]\Leftrightarrow[\, (M_k)^{\ast} \underset{R_k} \cong (N_k)^{\ast} \,]\, \forall \, k\in \mathrm{T}_n$. \quad $\mathrm{(ix.d)}\,\,[\, M^{\ast} =\bold{0} \,]\Leftrightarrow[\, (M_k)^{\ast} =\bold{0} \,]\, \forall \, k\in \mathrm{T}_n $. \\
$\mathrm{(ix.e)}\,\,[\, M^{\ast} \overset{_w} {\underset{R} \twoheadrightarrow} N^{\ast} \,]\Leftrightarrow[\, (M_k)^{\ast} \overset{_w}{\underset{R_k} \twoheadrightarrow} (N_k)^{\ast} \,]\, \forall \, k\in
 \mathrm{T}_n $. \space 
$\mathrm{(ix.f)}\,\,[\, M^{\ast}  \overset{_s} {\underset{R} \rightarrowtail} N^{\ast} \,]\Leftrightarrow[\, (M_k)^{\ast}  \overset{_s}{\underset{R_k} \rightarrowtail} (N_k)^{\ast} \,]\, \forall \, k\in \mathrm{T}_n $.  \\
$\mathrm{(x)}\, $ $\mu_R(M)=\mathrm{max}\{\mu_{R_i}(M_i)\,|\, i\in \mathrm{T}_n  \}$.          
 
\begin{proof}  We will prove  only $\mathrm{(ix.a)}$ and $\mathrm{(vii)}$  and leave the rest to the reader. \\
 \underline{$\mathrm{Proof\, of \ (vii)}$} :  $[\, g$  is  
 $R$-dense $\,]\overset{\text{def}}\Longleftrightarrow $   $[\, (\mathrm{Cok}(g))^{\ast}=\bold{0}\,]\overset{\mathrm{(iii.a)}}\Longleftrightarrow[\,\overset{n}{\underset{k=1} \boxplus}   (\mathrm{Cok}(g_k))^{\ast}=\bold{0} \,] $          
 $\overset{\mathrm{(i.d)}}\Longleftrightarrow  [\,  (\mathrm{Cok}(g_k))^{\ast}=\bold{0} \,]\, \forall k \in \mathrm{T}_n $   $\Longleftrightarrow[\, g_k$ is $R$-dense $\,]\, \forall \, k \, \in \mathrm{T}_n$.      \\
  \underline{$\mathrm{Proof\, of \ (ix.a)}$} :
$[\, M^{\ast} \underset{R} \rightarrowtail N^{\ast} \,]\overset{\mathrm{(iv.a)}}\Longleftrightarrow$    $[\, (M^{\ast})_{\pi_k} \underset{R_k} \rightarrowtail (N^{\ast})_{\pi_k} \,]\,\forall \, k\in \mathrm{T}_n \overset{\mathrm{(vii.b)}}\Longleftrightarrow$
 $[\, (M_{\pi_k})^{\ast} \underset{R_k} \rightarrowtail (N_{\pi_k})^{\ast}) \,]\,\forall \, k\in \mathrm{T}_n $
$\overset{\mathrm{(i.a)}}\Leftrightarrow$
 $[\, (M_{k})^{\ast} \underset{R_k} \rightarrowtail (N_{k})^{\ast}) \,]\,\forall \, k\in \mathrm{T}_n $ .  \end{proof}

We denote the class of reflexive $R$-modules by $\mathrm{Refl}(R)$.
\begin{lemma} \label{refllemma} Let  $R$ be a commutative ring and 
$\{M,\, N\}\subseteq \mathrm{Refl}(R)$. Then $ [\,  M^{\ast} \cong_R N^{\ast}\,]\Rightarrow [\,M   \cong_R N\,]$.
\end{lemma}\begin{proof} By assumption, we have that  $[\, \{ M,\, N \}\subseteq \mathrm{Refl}(R)\,]\Rightarrow$
$[\, M^{\ast\ast}\cong_R M,\, N^{\ast\ast}\cong_R N\,]\,(2) $. Therefore 
 $ [\,  M^{\ast} \cong_R N^{\ast}\,]\Rightarrow   [\,  M^{\ast\ast} \cong_R N^{\ast\ast}\,]\overset{(2)} \Rightarrow [\,M  \cong_R N\,]$. \end{proof}

\begin{definition} \label{vertical} Let $R$ be a commutative ring, 
$\{\kappa,\, \lambda\}\subseteq \mathrm{Card}$, $S \subseteq R^{[\lambda]}$ and $\{a,\,b\} \subseteq R^{[\lambda]}$. Then $:$ \\
$\mathrm{(i)}$ $\mathcal{h} a, b \mathcal{i}=\underset {i\in \lambda} \Sigma a_i \cdot b_i$ \quad\quad $\mathrm{(ii)}$ $ a \perp b\Leftrightarrow \, \mathcal{h} a, b \mathcal{i}=0$
\quad\quad $\mathrm{(iii)}\, S^{\perp}=\{ c\in R^{\lambda} \big|\, ( c \perp y)\, \forall \, y\in S \}$  \\
$\mathrm{(iv)}\, \xi_a \in \mathrm{Hom}_R(R^{[\lambda]},R)$ is defined by $\xi_a(x)=\mathcal{h} a, x \mathcal{i}$.\space $\mathrm{(v)}\, \hat{\xi}_a \in \mathrm{Hom}_R(R^{\lambda},R)$ is defined by $\hat{\xi}_a(x)=\mathcal{h} a, x \mathcal{i}$.\end{definition}

\begin{lemma}\label{refl-slenderness}  Let $R$ be a commutative ring and $\kappa\in \mathrm{Card}$. Then  $:$ \\
$\mathrm{(i)}$  $R^{[\kappa]}\in \mathrm{Refl}(R)\Leftrightarrow$                   
$\forall \, f\in (R^{\kappa})^{\ast} \exists \, a \in R^{[\kappa]} :
f=\hat{\xi}_a$ \space 
$\mathrm{(ii)}$
$R$ is slender  
$ \Leftrightarrow R^{[\omega]}\in \mathrm{Refl}(R)$. \\
$\mathrm{(ii)}$ If $\kappa \in  \mathcal{L}_1'$ then $:$  $\mathrm{(a)}$ $(R$ is slender$)$  
$ \Rightarrow( R^{[\kappa]}\in \mathrm{Refl}(R))$. \space
 $\mathrm{(b)}$ $(R$ is slender$)$  
$ \Rightarrow (R^{ \cdot\kappa})^{\ast}=\bold{0})$.\\
$\mathrm{(iii)}$ Assume that $\mathcal{L}_1' \neq \emptyset$. Then 
$\exists \, \kappa \in \mathcal{L}_1' : ( R^{[\kappa]}\notin \mathrm{Refl}(R))$.
\end{lemma}
\begin{proof}
See \cite{eklof}.
\end{proof}

\begin{lemma} \label{strangesum2}
Let $R_1,\, R_2\dots ,R_n$ be commutative rings, $R=\overset{n}{\underset{k=1}\Pi} R_k$ and $ \lambda,\,\mu$ be  cardinals
. Then \\
$\mathrm{(i)}\,$  
$ (R^{[\lambda]})_{\pi_k} \underset{R_k} \cong (R_k)^{[\lambda]}$  \quad
$\mathrm{(ii)}\,$  
$ (R^{\lambda})_{\pi_k} \underset{R_k} \cong (R_k)^{\lambda}$ \quad
$\mathrm{(iii)}\,$  
$ R^{\lambda} \underset{R} \cong  \overset{n}{\underset{k=1} \boxplus}  (R_k)^{\lambda}$  \quad\quad $\mathrm{(iv
)}\,$  $ R^{[\lambda]} \underset{R} \cong  \overset{n}{\underset{k=1} \boxplus}  (R_k)^{[\lambda]}$ \\
$\mathrm{(v)}\,$If $[\,  \dim_{R_k} (R_k)^{\lambda}=\mu\,]\, \forall\, k\in \mathrm{T_n}$ then
$[\, \dim_{R} R^{\lambda}=\mu\,]$. $($Reformulation of  \ref{strange} $)$. 
\end{lemma}
\begin{proof} $\mathrm{(i)}$ follows from $\ref{rh}\mathrm{.i}$ and  
 $\mathrm{(ii)}$ follows from $\ref{rh}\mathrm{.ii}$ since $\pi_k$'s are f.p.\\Also   $\mathrm{(iii)}$ follows from $\{\ref{strangesum}\mathrm{.i}, \mathrm{(ii)} \}$ and   $\mathrm{(iv)}$ follows from $\{\ref{strangesum}\mathrm{.i}, \mathrm{(i)} \}$ .  \\
  \underline{$\mathrm{Proof\, of \ (v)}$} :
$[\,  \dim_{R_k} (R_k)^{\lambda}=\mu\,]\, \forall\, k\in \mathrm{T_n}\Rightarrow$
$[\,  (R_k)^{\lambda}\cong_{R_k} (R_k)^{[ \mu]}\,]\, \forall\, k\in \mathrm{T_n}\Rightarrow$
$ [\, \overset{n}{\underset{k=1} \boxplus}  (R_k)^{\lambda}\underset{R}\cong   \overset{n}{\underset{k=1} \boxplus}  (R_k)^{[\mu]} \,]$ 
$\overset{\mathrm{(iii)},\mathrm{(iv)}}\Longrightarrow$
$[\, R^{\lambda} \cong_R \, R^{[\mu]}]\Rightarrow$
$[\, \dim_{R} R^{\lambda}=\mu\,]$.
\end{proof}
\begin{lemma} \label{strangesum3}   Let $[\, R_k$ be commutative ring$\,]\, \forall \, k\in  \mathrm{T}_n$, $R=\overset{n}{\underset{k=1}\Pi} R_k$ and $[\, M_k\in R_k$-$\mathrm{Mod}\,]\, \forall \, k\in \mathrm{T}_n$. Then
$\mathrm{(i)}\,[\,M \in \mathrm{Proj}(R)\,]\Leftrightarrow[\,  M_k \in \mathrm{Proj}(R_k)\,]\, \forall \, k\in \mathrm{T}_n $. 
$\mathrm{(ii)}\,$ $[\, M \in \mathrm{Refl}(R)\,]\Leftrightarrow[\,  M_k \in \mathrm{Refl}(R_k)\,]\, \forall \, k\in \mathrm{T}_n $.\\
$\mathrm{(iii)}\,[\,M \in \mathrm{Inj}(R)\,]\Leftrightarrow[\,  M_k \in \mathrm{Inj}(R_k)\,]\, \forall \, k\in \mathrm{T}_n $. 
$\mathrm{(iv)}\,  R$ is self-injective $\Leftrightarrow[\,  R_k$ is self-injective $\,]\, \forall \, k\in \mathrm{T}_n $.\\
$\mathrm{(v)}\, [\, R$ is reduced$\,]\Leftrightarrow[\,  R_k$ is reduced$\,]\, \forall \, k\in \mathrm{T}_n $ 
\end{lemma}
\begin{proof}\underline{$\mathrm{Proof\, of \ (i)}$} : If 
$ [\, M \in \mathrm{Proj}(R)\,]$ then $M_k \underset{R_k}\cong M_{\pi_k} \overset{\ref{remainprojective}}\in \mathrm{Proj}(R_k)$.
Conversely, assume that  $[\,  M_k \in \mathrm{Proj}(R_k)\,]\, \forall \, k\in \mathrm{T}_n $. Let $ k\in \mathrm{T}_n $. Then
$[\,M_k\oplus N_k \cong_{R_k} R^{[\lambda_k]}\,]\,(5) $ for some  $ N_k \in \mathrm{R_k}$-$\mathrm{Mod}$ and $ \lambda_k\in \mathrm{Card}$. Let $\lambda=\mathrm{max}\{\lambda_k\,|\,k\in \mathrm{T}_n \}$. Then $[\,\lambda_k +\mu_k=\lambda\,]\,(6)$ for some   
$\mu_k \in \mathrm{Card}$. Let $[\, L_k=(N_k \oplus R^{[\mu_k]})\,]\,(7)$. It follows from $\{(5),(6),(7) \}$ that $[\,M_k \oplus L_k \cong_{R_k} (R_k)^{[\lambda]} \,]\Rightarrow$
$[\, (\overset{n}{\underset{i=1} \boxplus}  M_k) \oplus (\overset{n}{\underset{i=1} \boxplus}  L_k) \underset{R}\cong \overset{n}{\underset{i=1} \boxplus} (R_k)^{[\lambda]}\overset{\ref{strange2}\mathrm{.iv}} {\cong}_R R^{[\lambda]}\,]\Rightarrow$
$[\, (\overset{n}{\underset{i=1} \boxplus}  M_k)\in \mathrm{Proj}(R)\,]\Rightarrow [\, M \in \mathrm{Proj}(R)\,]$ \\
\underline{$\mathrm{Proof\, of \ (ii)}$} : Consider the maps $\sigma_M: M\rightarrow M^{\ast\ast}$ and $\sigma_{M_k}: M_k\rightarrow (M_k)^{\ast\ast}$. It is not hard to check that $[\,\sigma_M \underset{R}\sim  \overset{n}{\underset{i=1} \boxplus} \sigma_{M_k}\,]\,(9)$. Hence  $[\,  M \in \mathrm{Refl}(R)\,]\Leftrightarrow[\,\sigma_M $ is bijection$\,]$ $\underset{(9)}{\overset{\ref{strangesum}}\Longleftrightarrow}$
$[\,\sigma_{M_k}$ is bijection$\,]\, \forall \,k\in \mathrm{T}_n \Leftrightarrow 
[\,  M_k \in \mathrm{Refl}(R_k)\,]\, \forall \, k\in \mathrm{T}_n $.\\
\underline{$\mathrm{Proof\, of \ (iii)}$} : Left to the reader. \quad\quad\quad\quad\quad\quad\quad\quad\quad
\underline{$\mathrm{Proof\, of \ (iv)}$} : Follows from $\mathrm{(iii)}$.
\end{proof}

\begin{theorem} \label{slenderproduct}
 Let $R_1,\, R_2\dots ,R_n$ be commutative rings, $R=\overset{n}{\underset{k=1}\Pi} R_k$ and $ \lambda$ be a  cardinal. Then \\
 $\mathrm{(i)}\,[\, R^{[\lambda]}\in \mathrm{Refl}(R)\,]\Leftrightarrow [\, (R_{k})^{[\lambda]}\in \mathrm{Refl}(R_k)\,]\, \forall \,  k \in \mathrm{T}_n$. \\ $\mathrm{(ii)}\,[\, R$ is slender$\,]\Leftrightarrow [\, R_k$ is slender$\,]\, \forall\, k\in \mathrm{T}_n$. \\ $\mathrm{(iii)}$   $[\,R$ is  
 $\lambda$-half-slender$\,]\Leftrightarrow [\, R_k$ is 
 $\lambda$-half-slender$\,]\, \forall\, k \in \mathrm{T}_n$.  \\
 $\mathrm{(iv)}$   $[\, \Phi^{\lambda}_R$ is $R$-left invertible $\,]\Leftrightarrow$
$[\, \Phi^{\lambda}_{R_k}$ is $R_k$-left invertible $\,]\, \forall \, k\in \mathrm{T}_n$. 
\end{theorem}
\begin{proof}
$\ref{slenderproduct}\mathrm{.i}$ follows from $\{\ref{strangesum3}\mathrm{.ii},\,  \ref{strangesum2}\mathrm{.iv} \}$ and $\ref{slenderproduct}\mathrm{.ii}$ follows from 
$\{\ref{refl-slenderness}, \, \ref{slenderproduct}\mathrm{.i}\}$.\\
\underline{$\mathrm{Proof\, of \ (iii)}$} : Consider the map $\Phi^{\lambda}_R:R^{[\lambda]}\rightarrow R^{\lambda}$ defined in 
$\ref{def.half.slender}$ and $\pi_k : R\rightarrow R_k$. It follows from  $\ref{rh}$ that 
$[\, (\Phi^{\lambda}_R)_{\pi_k} \underset{R_k}\sim \Phi^{\lambda}_{R_k})\,]\,(1)$. So  $[\,R$ is  $\lambda$-half-slender$\,]\overset{\ref{def.half.slender}}\Leftrightarrow$
$[\, \Phi^{\lambda}_R$ is $R$-dense$\,]$   
$\overset{\ref{strangesum}\mathrm{.vii}}\Longleftrightarrow$
$[\, (\Phi^{\lambda}_{R})_{\pi_k}$ is $R_k$-dense $\,]\, \forall \, k\in \mathrm{T}_n$ \\ $\underset{(1)}\Leftrightarrow$
$[\, \Phi^{\lambda}_{R_k}$ is $R_k$-dense $\,]\, \forall \, k\in \mathrm{T}_n$
$\overset{\ref{def.half.slender}} \Leftrightarrow [\, R_k$ is 
$\lambda$-half-slender$\,]\, \forall\, k \in \mathrm{T}_n$. \\
\underline{$\mathrm{Proof\, of \ (iv)}$} : 
$[\, \Phi^{\kappa}_R$ is $R$-left invertible$\,]$
$\overset{\ref{strangesum}\mathrm{.vi}}\Longleftrightarrow$
$[\, (\Phi^{\lambda}_{R})_{\pi_k}$ is $R_k$-left invertible $\,]\, \forall \, k\in \mathrm{T}_n$    
$\underset{(1)} \Leftrightarrow$
$[\, \Phi^{\lambda}_{R_k}$ is $R_k$-left invertible$\,]\, \forall \, k\in \mathrm{T}_n$.

\end{proof}

\begin{lemma} \label{ontodual}Let $M,\,N,$ be modules over a  ring $R$, $f\in \mathrm{Hom}_R(M,N)$ and $f^{\ast}=\mathrm{Hom}_R(f,R)$. Then $:$
$\mathrm{(i)}\,[\,(\mathrm{Cok}(f))^{\ast}=\bold{0}\,]\Rightarrow[\, f^{\ast}$  $ \text{is 1-1} \,]$ \quad $\mathrm{(ii)}\,[\,f $ is onto$\,]\Rightarrow[\, f^{\ast}$  $ \text{is 1-1} \, ]$ \quad
 $\mathrm{(iii)}\,[\, M \underset{R} \twoheadrightarrow N\,]\Rightarrow[\, N^{\ast} \underset{R}\rightarrowtail M^{\ast}\,]$ 
\end{lemma}
\begin{proof} \underline{$\mathrm{Proof\, of \ (i)}$} :
By assumption, we have that $[\,  (\mathrm{Cok}(f))^{\ast}=\mathrm{Hom}_R((N/ \mathrm{Im}(f)),R   )=\bold{0}\,]\,(1)$. Trivially  $[\, f^{\ast}\in \mathrm{Hom}_R(N^{\ast},
M^{\ast})\,]\,(2)$. Let $[\, h\in \mathrm{Ker}(f^{\ast})\,]\Rightarrow[\,f^{\ast}(h)=h\circ f=0\,]\Rightarrow $ 
$[\, (h(t)=0)\, \forall \, t\in \mathrm{Im}(f) \,]\overset{(1)}\Rightarrow[\, f=0\,]$. Hence                                               
$[\, \mathrm{Ker}(f^{\ast})=\bold{0}\,]\Rightarrow$
$[\, f^{\ast}$ is 1-1$\,]\,(3)$.   
Also $\mathrm{(ii)}$ follows $\mathrm{(i)}$ and $\mathrm{(iii)}$
follows from $\{(2).\, (3)  \}$. \end{proof}

\begin{theorem} 
Let $R$  be a Noetherian commutative ring and $\kappa,\, \lambda$ be cardinals.  Assume  that either   $\mathrm{(i)}\,[\,  \exists \,\mathfrak{m} \in \tmop{max(R)}: \big| \dfrac{R}{\mathfrak{m}} \big| \le \mathfrak{c}, \, \mathsf{ICF} \,] $ or  $ \mathrm{(ii)} \,[\, |R|,\, \lambda,\,  \kappa $ are all non-$\omega$-measurable  and  $R$ is $\mathrm{w}$-slender$\,]$ or 
 $ \mathrm{(iii)} \,[\,  |R|,\, 
 \lambda $  are both weakly-reachable and  $R$ is $\mathrm{w}$-slender$\,]$.  Then $[\,  R^{\lambda} \underset{R}\twoheadrightarrow R^{\kappa}\,]\,(1) \Leftrightarrow [\, \kappa \leq \lambda \,]\,(2)$
 \end{theorem}
 \begin{proof} $\bullet$ The direction $[\, (2)\Rightarrow (1)\,]$ is  $\ref{trivial1}\mathrm{.iv.a}$. \quad\quad  \\
 $\bullet$ Proof of   $[\, (1)\Rightarrow (2)\,]$  : \quad
 It follows from  (1) that $[\, \exists \, f\in \mathrm{Hom}_R(R^{\lambda},R^{\kappa}) $  such that $ f$ is onto$\,]$ (3).\\
 \underline{$\mathrm{Case \, (i)}$} :  Let  $\phi=\pi^R_\mathfrak{m}:R\rightarrow \dfrac{R}{\mathfrak{m}}$ . It follows by assumption that  $[\, \dfrac{R}{\mathfrak{m}}=F$ is a field$\,]\,(4)$, $[\, |F |\le \mathfrak{c} \,]$ (5) and    $[\, \phi$ is a f.p $]$ (6), since $R$ is Noetherian. So  $[\,  R^{\lambda} \underset{R}\twoheadrightarrow R^{\kappa}\,]\underset{(6)} {\overset{\ref{rh2}\mathrm{.vi.a}}\Longrightarrow} [\,   F^{\lambda} \underset{F}\twoheadrightarrow F^{\kappa}\,]\underset{(4)}\Rightarrow$
$[\, \dim_F(F^{\lambda})\geq \dim_F(F^{\kappa})  \,]$                         
$\underset{(5)}  {\overset{ \ref{EKextension}\mathrm{.iv}} \Longrightarrow }[\,2^{\, \lambda} \geq 2^{\, \kappa}\,] $ 
$\underset{\mathsf{ICF}}{\overset{\ref{isomorphicduals}\mathrm{.x}}\Longrightarrow}[\, \lambda \geq \kappa \,]\Rightarrow [\, \kappa \leq \lambda \,]$. \\
  \underline{$\mathrm{Case \, (ii)}$} :  There is a r.h $\psi:R\rightarrow S$  such that $[\, \psi$ is f.p$\,]$ (15) and $[\, S$                  is slender$\,]$ (16), since $R$ is w-slender.  Let $f_{\psi}=\mathrm{id}_F \underset{R}\otimes f$. It follows from
 $\{ \mathrm{(ii)},\, \ref{measurable}\mathrm{.ii}   \}$ that $[\, (S^{\kappa})^{\ast}\underset{S}\cong S^{[\kappa]}\,]\,(17)$ and  $[\, (S^{\lambda})^{\ast}\underset{S}\cong S^{[\lambda]}\,]\,(18)$.
 It follows from (3) that $[\, f_{\psi}$ is onto$\,]\Rightarrow$
 $[\,  S^{\lambda} \underset{S}\twoheadrightarrow S^{\kappa}\,] \overset{\ref{ontodual}\mathrm{.ii}}\Longrightarrow$ 
  $[\,  (S^{\kappa})^{\ast} \underset{S}  \rightarrowtail (S^{\lambda})^{\ast}\,] \overset{(17),(18)}\Longrightarrow$ 
 $[\,  S^{[\kappa]} \underset{S}  \rightarrowtail S^{[\lambda}]\,]\overset{\ref{innocuous}}\Rightarrow  [\, \kappa \leq 
 \lambda  \,]$.\\
 \underline{$\mathrm{Case \, (iii)}$} : The proof is similar to the proof of $\mathrm{Case \, (ii)}$. The only difference is the proof of $\{ (17),(18) \}$. \end{proof}
In the next lemma  we deal with left modules over rings that may not be commutative.
\begin{lemma} \label{noncommutativeIBN}Let $R$ be a   ring $($with identity$)$ and $\kappa,\lambda$ be cardinals. Assume that $[\, \lambda \geq \omega\,]\,(1)$. Then $:$ \\
$\mathrm{(i)}\, [\, R^{[\kappa]} {\underset{R}\twoheadrightarrow} R^{[\lambda]}\,]\Rightarrow [\, \kappa \geq \lambda\,]$ \quad
$\mathrm{(ii)}\, [\, R^{[\kappa]}  {\underset{R}\cong } R^{[\lambda}]\,]\Rightarrow [\, \kappa = \lambda\,]$ \quad
$\mathrm{(iii)}$ Both $R^{[\lambda]},\, R^{\lambda}$ are not f.g 
left $R$-modules.
\end{lemma}
\begin{proof}
\underline{$\mathrm{Proof\, of \, (i)}$} : $ [\, R^{[\kappa]}\underset{R}\twoheadrightarrow R^{[\lambda]}]\,]\Rightarrow$
$[\, \exists \, f\in \mathrm{Hom}_R(R^{[\kappa]},R^{[\lambda]}): f$
is onto$\,]\Rightarrow [\, R^{[\lambda]}=\underset{i\in \kappa}\Sigma R\cdot f(e_i)\,]\Rightarrow $ 
$[\,\lambda \subseteq \underset{i\in \kappa}  \cup \mathrm{supp}(f(e_i)) \,] \Rightarrow[\, \lambda \leq  \underset{i\in \kappa}  \Sigma | \mathrm{supp}(f(e_i))|\,\,]\,(2)$. Note that $ (3)\,[\, |\, \mathrm{supp}(f(e_i))|\, <\omega\,]\, \forall \, i\in \kappa $ .  If we assume that  $[\, \kappa < \omega\,]\,(4)$, then it follows from  $\{(2),(3),(4)  \}$ that $[\,\lambda <\omega  \,]$, which contradicts $(1)$. Therefore  
$[\, \kappa \geq \omega\,]\,(4)$. Hence 
$[\, \{ (2),(3) \}\,] \Rightarrow[\,  \lambda \leq \kappa \cdot \omega \overset{(4)}=\kappa \,]\Rightarrow [\, \lambda \leq  \kappa\,]\Rightarrow[\, \kappa \geq \lambda \,]$.\\
\underline{$\mathrm{Proof\, of \, (ii)}$} : $[\, R^{[\kappa]}\underset{R}\cong R^{[\lambda]} \,] \overset{(1)}\Rightarrow$
$[\, R^{[\kappa]}\underset{R}\twoheadrightarrow R^{[\lambda]}\,, \lambda\geq \omega]\overset{\mathrm{(a)}}\Rightarrow[\, \kappa \geq \lambda \geq \omega  \,]\,(5)$. \\ Also
$[\, R^{[\kappa]}\underset{R}\cong R^{[\lambda]}\,]\Rightarrow$
$[\, R^{[\lambda]}\underset{R}\twoheadrightarrow R^{[\kappa]}\,]\overset{\mathrm{(5),(a)}} \Longrightarrow[\, \lambda \geq \kappa  \,]\overset{(5)}\Rightarrow [\, \kappa=\lambda\,]$.\\
\underline{$\mathrm{Proof\, of \, (iii)}$} : Let $\mu=\mu_R(R^{[\lambda]})$. Then  $ [\, R^{[\mu]} {\underset{R}\twoheadrightarrow} R^{[\lambda]}\,]\overset{\mathrm{(i)}}\Rightarrow [\, \mu \geq \lambda\,]\overset{(1)}\Rightarrow [\, \mu \geq \omega\,]\Rightarrow [\, R^{[\lambda]}$  is non-f.g $R$-module$\,]$.\\
 Let $\xi=\mu_R(R^{\lambda})$. Then  $ [\, R^{[\xi]} {\underset{R}\twoheadrightarrow} R^{\lambda} {\underset{R}\twoheadrightarrow} R^{[\lambda]}\,] \Rightarrow [\,  R^{[\xi]} {\underset{R}\twoheadrightarrow} R^{[\lambda]}  \,] \overset{\mathrm{(i)}}\Rightarrow [\, \xi \geq \lambda\,]\overset{(1)}\Rightarrow [\, \xi \geq \omega\,]\Rightarrow [\, R^{\lambda}$  is a  non-f.g $R$-module$\,]$. \end{proof}
\begin{remark} .\\ $\mathrm{(i)}$ In \ref{noncommutativeIBN} the ring $R$ may be non-commutative. The commutative case has been covered by \ref{rankiswelldefined} and \ref{noninnocuous}.   \\
$\mathrm{(ii)}$ The proof of \ref{noncommutativeIBN} above has been adapted from  \cite{lam} and was the main motivation behind the proof of Theorem \ref{monoslender} below.\\
$\mathrm{(iii)}$ The assumption $[\, \lambda \geq \omega   \,]$ in
$\ref{noncommutativeIBN}\mathrm{.ii}$ is necessary, according to \cite{lam}, but it becomes redundant if $R$ is left Noetherian. See
\cite{lam} or Theorem \ref{noncommutative} below.
\end{remark} 
\begin{theorem} \label{noncommutative}Let $R$ be  left Noetherian ring   ring and $\kappa,\lambda$ be cardinals.  Then $ [\, R^{[\kappa]}   \cong_R  R^{[\lambda]}\,]\Rightarrow [\, \kappa = \lambda\,]$. \end{theorem}
\begin{proof}
$\cdot$ Let $[\, \lambda \geq \omega\,]$ then the claim follows from $\ref{noncommutativeIBN}\mathrm{.ii}$ . \\
$\cdot$ Let $[\, \lambda < \omega\,]$. It  follows from $\ref{noncommutativeIBN}\mathrm{.iii}$ that $[\,\kappa <\omega, \, \lambda <\omega  \,]$. Hence the claim follows from $\ref{dimensions1}\mathrm{.vi}$.
\end{proof}

\begin{theorem} \label{monoslender} $^{\ast}$
Let $R$ be a  commutative ring and $\kappa,\lambda$ be cardinals.
Asuume that either $\mathrm{(i)}\,  [\, \big| R \big| \leq \mathfrak{c},\, \mathsf{ICF} \,]$ or $\mathrm{(ii)} \,[\,  R $ is slender$\,]$. Then \space  $[\, R^{\kappa} \underset{R} \rightarrowtail  R^{\lambda}\,]\,(1) \Leftrightarrow[\, \kappa \leq \lambda \,]\,(2)$
\end{theorem}
\begin{proof} $\bullet$ The direction $[\, (2)\Rightarrow (1)\,]$ follows from $\ref{trivial1}$. \quad\quad \quad\quad \quad\quad  
$\bullet$ Proof of   $[\, (1)\Rightarrow (2)\,]$  : \quad \\
 \underline{$\mathrm{Case \, (i)}$} : Let $\mu=|R|$. It follows from  
 $\{\mathrm{(i)},\, \ref{cardinallemma}\mathrm{.iii}\}$ that $[\, \mu^{\,\kappa}=2^{\,\kappa}\,]\,(3)$ and  $[\, \mu^{\,\lambda}=2^{\,\lambda}\,]\,(4)$. Therefore 
   $[\, R^{\kappa} \underset{R} \rightarrowtail  R^{\lambda}\,]\Rightarrow$
 $[\, | R^{ \kappa}| \leq | R^{ \lambda}|  \,]\Rightarrow [\,\mu^{\kappa} \leq \mu ^{\, \lambda}\,]\overset{(3),(4)}\Longrightarrow[\,  2^{\,\kappa} \leq 2^{\, \lambda}\,]$
 $\underset{\mathsf{ICF}} {\overset{\ref{isomorphicduals}\mathrm{.x}} \Longrightarrow }$ 
 $[\, \kappa \leq \lambda \,]$  \\
 \underline{$\mathrm{Case \, (ii)}$} : 
  $[\, R^{\kappa} \underset{R} \rightarrowtail  R^{\lambda}\,]\Rightarrow [\, \exists \, f\in \mathrm{Hom}_R(R^{\kappa},R^{\lambda}):f $ is 1-1$\,]\,(6)$. But 
  $f=\underset{i\in \lambda}\Pi f_i$, for some $f_i \in  \mathrm{Hom}_R(R^{\kappa},R)$. From $\{\ref{slenderdef2},\mathrm{(ii)}\}$, we get that  $(7)\, [\, |\mathfrak{t}_{\kappa}(f_i)|\leq \omega\,]\,  \forall \, i \in \lambda$.  Let $[\,  A= \underset {i\in \lambda}\cup ( \mathfrak{t}_{\kappa}(f_i))\subseteq \kappa\,]\,(8)$. \\ Assume that  $[\, A\neq \kappa\,]\,(9)$. Let
  $x\in R^{\kappa}$ defined by $\{[\,  x(j)=0, \, \,  \text{if\, }  j\in \kappa  \,],\, [\,  x(j)=1, \,   \text{if\, }  j\in (\kappa \smallsetminus A ) \,] \}\, (10)$. It follows from $\{ (8),(10) \}$  that  $[\, f(x)=0\,] \,(11)$  and from (9) that  $[\, x\neq 0\,]\,(12)$. So we have reached a contradiction by $\{(6),(11),(12) \}$. So $[\,\kappa=A=  \underset {i\in \lambda}\cup ( \mathfrak{t}_{\kappa}(f_i))  \,]\Rightarrow$
$[\, \kappa \leq \underset{i\in \lambda}\Sigma   |\mathfrak{t}_{\kappa}(f_i)| \overset{(7)}\leq \lambda \cdot \omega  \,]\Rightarrow [\, \kappa \leq \lambda \cdot \omega\,]\,(13)$. \\
$\cdot$ Let $[\, \lambda\geq \omega\,]\,(14)$. Then 
$(13) \Rightarrow[\,  \kappa \leq \lambda \cdot \omega \overset{(14)}=\lambda \,]\Rightarrow [\, \kappa \leq \lambda\,]$. \\
$\cdot$ Let $[\, \lambda < \omega\,]\,(15)$. Then
$[\, R^{[\kappa]} \underset{R}\rightarrowtail R^{\kappa} \overset{(1)} {\underset{R}\rightarrowtail} R^{\lambda}\overset{(15)}=R^{[\lambda]}\,]\Rightarrow$ 
$[\, R^{[\kappa]} \underset{R}\rightarrowtail R^{[\lambda]}\,]\overset{\ref{innocuous}}\Rightarrow[\, \kappa \leq \lambda\,]$.\\
So we proved the claim in any case.  \end{proof} 

\begin{theorem} \label{w-slender4} $^{\ast}$ Let $R$  be a commutative ring and $\kappa,\, \lambda$
be cardinals.\\
$\mathrm{(i)}\,$  If $R$ is slender then 
 \space  $[\, R^{\kappa}  \cong_R  R^{\lambda}\,]  \Rightarrow[\, \kappa = \lambda \,]$ \quad
 $\mathrm{(ii)}\,$ If $R$ is w-slender then 
   $[\, R^{\kappa}  \cong_R  R^{\lambda}\,] \Leftrightarrow[\, \kappa = \lambda \,]$ \end{theorem} 
\begin{proof}\underline{$\mathrm{Proof}$ of $\mathrm{(i)}$} :
 $[\, R^{\kappa} \underset{R} \cong  R^{\lambda}\,] \Rightarrow[\, R^{\kappa} \underset{R} \rightarrowtail  R^{\lambda},\, R^{\lambda} \underset{R} \rightarrowtail  R^{\kappa} \,]\overset{\ref{monoslender}}\Rightarrow$
 $[\,\kappa \leq \lambda\,\,  , \lambda \leq \kappa \,]\Rightarrow [\, \kappa=\lambda\,]$. \\
\underline{$\mathrm{Proof}$ of $\mathrm{(ii)}$} : There is a r.h $\phi:R\rightarrow S$ such that $[\, S$ is slender$\,]\,(1)$ and $[\, \phi$ is f.p$\,]\,(2)$, since $R$
is a w-slender ring. Therefore  $[\, R^{\kappa}  \cong_R  R^{\lambda}\,] \underset{(2)}{\overset{\ref{rh2} } \Longrightarrow}[\, S^{\kappa}  \cong_S  S^{\lambda}\,] $
$\overset{\mathrm{(i)},(1)}\Longrightarrow[\, \kappa=\lambda \,]$\end{proof}
\begin{remark} \label{w-slenderremark} We can reformulate the above theorem like below: \\
Let $R$ be a w-slender commutative ring and $ \{ P, \, Q  \} \subseteq \mathrm{Free}(R)$. Then  $ [\,  P^{\ast} \cong_R Q^{\ast}\,]\Leftrightarrow [\,P  \cong_R Q\,] $ 
\end{remark}

\begin{theorem}Let $R$ be a commutative ring such that either
 $\mathrm{(i)}\,[\,  \exists \,\mathfrak{m} \in \tmop{max(R)}: \big| \dfrac{R}{\mathfrak{m}} \big| \le \mathfrak{c}, \, \mathsf{ICF}\,]  $  or $\mathrm{(ii)}[\, R$ is w-slender$\,]$. Then  $\mathsf{L}_R$ holds.
\end{theorem}
\begin{proof}Apply Definition \ref{definitions2}. So $\mathrm{case\, (i)}$ follows from \ref{secondtheorem2} and $\mathrm{case\, (ii)}$ folows from $\ref{w-slender4}\mathrm{.ii}$.
\end{proof}
 
\begin{lemma} \label{reflexivesum} Let $M_1,M_2\dots,M_n$ be modules over a commutative ring $R$. Then $:$ \\
$\mathrm{(i)}\,  (\overset{n}{\underset{k=1} \oplus}M_k)\in \mathrm{Proj}(R)\Leftrightarrow $
$[\, M_k \in \mathrm{Proj}(R_k)\,]\, \forall k\in \mathrm{T}_n$
 $\mathrm{(ii)}\,  (\overset{n}{\underset{k=1} \oplus}M_k)\in \mathrm{Refl}(R)\Leftrightarrow $
$[\, M_k \in \mathrm{Refl}(R_k)\,]\, \forall k\in \mathrm{T}_n$
\end{lemma}

\begin{lemma} \label{projective}Let $R$ be a Noetherian commutative ring, $ \{P, Q\} \subseteq \mathrm{Proj}(R)$ and $\kappa=\mu_R(P),\,\lambda=\mu_R(Q) $.   \end{lemma}
\begin{enumerate} [label=(\roman*)] 
\item     $ \exists\, N \in R$-$\mathrm{Mod} :[\,( P\oplus N )\cong_R R^{[\kappa]}\,] $\\
$\mathbf{Proof :}$
 Let $\kappa=\mu_R(P)$. Then  $[\, R^{[\kappa]} \underset{R}\twoheadrightarrow P,\,P\in \mathrm{Proj}(R)\,]\Rightarrow \exists\, Q \in R$-$\mathrm{Mod} :[\,( P\oplus Q)\cong_R R^{[\kappa]}\,]\, $.
\item If $\kappa < \omega$ (i.e $P$ is f.g) then : $\mathrm{(a)}\,[\,   P\in \mathrm{Refl}(R)\,]$ \quad   $\mathrm{(b)}\,[\,   P^{\ast}\in \mathrm{Proj}(R)\,]$  \\
$\mathbf{Proof :}$ $[\, \kappa <\omega \,]\Rightarrow [\, R^{[\kappa]}\in \mathrm{Refl}(R)\,] $
$\overset{\ref{reflexivesum}} {\underset{\mathrm{(i)}} \Longrightarrow } [\, P\in \mathrm{Refl}(R)\,]\, \mathrm{(a)}$.\\
 Also $\{\mathrm{(i),\, \kappa< \omega}\}\Rightarrow$
$[\, (P^{\ast}\oplus N^{\ast}) \cong_R R^{[\kappa]}\,]\Rightarrow [\, P^{\ast}\in \mathrm{Proj}(R)\,]\, \mathrm{(b)}$. 
\item  If $[ \, |R|, \mu_R(P)$  are not $\omega$-measurable, $R$ is slender$]\, (3)$. Then $ [\, P\in \mathrm{Refl}(R) \,]$   \\
$\mathbf{Proof :}$ It follows from $\{(3),\, \ref{measurablefacts}\mathrm{.ii} \}$ that   $[\, R^{[\kappa]}\in \mathrm{Refl}(R)\,] $
$\overset{\ref{reflexivesum}} {\underset{\mathrm{(i)}} \Longrightarrow } [\, P\in \mathrm{Refl}(R)\,]$. 
\item  If $[ \, |R|, \, \mu_R(P), \, \mu_R(Q)$  are not $\omega$-measurable, $R$ is slender$]\, (4)$. Then $ [\,  P^{\ast} \underset{R}\cong Q^{\ast}\,]\Rightarrow [\,P   \underset{R}\cong Q\,] $ \\
$\mathbf{Proof :}$ It follows from $\{(4), \mathrm{(iii)}  \}$
that $[\,  \{ P, \, Q  \} \subseteq \mathrm{Refl}(R)\,] \,(5)$. The claim follows from $\{(5), \mathrm{\ref{refllemma}}  \}$. 
\item   $\mathrm{(a)}\,$ Let $M\in R$-$ \mathrm{Mod}$. Then   $\, [\,\mu_R(M)< \omega \,]\Rightarrow [\, \mu_R(M^{\ast})<\omega\,]$ \space 
 $\mathrm{(b)}$  $[\, \mu_R(P)< \omega \,] \Leftarrow [\, \mu_R(P^{\ast})<\omega\,]$ \\
\underline{$\mathrm{Proof \,  of \,(a)}$} : Let $[\, \nu=\mu_R(M)< \omega\,] $ (6). Then $[\, R^{[\nu]}\underset{R}\twoheadrightarrow M\,]\overset{\ref{ontodual}}\Rightarrow$
$[\, M^{\ast} \underset{R} \rightarrowtail   (R^{[\nu]})^{\ast}\overset{(6)}{\underset{R}\cong}  R^{[\nu]} \,]\Rightarrow$
$[\, M^{\ast} \underset{R} \rightarrowtail R^{[\nu]}\,]\,(7) $.
The $R$-module $[\, R^{[\nu]}$ is Noetherian$\,]\,(8)$, since $[\, R$     
is Noetherian and $\nu <\omega\,]$. It follows from $\{ (7),\, (8) \}$
that $ [\, \mu_R(M^{\ast})<\omega\,]$. \\ 
\underline{$\mathrm{Proof \,  of \,(b)_{}}$}: Let $\kappa''=\mu_R(P^{\ast\ast})$. By assumption 
$[\, \mu_R(P^{\ast})<\omega\,]\overset{\mathrm{(a)}}\Rightarrow$
$[\,\kappa''= \mu_R(P^{\ast\ast})<\omega\,]\,(9)$. It follows from $\mathrm{(i)}$ that $[\, P$ is a torsionless $R$-module$\,]$ hence $ [\, P \underset{R} \rightarrowtail P^{\ast \ast}\,]\,(10)$. It follows from $\mathrm{(i)}  $ that $[\,  P^{\ast \ast} \underset{R} \rightarrowtail R^{[\kappa'']}\,] \overset{(10)} \Rightarrow$
$[\,  P \underset{R} \rightarrowtail R^{[\kappa'']}\,]\,(11)$.
The $R$-module $[\, R^{[\kappa'']}$ is Noetherian$\,]\,(12)$, since $[\, R$     is Noetherian and $\kappa'' <\omega\,]$. It follows from 
$\{ (11),\,  (12)  \}$ that  $[\, \mu_R(P)< \omega \,]$. 
\item Assume that $[\,  \kappa<\omega \,]\,(13)$. Then   $ [\,  P^{\ast} \cong_R Q^{\ast}\,]\Rightarrow [\,P   \cong_R Q\,] $  \\ 
$\mathbf{Proof :}$  $ [\,  P^{\ast} \underset{R}\cong Q^{\ast}\,]\Rightarrow [\, \mu_R(P^{\ast})=\mu_R(Q^{\ast})\,]\overset{(13),\mathrm{(v)}}\Longrightarrow[\, \omega  > \mu_R(Q^{\ast})\,]$
$\overset{\mathrm{(v.b)}}\Longrightarrow[\, \omega > \mu_R(Q)\,]$.
Hence $[\, \kappa < \omega,\, \, \lambda < \omega, \,  \{ P, \, Q  \} \subseteq \mathrm{Proj}(R)\,]\overset{\mathrm{(ii.a)}}\Rightarrow$
$[\,  \{ P, \, Q  \} \subseteq \mathrm{Refl}(R)\,]\,(14)$.
The claim follows from $\{(14), \mathrm{\ref{refllemma}}  \}$. 
\item Assume  $\{ \,[\, \kappa \geq \omega\,]\,(15)$,   $[\, R$ is connected$\,]\,(16)$, $[\, R$ is w-slender$\,]\,(17)\}$. Then   $ [\,  P^{\ast} \underset{R}\cong Q^{\ast}\,]\Rightarrow [\,P   \underset{R}\cong Q\,] $ \\
$\mathbf{Proof :}$ By assumption $ [\,  P^{\ast} \cong_R Q^{\ast}\,]\Rightarrow [\, \mu_R(P^{\ast})=\mu_R(Q^{\ast})\,]\,(18)$.
It follows from $\{\mathrm{(v)}, (15) \}$ that $[\, \mu_R(P^{\ast})\geq \omega \,]\underset{(18)}\Rightarrow$
$[\, \mu_R(Q^{\ast})\geq \omega \,]\underset{\mathrm{(v)}}\Rightarrow [\, \mu_R(Q)\geq \omega \,]$. So  $ [\,  \{ P, \, Q  \} \subseteq \mathrm{Proj}(R),\,   \kappa\geq \omega,\,  \lambda\geq \omega,      \, (16)\,]\overset{\ref{Bass}} \Rightarrow [\,  \{ P, \, Q  \} \subseteq \mathrm{Free}(R)\,]\,(19)$. The claim follows from 
$\{ (17),\, (19),\, \ref{w-slenderremark} \}$ 
\item Assume that $[\,  R$ is  connected, $R$ is w-slender$\,]$. Then   $ [\,  P^{\ast} \cong_R Q^{\ast}\,]\Rightarrow [\,P   \cong_R Q\,] $ \\
$\mathbf{Proof :}$ It follows from $\{\mathrm{vi},\, \mathrm{vii}  \}$
\end{enumerate}

\begin{definition} \label{s-slender} Let $R$ be a Noetherian commutative ring with connected components $R_1,R_2,\dots,R_n$.\\
$\mathrm{(i)}\, R$ is s-slender   $\Leftrightarrow [\, R_k$                is w-slender$\,]\, \forall \, k \in \mathrm{T}_n$. \\
$\mathrm{(ii)}\, R$  is t-slender   $\Leftrightarrow [\,  \exists\,  k \in \mathrm{T}_n\, : (\, R_k$   is slender$\,)\,]$   \\
$\mathrm{(iv)}\, R$  is QS (Quillen-Suslin) iff all f.g projective $R$-modules are free.  \\
$\mathrm{(v)}\, R$  is p-slender ring  $\Leftrightarrow [ R_k$ is slender,  $R_k$ is $\mathrm{QS}$ $\,]\, \forall \, k \in \mathrm{T}_n$ .          
\end{definition}
\begin{theorem}$^{\ast}$ \label{final} Let $P,\, Q$ be projective modules over a Noetherian  commutative ring $R$. Then $:$  \\
$ \mathrm{(i)}\,$ If $R$ is s-slender then    $ [\,  P^{\ast} \underset{R}\cong Q^{\ast}\,]\Rightarrow [\,P   \underset{R}\cong Q\,] $.
$\mathrm{(ii)\,\,}$If $R$ is p-slender then 
$[\,  P \underset{R} \twoheadrightarrow Q\,] \Leftarrow [\,Q^{\ast}   \underset{R} \rightarrowtail P^{\ast}\,]$. \\
$\mathrm{(iii)}$ If $[\, R$ is t-slender, $\{P,Q  \} \subseteq \mathrm{Free}(R)\,]$ then  $\mathrm{(a)}\,[\,  P \underset{R}  \twoheadrightarrow Q\,] \Leftarrow [\,Q^{\ast}   \underset{R} \rightarrowtail P^{\ast}\,]$ \space  $\mathrm{(b)}\, [\,  P^{\ast} \underset{R}\cong Q^{\ast}\,]\Rightarrow [\,P   \underset{R}\cong Q\,] $. \end{theorem}
\begin{proof}
It follows from that there are commutative rings $R_1, R_2,\dots, R_n$ 
such that  $[\, R=\overset{n}{\underset{k=1}\Pi} R_k\,]\,(2)$ and 
$\,(3)\,[ R_k$ is Noetherian and connected$\,]\, \forall \, k \in \mathrm{T}_n$.  Let $P_k=P_{\pi_k},\,Q_k=Q_{\pi_k}$. \\
\underline{$\mathrm{Proof_{}}$ of $\mathrm{(i)}$} : Let $k\in \mathrm{T}_n$. By assumption $[\, R$ is s-slender$\,]\overset{\ref{s-slender}}\Rightarrow [\, R_k$                is w-slender$\,]$ (4).\\  
 Assume that   $ [\,  P^{\ast}  \cong_R Q^{\ast}\,]\overset{\ref{strangesum}}\Longrightarrow$     $[\, (P_k)^{\ast}  \cong_{R_k} (Q_k)^{\ast}]$
$\underset{(3),(4)}{\overset{\ref{projective}\mathrm{.vii}}\Longrightarrow}$
 $[\, P_k  \cong_{R_k} Q_k]$. So $[\, P_k  \cong_{R_k} Q_k]\,  \forall \, k \in \mathrm{T}_n$
 $\underset{\ref{strangesum}} \Rightarrow  [\,  P  \cong_R Q\,]$. Therefore  $ [\,  P^{\ast} \cong_R Q^{\ast}\,]\Rightarrow [\,P   \cong_R Q\,] $.  \\
\underline{$\mathrm{Proof}$ of $\mathrm{(ii)}$} : 
 Let $k\in \mathrm{T}_n$. By assumption $[\, R$ is p-slender$\,]\overset{\ref{s-slender}}\Rightarrow \{ \, [\, R_k$                is slender$\,]$ (4) and   $[\, R_k$ is QS$\,]\,(5)\, \}$.\\ By assumption $\{ P, \, Q  \} \subseteq \mathrm{Proj}(R)\overset{\ref{strangesum}}\Rightarrow \{ P_k, \, Q_k  \} \subseteq \mathrm{Proj}(R_k) \overset{\ref{Bass}}{ \underset{(3),(5)}\Longrightarrow}[\, \{ P_k, \, Q_k  \} \subseteq \mathrm{Free}(R_k)\,]\,(7)$. \\
 Assume that   $ [\,  Q^{\ast} \underset{R} \rightarrowtail P^{\ast}\,]\overset{\ref{strangesum}}\Longrightarrow$     $[\, (Q_k)^{\ast} \underset{R_k} \rightarrowtail (P_k)^{\ast}]$
$\underset{(4),(7)}{\overset{\ref{monoslender}}\Longrightarrow}$
 $[\, Q_k \underset{R_k} \rightarrowtail P_k]$. So $[\, Q_k \underset{R_k} \rightarrowtail P_k]\,  \forall \, k \in \mathrm{T}_n$
 $\overset{\ref{strangesum}} \Rightarrow  [\,  Q \underset{R} \rightarrowtail P\,]$. Therefore  $ [\,  Q^{\ast} \underset{R} \rightarrowtail P^{\ast}\,]\Rightarrow [\,Q   \underset{R} \rightarrowtail P\,] $. Note also that  $[\,  P \underset{R} \twoheadrightarrow Q\,] \Rightarrow [\,Q^{\ast}   \underset{R} \rightarrowtail P^{\ast}\,]$ holds in general from  \ref{ontodual} \\
 \underline{$\mathrm{Proof}$ of $\mathrm{(iii.a)}$} :
 It follows from the assumption  $\{P,Q  \} \subseteq \mathrm{Free}(R)$ that $[\, Q\cong R^{[\mu]},\,  P\cong R^{[\lambda]} \,]\,(8)$, for some cardinals $\mu,\lambda$. Also it follows from the assumption that $R$ is t-slender that $\exists\, k\in \mathrm{T}_n$ such that  $[\,R_k$ is w-slender$\,]\,(9)$. We consider the r.h $\pi_k:R\rightarrow R_k$. Then  $[\, \pi_k$ is f.p$\,]\,(10)$ and  $[\, \pi_k$ is flat$\,]\,(11)$. So 
 $ [\,  Q^{\ast} \underset{R} \rightarrowtail P^{\ast}\,]\Rightarrow [\,R^{\mu} \underset{R} \rightarrowtail R^{\lambda}\,] \underset{(10),(11)}{\overset{\ref{rh2}\mathrm{.vi.b}} \Longrightarrow}$ 
$[\,(R_k)^{\mu} \underset{R}\rightarrowtail (R_k)^{\lambda}\,] \underset{(9)}{\overset{\ref{w-slender4}\mathrm{.ii}}\Longrightarrow}[\,\mu \leq \lambda \,]\Rightarrow [\,R^{[\lambda]} \underset{R}\cong R^{[\mu]}\,]\overset{(8)}\Rightarrow   [\,  P \underset{R} \twoheadrightarrow Q\,]$. 
Hence, we proved that $[\,  P \underset{R}  \twoheadrightarrow Q\,] \Leftarrow [\,Q^{\ast}   \underset{R} \rightarrowtail P^{\ast}\,]$. \\
\underline{$\mathrm{Proof}$ of $\mathrm{(iii.b)}$} :
The ring $R$ is trivially w-slender since $\pi_k$'s are f.p. So the claim follows from $\ref{w-slender4}\mathrm{.ii}$. 
\end{proof} 
\begin{remark} $:$ \\ 
$\mathrm{(i)}$ Note that the projective modules $P,\,Q$ of \ref{final} may not be reflexive (It would make the proof easier). \\
$\mathrm{(ii)}\, $ Theorem $\ref{final}\mathrm{.i}$ is not valid under the weaker assumption that $R$ is w-slender.  Take, for example, 
 $R=S\times F$, where $S$ is a slender commutative ring and $F$ be a field such that $[\, \exists \,  \lambda \in \mathrm{Card}_{\geq \mathfrak{c}} :|F|=2^{\, \lambda}\,]$. Note that the  \emph{Löwenheim–Skolem  Theorem}  guarantees the existence of such a field $F$. It follows from Theorem \ref{failure} that there  are $F$-modules $V,\, W$  such that $[\, V^{\ast}\cong_F W^{\ast}\,]\,(3)$ and 
 $[\, V \ncong_F W\,]\,(4)$. Consider the $R$-modules  
 $P=S \boxplus V, \, Q=S \boxplus W $. Then it follows from \ref{strangesum} that  
  $[\, P^{\ast} \cong_R (R\, \boxplus  V^{\ast})\,]\,(5) $ and 
   $[\, Q^{\ast}\cong_ R(R\boxplus W^{\ast})\,]\,(6) $.  It follows from $\{ (5),(6), \ref{strangesum} \}$ that $[\, P^{\ast}\cong_R Q^{\ast}\,] $, and  that   $[\, P \ncong_R Q \,] $. So  $[\, \{ P, \, Q  \} \subseteq 
   \mathrm{Proj}(R),\,  P^{\ast}\cong_R Q^{\ast}, \,   P \ncong_R Q \,]$. \\
$\mathrm{(iii)}\,$ It follows from  \ref{slenderproduct}  that the above ring  $R=S\times F$ is not slender since $F$ is not slender but $R$  is t-slender since $S$  is slender.   \\  
$\mathrm{(iv)}\,$ We do not know whether Theorem \ref{monoslender} is valid under the weaker assumption that $R$  is w-slender. \\
$\mathrm{(iv)}$ Let $F$ be a field, $n\in \mathbb{N}_{\geq 1}$  and
$R=F[x_1,x_2,\dots,x_n]$. Then $R$ is trivially connected and it  follows from $\ref{slendercollection}\mathrm{.vi}$ that $R$ is slender and from the \emph{Quillen-Suslin Theorem} that $R$ is a QS-ring. Therefore $R$ is both s-slender, t-slender and p-slender.\\
$\mathrm{(v)}\,$  Let $R$ be a Hilbert commutatibve ring. Then
$[\, R$ is s-slender$\,]\Leftrightarrow [\,$All connected components of $R$ are non-Artinian$\,] $.   Let $R=\dfrac{\mathbb{C}[x_1,x_2,\dots ,x_n]}{J} $ be  affine algebra and $\mathrm{Var}(J)\overset{\text{def}}=\{x\in \mathbb{C}^n\, |\,[\, f(x)=0\,]\, \forall \, f\in J \}$. Trivially   
$[\, R$ is w-slender$\,]$ $\Leftrightarrow[\, R $ is non-Artinian$\,]\Leftrightarrow [\, \mathrm{Var}(J)$ is an infinite set$\,]$.\\ Similarly 
$[\, R$ is s-slender$\,]\Leftrightarrow [\,$All components of $R$ are non-Artinian$\,]\Leftrightarrow [\, \mathrm{Var}(J)$ has no isolated points$\,]$ \\ 
\end{remark}

\section{On direct summands of duals of free modules}
\begin{definition}  \label{def.kdim}      Let $R$ be a commutative ring, $M\in R$-$\mathrm{Mod}$ and $\kappa,\, \lambda$ be cardinals.
 Then $:$ \\   $\mathrm{(i)}$ Let $f\in \mathrm{Hom}_R(M,N)$.                     Then $f$  is $R$-dense (or a weak $R$-epimorphism) iff $(\mathrm{Cok}(f))^{\ast}=\bold{0}$. \\
 $\mathrm{(ii)}$  $(\, M  \overset{_w}{\underset{R} \twoheadrightarrow}  N\,) \Leftrightarrow$                            $(\, \exists\,   f\in \mathrm{Hom}_R(M,N) :f $ is $R$-dense$\,)$.\\
 $\mathrm{(iii)} \,$   $\mathrm{E}_R(M)=\{\mu \in \mathrm{Card}\, \big|\, R^{[\mu]}  \overset{_w}{\underset{R} \twoheadrightarrow} M  \} $.
\quad \quad\quad \quad 
$\mathrm{(iv)} \,\, \mathrm{E.dim}_R(M)=\mathrm{min}
(\mathrm{E}_R(M))$                         \\
 $\mathrm{(v)} $     $(\, M \overset{_s} {\underset{R} \rightarrowtail}  N\,) \Leftrightarrow$                            $(\, \exists\,   f\in \mathrm{Hom}_R(M,N) :f $ is $R$-left invertible$\,)$.\\          
 $\mathrm{(vi)} $  $\mathrm{K}_R(M)=\{\mu \in \mathrm{Card}\, \big|\,$  $  R^{[\mu]} \overset{_s} {\underset{R} \rightarrowtail }  M\}$.\quad \quad\quad \quad   $\mathrm{(vii)} \,\, \mathrm{K.dim}_R(M)=\mathrm{sup}
(\mathrm{K}_R(M))$  \quad \quad\quad \quad\quad \quad 
\end{definition}

\begin{lemma} \label{good} Let $R$ be a commutative ring, $\{\kappa,\lambda\}\subseteq \mathrm{Card}$, 
$f\in \mathrm{Hom}_R (R^{[\kappa]},R^{\lambda})$ and $L=\mathrm{Im}(f)$.  Then $: $ \\  If $[\,\omega \leq \kappa <\lambda\,]$ then $:$
$\mathrm{(i)}$ If $L\subseteq R^{[\lambda]}$   $[\, (R^{[\lambda]}/L)^{\ast} \neq \bold{0},\, (R^{\lambda}/L)^{\ast} \neq \bold{0}\,]$
$\mathrm{(ii)\,} (R^{\lambda}/L)^{\ast} \neq \bold{0}$, if $R$ is a domain. 
\begin{proof}
We have assumed that  $[\, \kappa< \lambda\,]\,(1)$  and that $[\, \kappa\geq \omega\,]\,(2)$. \\   \underline{$\mathrm{Proof}$ of $\mathrm{(i)}$} :
 By assumption we have      $[\, \mathrm{Im}(f)= L\subseteq R^{[\lambda]}\,]\Rightarrow$  $ [\,(\, |\, \mathrm{supp}(f(e_i))|\, <\omega\,)\, \forall \, i\in \kappa\,]\,(3) $. Let
  $[\, A = \underset{i\in \kappa}  \cup \mathrm{supp}(f(e_i))\subseteq \lambda \,]\,(3)$. Then $[\, |A|\leq \underset{i\in \kappa}  \Sigma |\mathrm{supp}(f(e_i))| \underset{(3)}\leq  \underset{i\in \kappa}  \Sigma \omega \leq \kappa \cdot \omega\underset{(2)}=\kappa\,]\Rightarrow$
$[\, |A|\leq \kappa\,]\overset{(1)}\Rightarrow[\, |A|<\lambda\,]\Rightarrow[\, (\lambda \smallsetminus A)\neq \emptyset \,]\Rightarrow$
$[\, \exists \, j\, \in \lambda : (\, j
\notin  A=\underset{i\in \kappa}  \cup \mathrm{supp}(f(e_i))\,]\Rightarrow$
\underline{$[\, \exists\,  j \in \lambda : (R\cdot e_j) \cap\,  \mathrm{Im}(f)=\bold 0\,]\,(4)$}. Let $\hat{h}\in \mathrm{Hom}_R(R^{\lambda},R)$ defined by $[\, \hat{h}((x_i)_{i\in \kappa})=x_j\,]$ and $h\in \mathrm{Hom}_R(R^{[\lambda]},R)$ defined by $[\, h((x_i)_{i\in \kappa})=x_j\,]$. Trivially $[\,h \neq 0\,]\,(6)$  and $[\,\hat{h} \neq 0\,]\,(7)$. Also it follows from $(4)$ that  $[\, (\, h(t)=0\, )\, \forall \, t  \in \mathrm{Im}(f)\,]\,(8)$ and $[\, (\, \hat{h}(t)=0\, )\, \forall \, t  \in \mathrm{Im}(f)\,]\,(9)$.
It follows from $\{ (6),\, (8) \}$  that $[\, (R^{[\lambda]}/L)^{\ast} \neq \bold{0}\,]$  and from   $\{ (7),\, (9) \}$
that  $[\, (R^{\lambda}/L)^{\ast} \neq \bold{0}\,] $. \\
 \underline{$\mathrm{Proof}$ of $\mathrm{(ii)}$} : If 
 $[\, \exists \, j \in \lambda : (R\cdot e_j) \cap\,  \mathrm{Im}(f)=\bold 0\,]$ then the claim follows from the proof of the above. Otherwise \quad  $[\,  (R\cdot e_j) \cap\,  \mathrm{Im}(f)\neq  \bold 0\,]\,   \forall \, j \in \lambda $. So $[\, (\forall \, j \in \lambda)\, (\exists \, r_j \in (R\smallsetminus \bold{0}))(\exists\,  a_j\in R^{[\kappa]}):(\, r_j \cdot e_j=f(a_j)\,)]\,(10)$. The map  
$h:R^{\lambda}\rightarrow R^{[]\kappa}$, defined by the formula       
$[\, h((y_j)_{j\in \lambda)}= \underset{j\in \lambda}\Sigma y_j \cdot a_j  \,]\,(11)$ is $R$-linear and well defined by $(10)$. Let  $[\, 
(y_j)_{j\in \lambda}=y \in \mathrm{Ker}(h)\,]\Rightarrow[\, 0=h(y)=\underset{j\in \lambda}\Sigma y_j \cdot a_j\,]\Rightarrow$
$[\, 0=\underset{j\in \lambda}\Sigma y_j \cdot f( a_j)\,]\overset{(10)}=\underset{j\in \lambda}\Sigma y_j  r_j \cdot e_j=0\,]\Rightarrow$
$[\, (y_j r_j=0)\, \forall \, j\in \lambda,$ since $\{e_j\, |\, j\in \lambda\}\in \mathrm{Ind}_R(R^{[\lambda]})\,]\Rightarrow$
$[\, 0=y=(y_j)_{j\in \lambda}$, since $R$ is a domain$\,]\Rightarrow$
 $[\, \mathrm{Ker}(h)=\bold{0}\,]\Rightarrow$
$[\, h$ is 1-1$\,]\Rightarrow$
$[\, R^{[\lambda]} \underset{R} \rightarrowtail R^{[\kappa]} \,]\underset{\ref{innocuous}} \Longrightarrow [\, \lambda \leq \kappa \,] $, which contradicts $(1)$.  \end{proof}
\end{lemma}
\begin{remark} \label{kdim} Let $R$ be a commutative ring and $M,\,N$ be $R$-modules and $\kappa,\lambda$ be  cardinals. Trivially  $:$ \\
$\mathrm{(i)} \,[\, \kappa\in \mathrm{K}_R(M)\,]\Leftrightarrow [\, \exists \,  L\in R$-$\mathrm{Mod} : M  \cong_R ( R^{[\kappa]}\oplus L )\,]$. If $\kappa \geq \omega$ is such  that $\kappa \in \mathrm{K}_R(M))$ then \\
$[\, R^{[\kappa]}\underset{R} \rightarrowtail M\,]\,(1)$. Hence 
$[\, \mathrm{max}\{|R|,\kappa\} \underset{\ref{cardinallemma1}}=|R^{[\kappa]}| \underset{(1)}\leq |M|\,]\Rightarrow [\, \kappa \leq |M|\,]$. So $[\kappa \leq \mathrm{max} \{ \omega,\,|M|\}\,]\,(2)$. Also 
$[\,  0\in \mathrm{K}_R(M)\neq \emptyset\,]\,(3)$. It follows from 
$\{(2),\, (3)  \}$ that $ \mathrm{K.dim}_R(M)=\mathrm{sup}(\mathrm{K}_R(M))$ is well defined.\\
 $\mathrm{(ii)}$ $[\, M \cong_R N\,]\Rightarrow  [\, \mathrm{K.dim}_R(M)=\mathrm{K.dim}_R(N)\,]$ \quad
$\mathrm{(iii)} \,\mathrm{(a)}\,$   $  \mathrm{K}_R(R^{[\kappa]})=\mathrm{Card}_{\leq \kappa}$ \space
$\mathrm{(b)}\,$   $  \mathrm{K.dim}(R^{[\kappa]})=\kappa$\\
$\mathrm{(iv)}\, \mathrm{(a)}$ The class   $  \mathrm{E}_R(M)$ is upward closed.  \quad\quad   $\mathrm{(b)}$ The class is  $  \mathrm{K}_R(M)$ downward closed.                 \\
$\mathrm{(v)}$ Any onto homomorphism is dense and complosition of dense homomorphisms remain dense. \end{remark}

 \begin{remark} \label{edim}
 Let $R$ be a commutative ring and $M,\,N$ be $R$-modules and $\kappa,\lambda$ be  cardinals. Then  $:$ \\ 
$\mathrm{(i)}\, \, \mu_R(M)\in \mathrm{E}_R(M)\neq \emptyset$ hence
$[\, \mathrm{E.dim}_R(M) \geq  \mu_R(M)\,]$.        .
$\mathrm{(ii)}$ $[ M\cong_R N]\Rightarrow [\mathrm{E.dim}_R(M)=\mathrm{E.dim}_R(N)]$. \\
$\mathrm{(iii)}$ If $\lambda \geq \omega$ then $:$ $\mathrm{(a)}\,[\, R$ is $\lambda$-slender$\,]\underset{\ref{def.kdim}} {\overset{\ref{def.half.slender}}\Longrightarrow}$  
$[\, \lambda\in \mathrm{E}_R(R^{\lambda})\,]$. \space
  $\mathrm{(b)}\,\,[\, \kappa < \lambda,\ R\,\, \text{is}\,\, \text{domain}\,]\underset{\ref{def.kdim}} {\overset{\ref{good}\mathrm{.ii}}\Longrightarrow}$  
$[\, \kappa\notin \mathrm{E}_R(R^{\lambda})\,]$.\\
$\mathrm{(c)}\,$  If $R$ is a globally half-slender domain and $\lambda$ is not $\omega$-measurable then  $\mathrm{E.dim}_R(R^{\lambda})=\lambda $   \\\
$\mathrm{(iv)}\,\mathrm{(a)}$ $\kappa \in \mathrm{E}_R(R^{[\kappa]})$  \quad\quad\quad\quad\quad\quad  $\mathrm{(b)}\,\,[\, \kappa \geq \omega,] \underset{\ref{def.kdim}} {\overset{\ref{good}\mathrm{.ii}}\Longrightarrow}$  
$[\, \mathrm{E}_R(R^{[\kappa]})=\mathrm{Card}_{\,\geq \kappa}\,]\Rightarrow$   $[\, \mathrm{E.dim}_R(R^{[\kappa]})=\kappa\,]$.
\end{remark}
\begin{lemma} \label{edim3}Let $R$ be a commutative ring and $\kappa\in \mathrm{Card}$. Then  $ \mathrm{E.dim}_R(R^{[\kappa]})=\kappa$.\end{lemma} 
  \begin{proof}   \underline{Let $\kappa\geq \omega_{}$} : The claim follows from $\ref{edim}\mathrm{.iv.b}$  $\cdot$ \quad \underline{Let $\kappa < \omega$} : If $\lambda\in \mathrm{E}_R(R^{[\kappa]})$  then $[\, \mathrm{Cok}(f)^{\ast}=\bold{0}\,]\,(1)$, for some $[\, f\in \mathrm{Hom}_R(R^{[\lambda]},R^{[\kappa]})\,]\,(2) $. It follows from $(2)$ that $f=\mathfrak{T}^R_A$, for some $A\in R^{\kappa \times \lambda}$. $[\,\bold{0}= \mathrm{Cok}(f)^{\ast} \cong _R \mathrm{Ker}\mathfrak{T}^R_{A^{\top}}\,]\Rightarrow$
$[\, \mathfrak{T}^R_{A^{\top}}$ is $1$-$1\,]\Longrightarrow$
$[\,\lambda \geq \kappa\,]\overset {\ref{kdim}\mathrm{.iv.a}}\Longrightarrow$ 
$[\, \mathrm{E.dim}_R(R^{[\kappa]})=\kappa\,]$. 
\end{proof}
\begin{lemma} \label{weak}Let $R$ be a globally half-slender integral domain.  
 Then $\mathsf{L}_R^w$ holds. \end{lemma}
 \begin{proof}Let $\lambda,\, \mu$ be non $\omega$-measurable cardinals. It follows from $\{\ref{edim}\mathrm{.iii.c} ,\, \ref{edim3}  \}$ that $[\, \mathrm{E.dim}_R(R^{\kappa})=\kappa\,]\,(3)$ and  $[\, \mathrm{E.dim}_R(R^{\lambda})=\lambda\,]\,(4)$. So
 $[\,R^{ \kappa}\cong_R R^{ \lambda}\,]\overset{\ref{edim}\mathrm{.ii}} \Longrightarrow $ 
$[\, \mathrm{E.dim}_R(R^{\kappa})= \mathrm{E.dim}_R(R^{\lambda})\,]$
$\overset{(3),(4)}\Longrightarrow$ 
  $[\, \kappa=\lambda) \,]$. Therefore  $\mathsf{L}_R^w$
   holds. \end{proof}

\begin{lemma}  Let $R_1,\, R_2\dots ,R_n$ be commutative rings, $R=\overset{n}{\underset{i=1}\Pi} R_i$ and $[\, M_i=M_{\pi_i}\,]\, \forall \, i\in \mathrm{T}_n$.  Let also $\kappa,\, \lambda$ be cardinals. Then $:$  \\
$\mathrm{(i)}\,\,\mathrm{(a)}$ $[\, \kappa\in \mathrm{K}_R(M)\,]\Leftrightarrow$
$[\, \kappa\in \mathrm{K}_{R_i}(M_i)\,]\, \forall \, i \in \mathrm{T}_n$. \quad $\mathrm{(b)}$  $\mathrm{K.dim}_R(M)=\mathrm{max} \{ \mathrm{K.dim}_{R_i}(M_i) \,| \, i\in \mathrm{T}_n \}$.  \\
$\mathrm{(ii)\,\,\mathrm{(a)}}$ $[\, \kappa\in \mathrm{E}_R(M)\,]\Leftrightarrow$
$[\, \kappa\in \mathrm{E}_{R_i}(M_i)\,]\, \forall \, i \in \mathrm{T}_n$. \quad $\mathrm{(b)}$  $\mathrm{E.dim}_R(M)=\mathrm{min} \{ \mathrm{E.dim}_{R_i}(M_i) \,| \, i\in \mathrm{T}_n \}$.  \\
$\mathrm{(iv)}$  Let $[\, (R_i,\mathfrak{m}_i)$ be local Artinian,  
$\mu _i=| \dfrac{R_i}{\mathfrak{m}_i}|^{\, \lambda}$
$\,]\, \forall \, i
\in \mathrm{T}_n$ and $\Lambda=\{\mu_i \,|\, i\in \mathrm{T}_n  \}$.  If $[\, \lambda \geq \omega\,]\,(3)$ then $:$ \\ $\mathrm{(a)}\, [\, \mathrm{K.dim}_{R_i}(R_i^{\lambda})=\mathrm{E.dim}_{R_i}(R_i^{\lambda})=\mu_i\,]$. 
 $\mathrm{(b)}\, [\, \mathrm{K.dim}_R(R^{\lambda})=\mathrm{max}(\Lambda)\,]$.  $\mathrm{(c)}\,[\,  \mathrm{E.dim}_R(R^{\lambda})=\mathrm{min}(\Lambda)\,]$.  \end{lemma}
\begin{proof} \underline{$\mathrm{Proof}$ of $\mathrm{(ii)}$} : 
$[\, \kappa\in \mathrm{E}_R(M)\,]\overset{\ref{kdim}\mathrm{.iii}}\Longleftrightarrow[\,    R^{[\kappa]} \overset{_w} {\underset{R} \twoheadrightarrow }  M\,]$  $\overset{\ref{strangesum}\mathrm{.ix.e}}\Longleftrightarrow[\,   (R^{[\kappa]})_{\pi_i} \overset{_w} {\underset{R_i} \twoheadrightarrow }  M_{\pi_i}\,]\,  \forall \, i \in \mathrm{T}_n $
 $\overset{\ref{rh}\mathrm{.i}}\Longleftrightarrow[\,   R_i^{[\kappa]} \overset{_w}{\underset{R_i} \twoheadrightarrow}  M_i\,]\,  \forall \, i \in \mathrm{T}_n $
$\overset{\ref{kdim}\mathrm{.iii}}\Longleftrightarrow$
$[\, \kappa\in \mathrm{E}_{R_i}(M_i)\,]\, \forall \, i \in \mathrm{T}_n$. So $\mathrm{(ii.a)}$ was proven. Trivially $\mathrm{(ii.b)}$
follows from  $\{\mathrm{(ii.a)},\, \ref{def.kdim}\mathrm{.iv}  \}$
\quad \quad\quad \quad  \underline{$\mathrm{Proof}$ of $\mathrm{(i)}$} : Similar to the proof of $\mathrm{(ii)}$. We apply
 $\{\ref{kdim}\mathrm{.vi},\, \ref{strangesum}\mathrm{.ix.f}  \} $    instead. \\
\underline{$\mathrm{Proof}$ of $\mathrm{(iii)}$}: $\mathrm{iii.a}$
follows from  $ \ref{EKartinian}$, $\mathrm{iii.b}$ follows from 
$\{ \mathrm{iii.a},\, \mathrm{i.b} \}$ and  $\mathrm{iii.c}$ follows from $\{ \mathrm{iii.a},\, \mathrm{ii.b} \}$.   \end{proof}

\begin{lemma} \label{cohn1} Let $R$ be a commutative ring that  is not a field and $\kappa $ be an infinite cardinal.  Then $:$\\
$\mathrm{(i)} \, [\, R$ is a $\mathsf{UFD}\,]\Rightarrow[\, R^{\omega}\notin \mathrm{Proj}(R)\,]$ \quad\quad\quad\quad\quad $\mathrm{(ii)} \, [\, R$ is a $\mathsf{DVR}\,]\Rightarrow[\, R^{\omega}\notin \mathrm{Proj}(R)\,]$ \\
$\mathrm{(iii)} \, [\, R$ is  $w$-slender$\,]\Rightarrow[\, R^{\omega}\notin \mathrm{Free}(R)\,]$ \quad\quad\quad $\mathrm{(iv)} \, [\, R$ is $\kappa$-half-slender$\,]\Rightarrow[\, R^{\,\kappa}\notin 
\mathrm{Free}(R)\,]$  \end{lemma}
\begin{proof}
$\ref{cohn1}\mathrm{.i}$ follows from Theorem 3 of \cite{cohn}, 
$\ref{cohn1}\mathrm{.ii}$ follows from $\ref{cohn1}\mathrm{.i}$ and 
$\ref{cohn1}\mathrm{.iii}$ was proven in  \ref{slender12}. \\
\underline{$\mathrm{Proof}$ of $\mathrm{(iv)}$} : By assumption $[\kappa \geq \omega\,]\,(1)$ and  $ [\, R$ is  $\kappa$-half-slender$\,]\,(2)$ . Assume that 
$ [\, R^{\,\kappa}\in \mathrm{Free}(R)\,]\Rightarrow $ 
$[\, R^{\,\kappa }  \cong_R R^{[\lambda]}$, for  $\lambda \in \mathrm{Card}\,]\,(3)\overset{\ref{lost}} \Longrightarrow[\, \lambda \geq 2^{\, \kappa}> \kappa\,]\Rightarrow[\, \lambda > \kappa\,]\,(4)$. Also $\{(1),\, (2),\, \ref{edim}\mathrm{.iii.a}   \} \Rightarrow$     $[\, \kappa\in \mathrm{E}_R(R^{\kappa})\,]\Rightarrow
[\, \kappa\in \mathrm{E}_R(R^{\kappa})\underset{(3)}= \mathrm{E}_R(R^{[\lambda]})\underset{\ref{edim}\mathrm{.iv}}=\mathrm{Card}_{\geq \lambda} \,]\Rightarrow[\, \kappa \geq \lambda\,] $, which contradicts $(4)$. Therefore $ [\, R^{\,\kappa}\notin \mathrm{Free}(R)\,]$.     
 \end{proof}
\begin{lemma}\label{completelocalrings} Let $S$ be a Noetherian domain such that $\mathrm{dim}(S)=1$. Then $:$\\
$\mathrm{(i)}\, S$ is non-slender $\Leftrightarrow S$ is a compete local ring. \\$\mathrm{(ii)}$  If $S$ is a compete local ring then
there is a r.h \space $\psi:T\rightarrow S$ such that $T$ is a complete $\mathsf{DVR}$ and $\psi$ is f.p and 1-1.\end{lemma}
\begin{proof}  See Theorem 5.4 of \cite{jensen} for 
$\ref{completelocalrings}\mathrm{.i}$. Also $\ref{completelocalrings}\mathrm{.ii}$ follows from $\ref{completelocalrings}\mathrm{.i}$ and 
\emph{Cohen's structure theorem} (see Lemma 160.11 of the \emph{Stacks Project in Commutative Algebra}).\end{proof}

\begin{lemma} \label{nicelemma} Let $R$ be Noetherian  commutative ring. If $R$ is non-Arinian then 
$ R^{\omega} \notin \mathrm{Proj}(R)$.
\end{lemma} 
\begin{proof} Assume that $ [\, R^{\omega} \in \mathrm{Proj}(R)\,]\,(1)$. 
The ring $R$ is non-Artinian hence $\exists\, \mathfrak{p}\in \mathrm{Spec}(R)$ such that $[\, S=(R / \mathfrak{p})$ is a Noetherian domain$\,]\,(2)$ and $[\, \mathrm{dim}(S)=1\,]\,(3)$. Consider the r.h 
$\phi=\pi^R_{\mathfrak{p}}:R\rightarrow \ (R / \mathfrak{p})$. Then 
$[\, \phi$ is f.p$\,]\,(4)$ since  the ring $R$ is Noetherian.
It follows from $\{(1),\,(4),\,  \ref{remainprojective}\}$ that 
$ [\, S^{\omega} \in \mathrm{Proj}(S)\,]\,(5)$ and from
 $\{(5),\, (2),\, \ref{noncommutativeIBN},\, \ref{Bass}       \}  $
 that  $ [\, S^{\omega} \in \mathrm{Free}(S)\,]\,(6)$. Also it follows from $\{(6),\, \ref{slender12}       \}  $ that $[\, S$ is 
 non-slender$\,]\,(7)$ and from  $\{ (2),(3),(7),\, \ref{completelocalrings}       \}  $ that  there is a r.h \space $\psi:T\rightarrow S$ such that $[\, T$ is a complete
  $\mathsf{DVR}\,]\,(8)$ and $[\, \psi$ is f.p and 1-1$\,]\,(9)$. It follows from $(8)$ that $[\, T$ is a $\mathsf{PID}\,]\,(9)$  The ring $S$ becomes, under $\psi$, a torsion free (since $S$ is a domain and $\psi$ is 1-1) and f.g $T$-module.  $S$ is a f.g  free $T$-module.  Hence, it follows from $(9)$ that   $[\, \exists\, \kappa\in \mathrm{Card}_{\, < \omega}\,] (10)$ such that $[\, S\underset{R}\cong R^{[\kappa]}\,]\,(11)$. It follows from $(6)$ that $[\, S^{\omega} \cong_S S^{[\lambda]}\,]\Rightarrow [\, S^{\omega} \cong_T S^{[\lambda]}\,]\Rightarrow  [\, (T^{\kappa})^{\omega} \cong_T (T^{[\kappa]})^{[\lambda]}\,]\Rightarrow$
$ [\, T^{\kappa\cdot \omega }\cong_T T^{[\kappa \cdot \lambda]}\,]$   
$\overset{(10)}\Rightarrow   [\, T^{ \omega }\cong_T T^{[\kappa \cdot \lambda]}\,]     \Rightarrow[\, T^{\omega}\in \mathrm{Free}(T)\,]\Rightarrow [\, T^{\omega}\in \mathrm{Proj}(T)\,]$, which contradicts 
$\{ \ref{cohn1}\mathrm{.ii},\, (8) \}$. Therefore $ [\, R^{\omega} \notin \mathrm{Proj}(R)\,]$. \end{proof}
We now provide a simpler proof of Theorem \ref{relief}, using $\ref{cohn1}\mathrm{.ii}$  instead of the ubiquitous $\ref{oneil}\mathrm{.i}$. 
\begin{theorem} \label{oneil3} Let $R$ be Noetherian  commutative ring and $\lambda\in \mathrm{Card}$. Assume that  $\lambda \geq \omega$. Then \\
$[\, R^{\lambda} \in \mathrm{Proj}(R)\,]\,(1)\Leftrightarrow [\, R$ is Artinian$\,]\,(2)$. 
\end{theorem}

\begin{proof} The direction $[\, (2)\Rightarrow (1)\,]$ follows from  \ref{collection}.\quad\quad\ \underline{Proof of  $[\, (1)\Rightarrow (2)\,]$} : It follows from the  assumption 
$[\, \lambda \geq \omega\,]$ that $[\, \lambda =\lambda +\omega\,]\overset{\ref{trivial1}}\Rightarrow$
$[\, R^{\lambda}\cong R^{\lambda} \oplus R^{\omega}\,]$
 $\underset{(1)}   {\overset{\ref{reflexivesum}\mathrm{.i}} \Longrightarrow}[\, R^{\omega} \in \mathrm{Proj}(R)\,]$
 $\overset{\ref{nicelemma}}\Rightarrow [\, R$ is Arinian$\,]$.  
\end{proof}
\begin{lemma}\label{dvr-reduction} If $\mathsf{L}_T$ holds for for any  complete $\mathsf{DVR}$ $T$ then  $\mathsf{L}_R$ holds for any non-Artinian commutative ring $R$ that is Noetherian. \end{lemma}
\begin{proof} Let $R$ be  non-Artinian commutative ring and $S,\, T,\, \phi,\,  \psi$ be as in the proof of \ref{nicelemma}. By asumption   
$[\, \mathsf{L}_T$ holds$\,]$ since $T$ is a $\mathsf{DVR}$. Let also $\mu, \nu$ be cardinals. Assume that   $[\, R^{\mu} \cong_R R^{\nu}\,]\,(12)\Rightarrow [\, S^{\mu} \cong_S S^{\nu}\,] \Rightarrow[\, S^{\mu} \cong_S S^{\nu}\,]\Rightarrow$
$[\, (T^{\kappa})^{\mu} \cong_T (T^{\kappa})^{\nu}\,]\Rightarrow[\, (T^{\kappa \cdot \mu}) \cong_T (T^{\kappa \cdot \nu})\,]\,(13)$. \\
$\cdot$ Let $[\, \mu < \omega\,]\,(14)$. It follows from 
$\{(12),\, (14),\, \ref{noncommutativeIBN}\mathrm{.iii}  \}$ that 
$[\, \nu <\omega\,]\,(15)$ and from $\{(12),\, (14), (15), \, \ref{rankiswelldefined}  \}$ that $[\, \mu=\nu\,]$. \\
$\cdot$ Let $[\, \mu \geq \omega\,]\,(16)$. It follows from 
$\{(12),\, (16),\, \ref{noncommutativeIBN}\mathrm{.iii}  \}$ that 
$[\, \nu \geq \omega \,] \,(17)$. So 
$(13)\Leftrightarrow[\, (T^{\kappa \cdot \mu}) \cong_T (T^{\kappa \cdot \nu})\,]\overset{(16),(17)}\Longrightarrow$
$[\, T^{\mu}\cong_T  T^{\nu}\,]\overset{\mathsf{L}_T}\Longrightarrow[\, \mu=\nu  \,]$. \end{proof}

\begin{lemma} \label{ses} Let  $[\, \bold{0}\rightarrow M \overset{f} \rightarrow
N \overset{g}  \rightarrow L \rightarrow \bold{0}\,]\,(2)$   be a s.e.s, i.e  $[\, $f is 1-1, g is onto, $\mathrm{Ker}(f)=\mathrm{Im}(g)\,]$.  \\
$\mathrm{(i)}$ $[\, f$ is $R$-left invertible$\,]\Leftrightarrow$
 $[\, g$ is $R$-right  invertible$\,]$. In this case we say  the s.e.s  $(2)$ is $R$-split. \\ $\mathrm{(ii)}$
 $[\,(2)$ is $R$-split$\,]\Rightarrow$  $[\, N\underset{R}\cong (M \oplus L)\,]$ \space  $\mathrm{(iii)}$
 $[\, \mathrm{Ext}_R(L,M)=\bold{0} \,] \Rightarrow[\, f$ is $R$-left invertible$\,]$ .\\
$\mathrm{(iv)}\, [\, L\in \mathrm{Proj}(R)\,]\Rightarrow[\, (2)$ is $R$-split$\,]$.  \quad\quad\quad\quad\quad\quad
$\mathrm{(iv)}\,[\, M\in \mathrm{Inj}(R)\,]\Rightarrow[\, (2)$ is $R$-split$\,]$.
  
\end{lemma}

\begin{lemma} \label{nonexpected1} Let $T$ be a complete $\mathsf{DVR}$, $L \in T$-$\mathrm{Mod}$ and $\kappa$ be an infinite
 cardinal. Then $:$ \\
$\mathrm{(i)}$ If $L$ is  $T$-torsion free then $\mathrm{Ext}_T^1(L,T)=0$.  \quad\quad\quad\quad  $\mathrm{(ii)}\, \Phi^{\kappa}_T$ is $T$-left invertible. \\  $\mathrm{(iii)}\, T^{\kappa} \cong_T (T^{[\kappa]} \oplus T^{\cdot \kappa})$ \quad
$\mathrm{(iv)}\,\mathrm{K.dim}_T( T^{\kappa}) \geq \kappa$. \quad
$\mathrm{(v)} \, T^{\cdot \kappa}$ is  Whitehead  but not projective
$T$-module. \end{lemma}

\begin{proof} See of $\cite{eklof}$ for a proof of $\mathrm{(i)}$.
Consider the s.e.s $[\, \bold{0}\longrightarrow T^{[\kappa]}\overset{\Phi^{\kappa}_T} \longrightarrow
T^{\kappa}  \overset{\Pi^{\kappa}_T}  \longrightarrow T^{\cdot \kappa} \longrightarrow \bold{0}\,]$ and  $L=T^{\cdot \kappa}$. Trivially 
$[\, L$ is  $T$-torsion free$\,]\overset{\mathrm{(i)}}\Longrightarrow$  $[\, \mathrm{Ext}_T^1(L,T)=0\,]\,(2)  \Rightarrow$
$[\, \mathrm{Ext}_T^1(L,T^{[\kappa]}) \cong_T \underset{i\in \kappa }\Pi \mathrm{Ext}_T^1(L,T)  \cong_T \underset{i\in \kappa }\Pi \bold{0}  =\bold{0}\,]\Rightarrow$
$[\, \mathrm{Ext}_T^1(L,T^{[\kappa]})=\bold{0}\,]\overset{\ref{ses}}\Longrightarrow$
 $[\,  \Phi^{\kappa}_T$ is $T$-left invertible$\,]\,\mathrm{(ii)}\overset{\ref{ses}} \Longrightarrow$
 $[\, T^{\kappa} \cong_T (T^{[\kappa]} \oplus T^{\cdot \kappa})\,]\,\mathrm{(iii)}\Rightarrow [\, \kappa \in \mathrm{K}_T(T^{[\kappa]})\,]\Rightarrow $
$[\, \mathrm{K.dim}( T^{\kappa}) \geq \kappa\,]\,(\mathrm{iv})$.
It follows from   $(2)$ that $[\, T^{\cdot \kappa}$ is a Whitehead $T$-module$\,]\,(3)$. Assume that  $[\, T^{\cdot \kappa}\in \mathrm{Proj}(T)\,]\overset{(2)}\Longrightarrow$
$[\, T^{ \kappa}\in \mathrm{Proj}(T)\,]\overset{\ref{oneil3}}\Longrightarrow [\, T$ is Artinian$\,]$, which is a contradiction. 
Therefore $[\, T^{\cdot \kappa} \notin \mathrm{Proj}(T)\,]\,(4)$. Note that 
$\mathrm{(iv)}$ follows from $\{ (3),\, (4) \}$ \end{proof} 
\begin{remark}  \label{accident} Let $R$ be a Noetherian commutative ring that is not Artinian . Then $:$  \\
$\mathrm{i)}$   \underline {Conjecture} : If $T$ is a complete $\mathsf{DVR}$ then
 $:$ $\mathrm{K.dim}( T^{\kappa}) = \kappa$ \\
$\mathrm{(ii)}$  $[\, $Conjecture $\mathrm{\ref{accident}\mathrm{.i}}\,]\overset{\ref{kdim}\mathrm{.ii}}\Longrightarrow$
$[\, \mathsf{L}_T$ holds fon any complete $\mathsf{DVR}$ $T\,]$ 
$\overset{\ref{dvr-reduction}}\Longrightarrow$   $[\, \mathsf{L}_R$ holds$\,]$.\\ 
$\mathrm{iii)}$ It was proven in that $[$\,If $R$ is a  that is non-local domain then $R$ is  half-slender$\,]$. \\
$\mathrm{(iv)}$ \underline{$\text{Conjecture}$}: If $R$ is a Noetherian domain that is non-local then $R$
 is globally half-slender. \\  
$\mathrm{(iv)}$ Let  $T$ is a complete $\mathsf{DVR}$.  If $T$ was globally half-slender and Conjecture $\ref{accident}\mathrm{.iv}$ holded then 
$\mathsf{L}_T^w$ would hold. Unfortunately, it is proven  in $\cite{jensen}$    that  $T$ is not even half-slender. Despite it 
 Conjecture $\ref{accident}\mathrm{.iv}$ would still be useful since it would prove the validity of  $\mathsf{L}_T^w$ without invoking the (hard to prove) slendereness results.  
\end{remark}

\begin{theorem} \label{steinitz8}
Let $R$ be a commutative Artinian ring, $\{ P, \, Q  \} \subseteq \mathrm{Proj}(R)$ and  $f\in \mathrm{Hom}_R(P,Q)$. Then  
$\mathrm{(i)\, }\, [\,f$ is $1$-$1\,] \Leftrightarrow$  $[\, f$ is $R$-left invertible$\,]$. \quad \quad\quad \quad\quad \quad $\mathrm{(ii)\, }\, [\,f$ is onto$\,]  \Leftrightarrow$  $[\, f$ is $R$-right invertible$\,]$.
\end{theorem}
\begin{proof}It follows $\ref{collection}\mathrm{.v}$ from  that there is a factorizarion
 $R=\overset{n}{\underset{k=1}\Pi} R_k $, for some  $[\, ($local and Artinian  commutative rings $R_k)\, \forall \,
  k\in \mathrm{T}_n\,]\,(1)$. Let $k\in \mathrm{T}_n$. Consider the
  r.h $\pi_k: R\rightarrow R_k$ and the  $R_k$-modules $P_k=P_{\pi_k},\, Q_k=Q_{\pi_k}$.  Then     $[\, \{ P, \, Q  \} \subseteq \mathrm{Proj}(R)\,]\overset{\ref{remainprojective}}\Longrightarrow[\, \{ P_k, \, Q_k  \} \subseteq \mathrm{Proj}(R_k)\,]$  
 $\overset{\ref{collection}\mathrm{.i}} {\underset{(1)}\Longrightarrow}$   $[\, \{ P_k, \, Q_k  \} \subseteq \mathrm{Free}(R_k)\,]\,(3)$ \\
\underline{$\mathrm{Proof}$ of $\mathrm{(i)}$} : Let $k\in \mathrm{T}_n $. It follows from $\{ (3),\, \ref{steinitz7}\mathrm{.i} \}$ that 
  $[\, (\, f_k$ is 1-1$\,)]$
$\Leftrightarrow (\, f_k$ is $R_k$-left invertible$)\,]\,(4)$.
Therefore  $[\,f$ is 1-1$\,]\overset{\ref{strangesum}\mathrm{.ii}} \Longleftrightarrow$ 
$[\,f_k$ is 1-1$\,]\,  \forall \, k \in \mathrm{T}_n$ $\overset{(3)} \Leftrightarrow$ 
$[\,f_k$ is $R_k$-left invertible$\,]\,  \forall \, k \in \mathrm{T}_n$   $\overset{\ref{strangesum}\mathrm{.vi}}  \Longleftrightarrow$
$[\,f$ is $R$-left invertible$\,]$.
\\
\underline{$\mathrm{Proof}$ of $\mathrm{(ii)}$} : Similar to the above. 
The only difference is that the relation   $[\, (\, f_k$ is onto$\,)]$
$ \overset{\ref{steinitz7}\mathrm{.ii}}\Leftrightarrow (\, f_k$ is $R_k$-right invertible$)\,]\,(4')$ that is used instead of $(4)$ is  simple since it relies on the trivial Lemma $\ref{steinitz7}\mathrm{.ii} $.\end{proof}

\begin{theorem} \label{left.invertible} $^{\ast}$
Let $R$ be a commutative Artinian ring and $\kappa$ be an infinite cardinal. Then $:$ \\ 
$\mathrm{(i)\, }\,  \Phi^{\kappa}_R $   $R$-left invertible. \quad\quad \quad\quad 
$\mathrm{(ii)}\, R^{\kappa} \cong_ R(R^{[\kappa]} \oplus R^{\cdot \kappa})$  \quad\quad \quad 
$\mathrm{(iii)}\, R^{\cdot \kappa} \in \mathrm{Proj}(R)$  \\ 
$\mathrm{(iv)}\, [\, R^{\cdot \kappa} \in \mathrm{Free}(R)\,]\Leftrightarrow[\, 
 R^{ \kappa} \in \mathrm{Free}(R)\,]$  \space
$\mathrm{(v)}\, [\,|R| \leq   \mathfrak{c}\,]\Leftrightarrow[\, 
 R^{\cdot \kappa} \in \mathrm{Free}(R)\,]$ 
\space $\mathrm{(vi)}\,R $ is not $\kappa$-half-slender.     
\end{theorem}

\begin{proof} The map $ \Phi^{\kappa}_R :R^{[\kappa]}\rightarrow R^{\kappa} $ was defined in $\ref{def.half.slender}$. Trivially
$[\,  \Phi^{\kappa}_R$ is 1-1$\,[\,(1)$
  Also
$[\, \{ R^{[\kappa]},R^{\kappa}  \} \underset{\ref{collection}}\subseteq \mathrm{Proj}(R)\,]\,(2)$. So it follows from 
$\{(1),\, (2),\, \ref{steinitz8}\mathrm{.i}  \}$ that
$[\, \Phi^{\kappa}_R $   $R$-left invertible$\,]\,\mathrm{(i)}$
$\overset{\ref{ses}}\Longrightarrow$
$[\, R^{\kappa} \cong_R (R^{[\kappa]} \oplus R^{\cdot \kappa})\,]\,\mathrm{(ii)}$
$\overset{(2)}\Rightarrow$ 
$[\, R^{\cdot \kappa} \in \mathrm{Proj}(R)\,]\, \mathrm{(iii)}$.\\
\underline{$\mathrm{Proof}$ of $\mathrm{(iv)}$} :
 It follows from $\ref{important}$ that 
there are $[\,$Artinian local commutative rings $(R_i,\mathfrak{m}_i)\,]\,(3)$ such that $[\,  R=\overset{n}{\underset{i=1}\Pi} R_i\,]\,(4)$. Let $i \in \mathrm{T}_n$.  It also follows that   $[\, R_i^{\kappa}\in \mathrm{Free}(R)\,]$ and  $[\,\dim_{R_i} (R_i)^{\kappa}=\mu_i^{\, \kappa}\,]\,(5)$, where  $[\mu_i=|(R_i/ \mathfrak{m}_i)|\,]\,(6)$. It also follows from  $\{ \ref{collection},\, \mathrm{(iii} \}$  that  $[\, (R_i)^{\kappa}\in \mathrm{Free}(R_i)\,]\,(7)$ and from $\{\mathrm{(ii)},\,(5)\}$ that $[\,\kappa+ \dim _{R_i}(R_i)^{\cdot\kappa}=\mu_i^{\kappa}\,]\,(8)$. Also  $[\, \mu_i^{\kappa}\geq 2^{\, k}> \kappa\,]\,(9)$ and by assumption  $[\, \kappa \geq \omega, \,]\,(10)$. It follows from $\{(5),\,  (8),\, (9),\, (10) \}$ that  $[\,\dim_{R_i} (R_i)^{\kappa}=\dim_{R_i} (R_i)^{\cdot \kappa}\,]\,(11)$. It follows from $\ref{strangesum}\mathrm{.iii.a}$ that 
  $[\, \mathrm{Cok}( \Phi^{\kappa}_R)\underset{R}\cong  \overset{n}{\underset{i=1} \boxplus}   \mathrm{Cok}( \Phi^{\kappa}_{R_i})\,]\,(12)$. Therefore $ [\, R^{\cdot \kappa} \in \mathrm{Free}(R)\,]\Leftrightarrow[\,\mathrm{Cok}( \Phi^{\kappa}_R)\in \mathrm{Free}(R)  \,]\overset{(12)}\Leftrightarrow$
$ [\,(\overset{n}{\underset{i=1} \boxplus}   \mathrm{Cok}( \Phi^{\kappa}_{R_i}))\in \mathrm{Free}(R)\,] $
$\Longleftrightarrow[\,(\dim_{R_i} (R_i)^{\cdot \kappa}= \dim_{R_1} (R_1)^{\cdot \kappa})\, \forall \, i \in \mathrm{T}_n\,]\overset{(11)}\Leftrightarrow$ 
$[\,(\dim_{R_i} (R_i)^{ \kappa}= \dim_{R_1} (R_1)^{ \kappa})\, \forall \, i \in \mathrm{T}_n\,]\overset{\ref{important}\mathrm{.ii}}\Longleftrightarrow$
 $ [\, R^{ \kappa} \in \mathrm{Free}(R)\,]$. Hence 
$ [\, R^{\cdot \kappa} \in \mathrm{Free}(R)\,]\Leftrightarrow[\, 
 R^{ \kappa} \in \mathrm{Free}(R)\,]$  \quad\quad\quad\quad \underline{$\mathrm{Proof}$ of $\mathrm{(v)} $} : 
 It follows from $ \{ \ref{artinianfree},\,  \mathrm{(iv)}\}$. 
 \underline{$\mathrm{Proof}$ of $ \mathrm{(vi)}$} : 
 Let $i \in \mathrm{T}_n$. It follows from $(11)$ that
 $(R_i)^{\cdot \kappa}\cong _{R_i} (R_i)^{[\lambda]}$, for some $\lambda \in \mathrm{Card}\,]\Rightarrow$
 $[\, ((R_i)^{\cdot \kappa})^{\ast}\cong _{R_i} (R_i)^{\lambda}\neq \bold{0}\,]\,(14)\Rightarrow[\, R_i$ is not $\kappa$-half-slender$\,]\overset{\ref{slenderproduct} \mathrm{.iii}}\Longrightarrow$
 $[\,  R$ is not $\kappa$-half-slender]. \end{proof} 

\begin{lemma} \label{invertible} Let $R$ be a  commutative ring and $\kappa$ be any cardinal. Then $:$ \\
$\mathrm{(i)}$ If $R$ is $w$-slender then $ \Phi^{\kappa}_R$  is not  $R$-left invertible. \\
$\mathrm{(ii)}$ If $R$ is Hilbert then $[\,  \Phi^{\kappa}_R$  is  $R$-left invertible$\,]\,\mathrm{(a)}\Leftrightarrow[\, R$ is non-Artinian$\,]\,\mathrm{(b)}$. 

\end{lemma} 
\begin{proof} \underline{$\mathrm{Proof}$ of $\mathrm{(i)}$} : There is a r.h $\phi: R\rightarrow S$ such that  $[\, \phi: R\rightarrow S$
is f.p$\,]\,(1)$ and   $[\, S$ is a slender$\,]\,(2)$. It follows from $\{(2),\, \ref{rh}\}$ that $[\, ( \Phi^{\kappa}_R)_{\phi}\underset{S}\sim  \Phi^{\kappa}_S\,]\,(3) $.
Assume that  $[\,  \Phi^{\kappa}_R$  is  $R$-left invertible$\,]\Rightarrow$  $[\, ( \Phi^{\kappa}_R)_{\phi}$  is  $S$-left invertible$\,]\underset{(3)}\Rightarrow$
$[\,  \Phi^{\kappa}_S$  is  $S$-left invertible$\,]\,(4)\Rightarrow$
$[\, S^{\kappa} \cong_ S(S^{[\kappa]} \oplus S^{\cdot \kappa})\,]\Rightarrow$  $[\, (S^{\kappa})^{\ast} \cong_ S((S^{[\kappa]})^{\ast} \oplus (S^{\cdot \kappa})^{\ast})\,] \underset{(2)}\Rightarrow$
$[\, S^{[\kappa]} \cong_ S(S^{\kappa} \oplus (S^{\cdot \kappa})^{\ast})\,]\Rightarrow$
$[\, S^{\kappa}\in \mathrm{Proj}(S)\,]\underset{\ref{oneil3}}\Rightarrow[\, S$ is Artinian$]$, which contradicts
$\{\ref{freeartinian}\mathrm{.iv},\, (2)   \}$. Therefore
$ \Phi^{\kappa}_R$  is not  $R$-left invertible. \\
\underline{$\mathrm{Proof}$ of $\mathrm{(ii)}$} : $\cdot$The direction 
$\mathrm{(a)}\Rightarrow \mathrm{(b)} $ follows from  $\ref{left.invertible}\mathrm{.i}$ and the converse follows from 
$\{ \mathrm{(i)},\, \ref{hilbert}  \mathrm{.xi}\}$.  \end{proof}

\begin{theorem} \label{direct.cofree} $^{\ast}$
Let $R$ be a Noetherian commutative ring such that  $|R|\in  \mathcal{L}_1'$, $M\in R$-$\mathrm{Mod}$ and  $\{\kappa,\, \lambda\} \subseteq \mathrm{Card}$. Then $:$
$\mathrm{(i)}$ Assume that $ R$ is w-slender.
Then $[\, R^{[\lambda]} \overset{_s}{\underset{R}\rightarrowtail} R^\kappa\,]\Rightarrow[\,\lambda < \omega \,]$.  \\
$\mathrm{(ii)}$ If $[\, R$ is slender and connected,  $\kappa \in \mathcal{L}_1',\, \mu_R(M)\geq \omega\,]$ then $[\, M \overset{_s}{\underset{R}\rightarrowtail} R^\kappa\,]\Rightarrow[(\exists \, \mu \in \mathrm{Card}):(M \cong_R R^{\mu}) \,]$. \end{theorem}
\begin{proof} $\mathcal{L}_1'$ is the class of non-$\omega$-measurable cardinals. We have  assumed that  
$[\,|R|\in \mathcal{L}_1'\,]\,(1)$. \\
\underline{$\mathrm{Proof}$ of $\mathrm{(i)}$}: By assumption
 $[ R$ is w-slender$]\Rightarrow$ There is a commutative ring $S$ and a r.h $\phi:R\rightarrow S$ such that $[S$ is slender$]\,(3)$ and $[ \phi$ is f.p$]\,(4)$. It follows from $\{\ref{fpcardinal},\,(1),\, (4)\}$ that   $[\,|S|\in \mathcal{L}_1'\,]\,(5)$ and from  $\{ \ref{measurablefacts}\mathrm{.i},\,(5)  \}$ that $[\, (S^{\lambda})^{\ast} \in \mathrm{Free}(R)\,]\,(6)$. Also by assumption 
 $[\, R^{[\lambda]} \overset{_s}{\underset{R}\rightarrowtail} R^\kappa\,]\overset{\ref{def.kdim}}\Rightarrow$
$[\, R^{[\kappa]} \oplus L \cong_R R^{\lambda}$, for some $L\in R$-$\mathrm{Mod}\,]\underset{(4)}{\overset{\ref{rh}}\Longrightarrow}$
$[\, S^{[\kappa]} \oplus L \cong_S S^{\lambda}\,]
\Rightarrow$  $[\, (S^{[\kappa]})^{\ast} \oplus L^{\ast} \cong_S (S^{\lambda})^{\ast}\,]\overset{(6)}\Rightarrow$
 $[\, S^{\kappa} \oplus L^{\ast} \in \mathrm{Free}(S)\,]\Longrightarrow$
$[\, S^{\kappa}  \in \mathrm{Proj}(S)\,]\underset{\ref{oneil3}} \Longrightarrow[\, R$ is Artinian$]\Longrightarrow [R$ is not w-slender], which contradicts the assumption. \\
\underline{$\mathrm{Proof}$ of $\mathrm{(ii)}$}:  By assumption $[\, R$ is slender$\,]\,(8)$,  $[\, \kappa \in \mathcal{L}_1'\,]\,(9)$
and  $[\,  \mu_R(M)\geq \omega\,]\,(10)$. It follows from $\{\ref{measurablefacts}\mathrm{.ii},\, (8),\,(9)   \}$ that $[\, R^{\kappa}\in \mathrm{Relf}(R)\,]\,(11)$  and  $[\, (R^{\kappa})^{\ast}\cong_R R^{[\kappa]}\,]\,(12)$.  Also, by assumption we have   $[\, M \overset{_s}{\underset{R}\rightarrowtail} R^\kappa\,]\overset{\ref{def.kdim}}\Longrightarrow$
$[\,M \oplus L\cong_R R^{\kappa}$, for some  $L\in R$-$\mathrm{Mod}\,]\,(13)\Rightarrow$
$[\,M^{\ast} \oplus L^{\ast}\cong_R (R^{\kappa})^{\ast} \overset{(12)} {\cong}  _R R^{[\kappa]}\,]\Rightarrow$
$[\,M^{\ast} \oplus L^{\ast}\cong_R  R^{[\kappa]}\,]\,(14)\Rightarrow$
$[\, M^{\ast}\in \mathrm{Proj}(R)\,]\,(15)$. It follows from $\{\ref{strangesum3}\mathrm{.ii},\, (11),\,(13)  \}$ that
$[\, M\in \mathrm{Refl}(R)\,]\,(16)\Rightarrow[\, (M^{\ast})^{\ast}\cong_R M\,]\,(17)$. Assume that $[\,  \mu_R(M^{\ast})< \omega\,]\,(18)$. It follows from $\{\ref{projective}\mathrm{.v.a},\,  (17),\, (18)  \}$  that $[\,  \mu_R(M) <\omega\,]$, which contradicts $(10)$. Hence $[\,  \mu_R(M^{\ast}) \geq \omega\,]\,(19)$. Also,we have  by assumption  that 
$[R$ is connected$]\,(20)$.  It follows from $\{\ref{Bass},\,
 (19),\,(20)$ that $[\, M^{\ast}\in \mathrm{Free}(R)\,]\Rightarrow$
$[\, M^{\ast}\cong_R (R^{[\mu]}  \,]\,(21)$, for some $\mu \in \mathrm{Card}\,]$. 
Note that we could have derived $(21)$ by adding the assumption that $R$ is a $\mathrm{QS}$-ring without assuming $(10)$. Finally it follows from $\{(17),\, (21) \}$ that $[\, M \cong_R R^{\mu}\,]$.   \end{proof}

\begin{remark} \label{remarkondirectsummands} Let $R$ be a Noetherian commutative ring that is not Artinian. Then$:$ \\
$\mathrm{(i)}$ Assume that  $\mathcal{L}_1'=\emptyset$.  It follows from $\ref{hilbert}$ that  Theorem $\ref{direct.cofree}\mathrm{.i}$ is applicable when $R$    is a Hilbert ring. So any free direct summand of a dual of a free $R$-module is finitely generated. \\
$\mathrm{(ii)}$ Assume that  $\mathcal{L}_1'=\emptyset$. It follows from $\ref{hilbert}$ that Theorem $\ref{direct.cofree}\mathrm{.ii}$  is applicable when $R$    is a Hilbert domain. So any non-f.g direct summand of a cofree  $R$-module is cofree.   \\
$\mathrm{(iii)}$ The assumption   $ \mu_R(M)\geq \omega $  on  
 $\ref{direct.cofree}\mathrm{.ii}$ is redundant when $R$ is a $\mathrm{QS}$-ring (e.g $R=\mathbb{Z}$). The special case of Theorem $\ref{direct.cofree}\mathrm{.ii}$, where $R=\mathbb{Z}$  is proven in $\cite{eklof}$ in Theorem 1.2 of page 294. The special case of Theorem $\ref{direct.cofree}\mathrm{.i}$, where $R=\mathbb{Z}$,  can be derived by simply applying
 \it{Specker's Theorem} and  \it{Baer's Theorem}.\\
$\mathrm{(iv)}$ The  connectedness assumption about $R$  on  Theorem 
$\ref{direct.cofree}\mathrm{.ii}$   is not  redundant. Let $R_1,\, R_2$ be slender Noetherian commutative rings and $R=R_1 \times R_2$. Let     $ \kappa=2^{\, \mathfrak{c}}$, 
$\lambda_1=\omega$, $\lambda_2=\mathfrak{c}$ and $[\,M=\overset{2}{\underset{i=1} \boxplus}  R_i^{\lambda_i}\,]\,(3)$. It follows from $\ref{strangesum}$ that $[\, M\overset{_s} {\underset{R}\rightarrowtail} R^{\lambda}\,]\,(4)$, since $[\, R_i^{\lambda_i}\overset{_s} {\underset{R}\rightarrowtail} R_i^{\lambda}\,]\,\forall\, i\in \{ 1,2 \}$. So if $[\, M\cong_R R^{\mu}$, for $\mu \in \mathrm{Card}\,]$   it would follow from $\ref{strangesum}$ that  $[\, R_i^{\lambda_i}\cong_{R_i} R_i^{\mu}\,]\,\forall\, i\in \{ 1,2 \}\Longrightarrow[\,\omega= \lambda_1=\mu,\,\mathfrak{c} =\lambda_2=\mu\,]\Rightarrow[\, \omega=\mathfrak{c}\,]$, which is a contradiction. Therefore $M$  is not a cofree $R$-module although $M$ is an $R$-direct summand of $R^{\kappa}$. \\
$\mathrm{(v)}$  Theorem $\ref{direct.cofree}\mathrm{.i}$ extends
 Lemma  $\ref{invertible}\mathrm{.i}$. \\$\mathrm{(vi)}$  Theorem $\ref{steinitz7}$ was essential for the proof of $\ref{invertible}$ since it waas used in the proof of  Theorem $\ref{steinitz7}$. 
 \end{remark}

I am grateful to my academic adviser at the  \emph{University of Kentucky}, Professor Edgar Enochs, for reading my earlier manuscript and making suggestions. For the same reason, i am grateful to Professor David Leep (\emph{University of Kentucky}) and to Professor 
Ioannis Emmanouil (\emph{National and Kapodistrian University of Athens}).  I am also grateful to Professor Athanasios Tsarplalias, for teaching me cardinal arithmetic at the \emph{National and Kapodistrian University of Athens}

\end{document}